\documentclass[11pt, oneside]{mnthesis}
\usepackage{epsfig,epic,eepic,subfigure}
\usepackage{hyperref}
\usepackage{url}
\usepackage{longtable}
\usepackage{multirow}
\usepackage{bigstrut}
\usepackage{tikz}
\usepackage{graphicx}
\usepackage{xypic}
\usepackage{amsmath,amssymb,amsthm,bm,mathtools}
\usepackage{mathrsfs}

\newcommand{\ud}{\, \mathrm{d}}

\newtheorem{theorem}{Theorem}[chapter]
\newtheorem{lemma}[theorem]{Lemma}
\newtheorem{proposition}[theorem]{Proposition}
\newtheorem{definition}[theorem]{Definition}
\newtheorem{corollary}[theorem]{Corollary}
\newtheorem{remark}[theorem]{Remark}

\def\diam{\mathop{\rm diam}}

\def\lim{\mathop{\rm lim}}

\def\pf{{\it Proof:}~}
\def\var{\mathop{\rm Var}}

\providecommand{\abs}[1]{\lvert#1\rvert}
\providecommand{\Abs}[1]{\Bigr\lvert#1\Bigl\rvert}
\providecommand{\norm}[1]{\lVert#1\rVert}
\providecommand{\Bnorm}[1]{\Bigl\lVert#1\Bigr\rVert}

\begin{document}
\bibliographystyle{hunsrt} 

\phd 

\title{\bf Diffusivity and Ballistic Behavior of Random Walk in Random
Environment}
\author{Xiaoqin Guo}
\campus{University of Minnesota} 
\program{Mathematics} 
\director{Adviser: Ofer Zeitouni} 

\submissionmonth{September} 
\submissionyear{2012} 

\copyrightpage 
\acknowledgements{

First of all, I would like to thank my advisor, Ofer Zeitouni, for his advice
and tremendous help during my graduate study. Whenever I asked a
question, Professor Zeitouni answered quickly with insightful
responses, no matter how busy he was. 
His devotion to research, deep understanding of mathematics, constant
encouragement and warm personality are
among the treasures in my memory.
I owe him the deepest gratitude from my heart.

I owe respect and gratitude to Uncles Chen Guang and Peng Zhai who helped me
finish my high school education when my family was in financial difficulty.

I would like to thank Professors John Baxter, Maury Bramson, and Nicolai Krylov
for being on my thesis defense committee and Professor Gilad Lerman whose
suggestions improved my teaching.

I would also like to thank my friends in college and graduate school who helped
make my experience more enjoyable.

And finally, I wish to express my sincerest gratitude to my wife, for her
support and understanding while I was pursuing my academic studies. 

}
\dedication{\begin{center}{\Large To my grandpa and my parents}\end{center}}


\beforepreface 


\afterpreface         


\chapter{Introduction}
\label{intro_chapter}
\section{An introduction to RWRE}\label{Section II}

Let $\mathcal{M}=\mathcal{M}_1(V)$ be the space of all probability measures on 
$V=\{v\in\mathbb{Z}^d: |v|\le 1\}$,
 where $|\cdot|$ denotes the $l^2$-norm. We equip  $\mathcal{M}$ with the weak topology on probability measures, which makes it into a Polish space, and equip $\Omega=\mathcal{M}^{\mathbb{Z}^d}$ with the induced Polish structure. Let $\mathcal{F}$ be the Borel
$\sigma$-field of $\Omega$ and $P$ a probability measure on $\mathcal{F}$.

A random \textit{environment} is  an element $\omega
=\{\omega(x, v)\}_{x\in{\mathbb{Z}^d}, v\in V}$ of $\Omega$. The random environment is called \emph{balanced} if \[P\{\omega(x, e_i)=\omega(x,-e_i) \mbox{ for all $i$ and all $x\in\mathbb{Z}^d$}\}=1,\]
and \textit{elliptic} if $P\{\omega(x,e)>0 \mbox{ for all $|e|=1$ and all $x\in\mathbb{Z}^d$}\}=1$. We say that the random environment is \emph{uniformly elliptic} with ellipticity constant $\kappa$ if $P\{\omega(x,e)>\kappa \mbox{ for all $|e|=1$ and all $x\in\mathbb{Z}^d$}\}=1$.

The random walk in the
random environment $\omega\in\Omega$ (RWRE)
started at $x$ is the Markov
chain $\{X_n\}$ on $(\mathbb{Z}^d)^\mathbb{N}$,
with state space $\mathbb{Z}^d$ and law $P_\omega^x$ specified by
\begin{align*}
&P_\omega^x\{X_0=x\}=1,\\
&P_\omega^x\{X_{n+1}=y+v | X_n=y\}=\omega(y, v), \quad v\in V.
\end{align*}
Let $\mathcal{G}$ be the $\sigma$-field generated by cylinder functions. The probability distribution $P_\omega^x$ on $((\mathbb{Z}^d)^\mathbb{N}, \mathcal{G})$ is called the \textit{quenched law}.
Note that for each $G\in\mathcal{G}$, $P_\omega^x (G) : \Omega\to [0,1]$ is a $\mathcal{F}$-measurable function.
The joint probability distribution $\mathbb{P}^x$ on $\mathcal{F}\times\mathcal{G}$:
\[
\mathbb{P}^x (F\times G)=\int_{F} P_\omega^x (G)P(\ud\omega), \qquad F\in\mathcal{F},\, G\in\mathcal{G},
\]
is called the \textit{annealed} (or \textit{averaged})
law. Expectations with respect to $P_\omega^x$ and $\mathbb{P}^x$ are denoted
by $E_\omega^x$ and $\mathbb{E}^x$, respectively. We also write $\mathbb{P}^o$ as $\mathbb{P}$, where
$o=(0,\cdots, 0)$ is the origin.


For $\omega\in\Omega$, set 
\[\omega_x=\big(\omega(x,e)\big)_{|e|=1}.\]
Define the spatial shifts $\{\theta^y\}_{y\in\mathbb{Z}^d}$ on $\Omega$ by
$(\theta^y\omega)_x=\omega_{x+y}$. We say that the random environment is \textit{ergodic}
if the measure $P$ is ergodic with respect to the group of shifts $\{\theta^y\}$.
A special case is when the probability vectors $(\omega_x)_{x\in\mathbb{Z}^d}$
are independent and identically distributed (\textit{iid}).

Setting $\bar{\omega}(n)=\theta^{\mathrm{X}_n}\omega$,
then the process $ \bar{\omega}(n) $ is a Markov chain under
$ \mathbb{P}^{o} $ with state space $ \Omega $ and transition kernel \[ M(\omega',\ud\omega)=\sum_{i=1}^{d}[\omega'(o,e_i)\delta_{\theta^{e_i}\omega'}+
\omega'(o,-e_i)\delta_{\theta^{-e_i}\omega'}]+\omega'(o,o)\delta_{\omega'}. \]
$\big(\bar{\omega}(n)\big)_{n\in\mathbb{N}}$ is often referred as the ``\textit{environment viewed from the point of view of the particle}'' process.

For $t\ge 0$, let
\[
X_t=X_{\lfloor t\rfloor}+(t-\lfloor t\rfloor)(X_{\lfloor t\rfloor+1}-X_{\lfloor t\rfloor}).
\]
 We say that the \textit{quenched invariance principle} of the RWRE holds if,
for $P$-almost every $\omega\in\Omega$ and some deterministic vector
$v\in\mathbb{R}^d$ (called the {\it limiting velocity}),
the $P_\omega^o$ law of the path $\{(X_{tn}-tnv)/\sqrt{n}\}_{t\geq 0}$ converges
weakly to a Brownian motion,
as $n\to \infty$. For $\ell\in S^{d-1}$, we say that the RWRE is
\emph{ballistic} in the direction $\ell$ if
\[
 \varliminf_{n\to\infty}\frac{X_n\cdot\ell}{n}>0,\quad\mathbb{P}\mbox{-a.s.}
\]

\section{Structure of the thesis}
In this thesis, we study the diffusive and ballistic behaviors of random walks in random environment in $\mathbb{Z}^d, d\ge 2$. 

The organization of the thesis is as follows. 

Section~\ref{IOverview} gives an overview of the previous results in the study of the ballisticity, the central limit theorems (CLT), and the Einstein relation of RWRE. The three subsections in Section~\ref{Iresults} state the main results in this thesis and discuss the ideas of their proofs. 

Chapters \ref{LV chapter}, \ref{CLT chapter} and \ref{ER chapter} are devoted to the proofs of our three main results:

In Chapter~\ref{LV chapter}, we consider the limiting velocity of random walks in strong-mixing random Gibbsian environments in $\mathbb{Z}^d, d\ge 2$.
Based on regeneration arguments, we will first provide an alternative proof of Rassoul-Agha's conditional law of large numbers (CLLN) for mixing environment \cite{R-A3}.
Then, using coupling techniques, we show that there is at most one nonzero limiting velocity in high dimensions ($d\ge 5$).

Chapter~\ref{CLT chapter} proves the quenched invariance principles (Theorem~\ref{CLT1} and Theorem~\ref{CLT2}) for random walks in elliptic and balanced environments.
We first prove an invariance principle (for $d\ge 2$) and the
transience of the random walks when $d\ge 3$ (recurrence when $d=2$)
in an ergodic environment which is not uniformly elliptic but satisfies
certain moment condition. Then, using percolation arguments, we
show that under (not necessarily uniform) ellipticity, the above results hold
for random
walks in iid balanced environments.

Chapter~\ref{ER chapter} gives the proof of the Einstein relation in the context
of random walks in a balanced uniformly elliptic iid random environment. Our
approach combines a change of measure argument of Lebowitz and Rost \cite{Le}
and the regeneration argument of Gantert, Mathieu and Piatnitski \cite{GMP}. The
key step of our proof is the construction of a new regeneration structure.




\section{Overview of previous results}\label{IOverview}
\subsection{Ballisticity}
The ballistic behavior of the RWRE in dimension $d\ge 2$ has been extensively
studied.
For random walks in iid random environment in dimension $d\ge 2$,
the Kalikow's 0-1 law \cite{Ka81} states that for any direction $\ell\in S^{d-1}$, 
\[\mathbb{P}(A_\ell\cup A_{-\ell})\in\{0,1\}\]
where 
$A_{\pm\ell}=\{\lim_{n\to\infty}X_n\cdot\ell=\pm\infty\}$.
It is believed that for any direction $\ell$ and any $d\ge 2$,
a stronger 0-1 law is true:
\[
P(A_\ell)\in\{0,1\} \tag{0-1 Law}.
\]
When $d=2$, this $0$-$1$ law was proved by Zerner and Merkel \cite{ZM}. The question whether the 0-1 law holds for iid random environment in dimensions $d\ge 3$ is still open. (It is known that some strong mixing condition is necessary for the 0-1 law to hold, as the counterexample in \cite{BZZ} shows.)

Much progress has been made in the study of the limiting velocity $\lim_{n\to\infty}X_n/n$ of random walks in iid environment, see \cite{ZO} for a survey.  
For one-dimensional RWRE, the law of large numbers (LLN) was proved in \cite{So}.
For $d\ge 2$, a conditional law of large numbers (CLLN) was proved in \cite{SZ, Ze} (see \cite[Theorem 3.2.2]{ZO} for the full version). It states that $\mathbb{P}$-almost surely, for any direction $\ell$,
\[
\lim_{n\to\infty}\frac{X_n\cdot\ell}{n}=v_+ 1_{A_\ell}-v_-1_{A_{-\ell}}
\tag{CLLN}
\]
for some deterministic vectors $v_\ell$ and $v_{-\ell}$ (we set $v_\ell=o$ if
$\mathbb{P}(A_\ell)=0$). 
This was achieved by considering the regenerations of the random walk path.
Hence for $d\ge 2$, the 0-1 law would imply the LLN. Recall that when $d\ge 3$,
the 0-1 law is one of the main open questions in the study of RWRE. 
Nevertheless, in high dimensions ($d\ge 5$), Berger \cite{Be} showed
that the limiting velocity can take at most one non-zero value, i.e.,
\begin{equation}\label{Berger}v_\ell v_{-\ell}=0.\end{equation}

It is of interests to consider environments whose law $P$ is not iid but rather
ergodic (under possibly appropriate mixing conditions). Of special interest is
the environment that is produced by a Gibbsian particle system (which we call
the \textit{Gibbsian environment}) and satisfies Dobrushin-Shlosman's
strong-mixing condition IIIc in \cite[page 378]{DS}, see
\cite{R-A1,R-A2,CZ1,CZ2, R-A3} for related works.  
 An important feature of this model is that the influence of the environments in
remote locations decays exponentially as the distance grows. (We won't give the
definitions of the Gibbsian environment and the strong-mixing
condition in this thesis. For their definitions, we refer to \cite[pages
1454-1455]{R-A1}. We remark that our results only assume a mixing condition
(G), which is defined in page \pageref{LVdef1}. It is known that (G) is a
property of the strong-mixing Gibbsian environment, cf. \cite[Lemma 9]{R-A1}.)

In \cite{R-A1}, assuming a ballisticity condition (Kalikow's condition) which implies that the 
event of escape in a direction has probability $1$,  Rassoul-Agha proved the LLN for the strong-mixing Gibbsian environment, using the invariant measure of the ``environment viewed from the point of view of the particle" process $\big(\bar\omega(n)\big)$. 
In \cite{R-A3}, Rassoul-Agha also obtained the CLLN for the strong-mixing Gibbsian environment, under an analyticity condition (see Hypothesis (M) in \cite{R-A3}).
Comets and Zeitouni proved the LLN for environments with a weaker cone-mixing assumption ($\mathcal{A}1$) in \cite{CZ1}, but under some conditions about ballisticity and the uniform integrability of the regeneration times (see ($\mathcal{A}5$) in \cite{CZ1}).

\subsection{Central Limit Theorems}
In recent years, there has been much interest in the study
of invariance principles and transience/recurrence
for random walks in random environments (on the
$d$-dimensional lattice $\mathbb{Z}^d$)
with non uniformly
elliptic transitions probabilities. Much of this work has 
been in the context of reversible models, either for walks on percolation
clusters or for the random conductance model, see
\cite{Bar04,SS05,MR05,BB,MaP07,Ma08,BarDe10}.  In those cases,
the main issue is the transfer of annealed estimates (given e.g.
in \cite{DFGW89} in great generality) to the quenched setting, and the control
of the quenched mean displacement of the walk. 
On the other hand, in these models the reversibility of the walk
provides for explicit expressions for certain invariant measures
for the environment viewed from the point of view of the particle.

The non-reversible setup has proved to provide many additional, and
at this point insurmountable, 
challenges, even in the uniformly elliptic iid setup, see 
\cite{Zrev} for a recent account. 
In \cite{Sz3}, Sznitman shows that his condition (T') implies ballisticity and
LLN and a directional annealed central limit theorem. The proof uses
regeneration times and a renormalization argument and does not employ the
process of the environment viewed from the point of view of particle. (We
remark that weaker forms of the condition (T') exist, see \cite{Sz3, DR1, DR2,
BDR}. Recently it was shown in \cite{BDR} that polynomial decay of some exit
probabilities implies (T').)
Further, it was shown by Berger and Zeitouni \cite{BZei} and Rassoul-Agha and Sepp\"{a}l\"{a}inen \cite{R-AS} that in the ballistic case, an annealed invariance principle is equivalent to a quenched invariance principle, under appropriate moment conditions on the regeneration times (these conditions are satisfied in all cases where a ballistic annealed CLT has been proved).

When the walk is not ballistic, the regeneration structure employed in \cite{Sz2} is not available. Several classes of non-ballistic models were considered in the literature: balanced environment (see the definition in Section~\ref{Section II}),
environment whose sufficiently high-dimensional projection is a simple random walk \cite{BSZ}, and isotropic environment which is a small perturbation of the simple random walk \cite{BK, BolZei, SZei}. Historically, the first to be considered was the balanced environment, first investigated by Lawler \cite{La}, which we describe next as a good part of the thesis deals with that environment:
\begin{theorem}[\cite{La},\cite{ZO}]\label{LaThm}
	Assume the random environment is ergodic, balanced and uniformly elliptic. Then $P$-almost surely, the $P_\omega$ law of the rescaled path $\lambda X_{\cdot/\lambda^2}$ converges weakly to a Brownian motion on $\mathbb{R}^d$ with a non-degenerate diagonal covariance matrix. Moreover, the RWRE is recurrent for $d=2$ and transient for $d\ge 3$, $P$-almost surely.
\end{theorem}
In this case, a-priori estimates of the Alexandrov-Bakelman-Pucci
type give enough control that allows one to prove the existence
of invariant measures (for the environment viewed from
the point of view of the particle), and the fact that the walk
is a (quenched) martingale together with ergodic arguments yield
the invariance principle (obviously, control of the quenched 
mean displacement,
which vanishes, is automatic). The establishment of recurrence (for $d=2$)
and transience (for $d\geq 3$) requires some additional
arguments, due to Kesten and Lawler, respectively, see
\cite{ZO} for details. 
\subsection{Einstein relation}
In 1905, Einstein \cite[pp. 1-18]{Einstein} investigated the movement of suspended particles in a liquid under the influence of an 
external force. He established the following linear relation between the diffusion constant $D$ and the
\textit{mobility} $\mu$:
\[
D\sim T\mu,
\]
where $T$ is the absolute temperature, and $\mu$ is defined as the limiting
ratio between the velocity (under the external force) and the force, as the
force goes to zero.

More precisely, the Einstein relation (ER) describes the relation between
the response of a system to a perturbation and its diffusivity at equilibrium. It states that the derivative of the velocity (with respect to the strength of the perturbation) equals the diffusivity:
\[
    \lim_{\lambda\to 0}\lim_{t\to\infty}\frac{E_\lambda X_t/t}{\lambda}=D, \tag{ER}
\]
where $(X_t)_{t\ge 0}\in (\mathbb{R}^d)^{\mathbb{R}_+}$ denotes the random
motion of the particle, $\lambda$ is the size of the perturbation, $D$ is the
diffusion constant of the equilibrium state, and $E_\lambda$ is the annealed
measure of the perturbed media. 
 General derivations of this principle assume reversibility.

Recently, there has been much interest in studying the Einstein relation for reversible motions
in random media, see \cite{Le,KO,GMP,BHOZ}. 
In \cite{Le}, Lebowitz and Rost proved a weak form of the Einstein relation for a wide class of random motions in random media:
\[
	\lim_{\lambda\to 0}E_\lambda \frac{X_{t/\lambda^2}}{t/\lambda}=D \quad \forall t>0.
\]
In \cite{KO}, the ER is verified for random walks in random conductance, where the conductance is only allowed to take two values. The approach of \cite{KO} is an adaption of the perturbation argument and transience estimates in \cite{Loulakis}. For random walks on Galton-Watson trees, the ER is proved by \cite{BHOZ}. Their approach uses recursions due to the tree structure and renewal arguments. Recently, Gantert, Mathieu and Piatnitski \cite{GMP} established the ER for random walks in random potential, by combining the argument in \cite{Le} with good moment estimates of the regeneration times.

The Einstein relation for random motions in the non-reversible zero speed
set-up, e.g., random walks in balanced random environments (RWBRE), is a
challenging problem. (In general one expects correction terms in (ER) due to the
non-reversibility of the walk.)

\section{Our results}\label{Iresults}
In this section we will state the main results in the thesis and explain the
ideas of their proofs. The actual proofs will be presented in the following
chapters.

 Our contributions are in three directions: CLLN and regeneration structures for
RWRE in Gibbsian environments, quenched invariance principles for balanced
elliptic (but non uniformly elliptic) environments, and ER for balanced iid
uniformly elliptic environments.

\subsection{Limiting velocity for mixing random environment}\label{ILV}
Recall first the definition of an $r$-Markov environment (see \cite{CZ2}). 
\begin{definition}\label{LVdef1}
For $r\ge 1$, let $\partial_r V=\{x\in\mathbb{Z}^d\setminus V: d(x, V)\le r\}$ be the $r$-boundary of $V\subset\mathbb{Z}^d$.
A random environment $(P,\Omega)$ on $\mathbb{Z}^d$ is called $r$-Markov if 
for any finite $V\subset\mathbb{Z}^d$,
\[
P\big((\omega_x)_{x\in V}\in \cdot|\mathcal{F}_{V^c}\big)
=P\big((\omega_x)_{x\in V}\in \cdot|\mathcal{F}_{\partial_r V}\big), \text{ $P$-a.s.,}
\]
where $d(\cdot,\cdot)$ denotes the $l^1$-distance and $\mathcal{F}_{\Lambda}:=\sigma(\omega_x:x\in\Lambda)$.
\end{definition}
We say that an $r$-Markov environment $P$ {\it satisfies condition (G)} if there
exist 
constants $\gamma , C<\infty$ such that for all finite subsets $\Delta\subset V\subset\mathbb{Z}^d$ with $d(\Delta,V^c)\ge r$, and $A\subset V^c$,
\[
\frac{\ud P\big((\omega_x)_{x\in\Delta}\in\cdot|\eta\big)}
{ \ud P\big((\omega_x)_{x\in\Delta}\in\cdot|\eta'\big)}
\le 
\exp{(C\sum_{x\in A,y\in\Delta}e^{-\gamma d(x,y)})}\tag{G}
\]
for $P$-almost all pairs of configurations $\eta,\eta'\in\mathcal{M}^{V^c}$ which agree on $V^c\setminus A$.
Here 
\[
P\big((\omega_x)_{x\in\Delta}\in\cdot|\eta\big)
:=P\big((\omega_x)_{x\in\Delta}\in \cdot|\mathcal{F}_{V^c}\big)\big|_{(\omega_x)_{x\in V^c}=\eta}.
\]
We remark that $r$ and $\gamma$ are used as parameters of the environment throughout the article.

Recall that by Lemma 9 in \cite{R-A1}, the strong-mixing Gibbsian environment
satisfies (G).
Obviously, every finite-range dependent environment also satisfies (G).

Our main theorem concerning the mixing environments is:
\begin{theorem}\label{LVthm2}
Assume that $P$ is uniformly elliptic and satisfies \emph{(G)}.
Then there exist two deterministic constants $v_+, v_-\ge 0$
and a vector $\ell$ such that
\begin{equation}\label{ICLLN}
\lim_{n\to\infty} \frac{X_n}{n}=v_+\ell 1_{A_\ell}-v_-\ell 1_{A_{-\ell}},
\end{equation}
and $v_+=v_-=0$ if $\mathbb{P}(A_\ell\cup A_{-\ell})<1$.
Moreover, if $d\ge 5$, then
there is at most one non-zero velocity. That is, 
\begin{equation}\label{Iunique}
v_+ v_-=0.
\end{equation}
\end{theorem}
We remark here that for the finite-range dependent
case, the CLLN is proved in \cite{ZO}.

\eqref{ICLLN} is a minor extension of Rassoul-Agha's CLLN in \cite{R-A3}. He
assumes slightly more than strong-mixing, which in turn is slightly stronger
than our condition (G). Our proof is very different from the proof in 
\cite{R-A3} , which is based on a large deviation principle in \cite{R-A2}. The
main contribution of our proof of \eqref{ICLLN} is a new definition of the
regeneration structure, which enables us to divide a random path in the mixing
environment into ``almost iid" parts.
With this regeneration structure, we will use the ``$\epsilon$-coins" introduced
in \cite{CZ1} and coupling arguments to prove the CLLN. This regeneration
structure will also be used in the proof of \eqref{Iunique}.

\begin{figure}[h]
\centering
\includegraphics[width=0.6\textwidth]{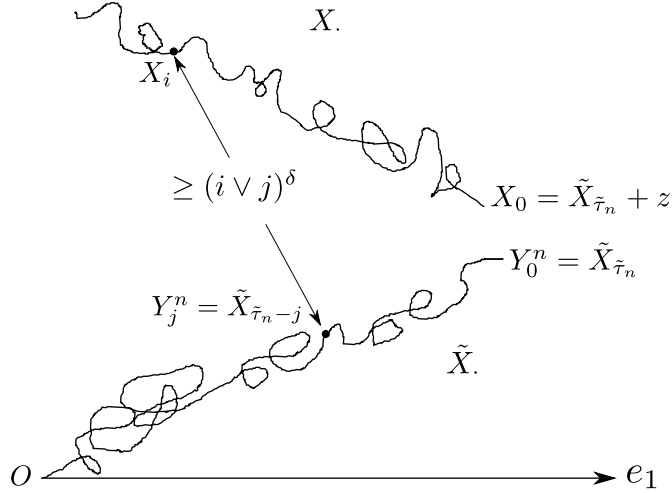}
\caption{$X_\cdot\in A(Y_\cdot^n, z)$. 
When $i\vee j>|z|/3$, the distance between $Y_j^n$ of the ``backward path" and $X_i$ is at least $(i\vee j)^\delta$.}
\label{LVfig:1}
\end{figure}

Display \eqref{Iunique} is an extension of Berger's result \eqref{Berger} from
the iid case to our case (G), which includes the strong-mixing case. 
In \cite{Be}, assuming that $\mathbb{P}(A_\ell)>0$ for a direction $\ell$, Berger coupled the iid environment $\omega$ with a transient (in the direction $\ell$) environment $\tilde\omega$ and a ``backward path", such that $\tilde\omega$ and $\omega$ coincide in the locations off the path. 
Using heat kernel estimates for random walks with iid increments, he showed that if $v_\ell v_{-\ell}>0$ and $d\ge 5$, then with positive probability, the random walks in $\tilde\omega$ is transient to the $-\ell$ direction without intersecting the backward path, which contradicts $\tilde{\omega}$ being transient in the direction $\ell$.
The difficulties in applying this argument to mixing environments are that the regeneration slabs are not
iid, and that unlike the iid case, the environments visited by two disjoint paths are not independent. To overcome these difficulties,
we will construct an environment (along with a path) that is ``very transient" in $\ell$,  and show that the ballistic walks in the opposite direction $-\ell$ will move further and further away from the given path (see Figure \ref{LVfig:1} in Section \ref{secunique}). The key ingredient here is a heat kernel estimate, which we will obtain in Section \ref{sechke} using coupling arguments. 
\subsection{Invariance principle for RWBRE}
As mentioned above, Lawler \cite{La} proved the invariance principle under the uniform ellipticity assumption.
We explore the extent to which the uniform ellipticity assumption can be dropped. 
Surprisingly, in the iid case, we can show that no assumptions of uniform ellipticity are needed at all.

Let
\begin{equation}
	\label{CLTepsdef}
	\varepsilon(x)=\varepsilon_{\omega}(x):=
	[\prod_{i=1}^{d}\omega(x,e_i)]^{\frac{1}{d}}.
\end{equation}
Our first main result
is that
if $\mathrm{E}\varepsilon(o)^{-p}< \infty$ for some $p>d$, then
the quenched invariance principle
holds and moreover, the RWRE is transient $P$-almost surely if
$d\geq 3$. (Recurrence for $d=2$ under 
the condition $E\varepsilon(0)^{-p}<\infty$
follows from the quenched invariance principle and ergodicity by an 
unpublished argument of Kesten 
detailed in
\cite[Page 281]{ZO}. Note that this argument cannot be used to
prove transience in dimensions $d\geq 3$, 
even given an invariance principle, since in higher dimensions
the invariance principle does not give useful
information on the range of the
random walk; the behavior of the
range is a crucial element in Kesten's argument.)

\begin{theorem}\label{CLT1}
Assume that the random environment is ergodic, elliptic
 and balanced.
\begin{enumerate}
\item[(i)] If $E\varepsilon(o)^{-p}< \infty$ for some $p>d\ge 2$,
then
 the quenched invariance principle holds with a
non-degenerate diagonal limiting covariance matrix.
\item[(ii)] If $E[(1-\omega(o,o))/\varepsilon(o)]^q< \infty$ for some $q>2$ and $d\ge 3$, then the RWRE is transient $P$-almost surely.
\end{enumerate}
\end{theorem}
\noindent
That some integrability condition on the tail of $\varepsilon(o)$ is needed
for part (i) to hold
is made clear by the (non-Gaussian) scaling limits of random walks in
Bouchaud's trap model, see \cite{Bou,BAC}.
In fact, it follows from that example that Theorem \ref{CLT1}(i), or even an annealed
version of the CLT, cannot hold in general with  $p<1$.

 The proof of Theorem \ref{CLT1} is based on a sharpening of the arguments in
\cite{La,Sz1,ZO}; in particular, refined versions of the maximum
principle for walks in balanced environments (Theorem \ref{CMP})
and of a mean value inequality (Theorem \ref{Cmvi}) play a crucial role.

When the environment is iid and elliptic, our second main result
is that if $|X_{n+1}-X_n|=1$ a.s., then the quenched invariance principle holds.
Moreover, the RWRE is $P$-almost surely transient when $d\ge 3$.
The proofs combine percolation arguments with Theorem \ref{CLT1}.
\begin{theorem}\label{CLT2}
Assume that the random environment is iid, elliptic and
balanced.
\begin{enumerate}
\item[(i)] If $P\{\max_{|e|=1}\omega(o,e)\ge \xi_0\}$=1 for some positive constant $\xi_0$, then the 
quenched invariance principle holds with a non-degenerate limiting covariance.
\item[(ii)] When $d\ge 3$, the RWRE is transient $P$-almost surely.
\end{enumerate}
\end{theorem}
Because the 
transience or recurrence of the random walks does not change
if one considers the walk restricted
to its jump times,
one concludes, using Kesten's argument and the invariance
principle, comparing with
Theorem \ref{CLT1}, that
for $d=2$, a random walk in a balanced elliptic iid
random environment is recurrent $P$-a.s.

Our proof of the invariance principles, like that of \cite{La}, is based
on the approach of the
``environment viewed from the point of view of the particle".

Since $\{X_n\}$ is a (quenched) martingale,
standard arguments (see the proof of Theorem 6.2 in \cite{BB}) show that
the quenched invariance principle holds
whenever  an invariant
measure $Q\sim P$ of $\{\bar{\omega}(n)\}$ exists.
 The approach of Lawler \cite{La}, which is a discrete version of the argument of Papanicolaou and Varadhan \cite{PV}, is to construct such a measure as the limit of invariant measures of periodized environments. We will
 follow this strategy using, as in \cite{Sz1,ZO}, variants of \cite{KT} to derive estimates on solutions of linear elliptic difference
 equations. In the iid setup of
Theorem~\ref{CLT2}, percolation estimates are used to control
pockets of the environment where those estimates are not strong enough.

For the proof of the transience in the ergodic case, 
we use a mean value inequality and follow \cite{ZO}.
To prove the transience in the iid case, we employ percolation
arguments together with
a new maximum principle (Theorem \ref{Cmp2}) for walks with (possibly)
big jumps.
\begin{remark}
Recently, Berger and Deuschel \cite{BD} have generalized our ideas and extended
the quenched invariance principle to the general non-elliptic case where the environment is only required to be iid and “genuinely $d$-dimensional”. 
\end{remark}

\subsection{Einstein relation for RWBRE}\label{IER}
In this subsection we will present the Einstein relation for random walks in uniformly elliptic balanced iid random environment. Recall that by Theorem~\ref{LaThm}, for $P$-almost every $\omega$, $(\lambda X_{t/\lambda^2})_{t\ge 0}$ converges weakly (as $\lambda\to 0$) to a Brownian motion with a non-degenerate covariance matrix, which we denote by $\bm{D}$.

For $\lambda\in (0,1)$ and a fixed direction
\[
\ell=(\ell_1,\ldots,\ell_d)\in S^{d-1},
\] 
define the perturbed environment $\omega^\lambda$ of $\omega\in\Omega$ by
\[
\omega^\lambda(x,e)=(1+\lambda\ell\cdot e)\omega(x,e).
\]
Since $\omega^\lambda$ satisfies Kalikow's condition (see (0.7) in \cite{SZ}),
it follows from \cite[Theorem 2.3]{SZ} that there exists a deterministic constant $v_\lambda\in\mathbb{R}^d$ such that 
\[
\lim_{t\to\infty} \frac{X_t}{t}=v_\lambda,
\quad \text{ $P\otimes P_{\omega^\lambda}^o$-almost surely}.
\]

Our main result is the
following mobility-diffusivity relation:
\begin{equation}\label{Einstein relation}
\lim_{\lambda\to 0}\frac{v_\lambda}{\lambda}=D_\ell,
\end{equation}
where
\[
D_\ell:=\bm{D}\ell=(2E_Q\omega(o,e_i)\ell_i)_{1\le i\le d}\in \mathbb{R}^d.
\] 
Our proof of the Einstein relation \eqref{Einstein relation} consists of proving the following two theorems:
\begin{theorem}\label{ER1}
Assume that the environment $P$ is iid, balanced and uniformly elliptic. Then for $P$-almost every $\omega$ and for any $t\ge 1$,
\begin{equation*}
\lim_{\lambda\to 0}
E_{\omega^\lambda}\frac{X_{t/\lambda^2}}{t/\lambda}=D_\ell.
\end{equation*}
\end{theorem}

\begin{theorem}\label{ER2}
Assume that the environment $P$ is iid, balanced and uniformly elliptic. Then for all sufficiently small $\lambda\in (0,1)$ and any $t\ge 1$,
\[
\left|E_PE_{\omega^\lambda}\frac{X_{t/\lambda^2}}{t/\lambda}-\frac{v_\lambda}{\lambda}\right|
\le
\frac{C}{t^{1/5}}.
\]
\end{theorem}

Our proof of Theorem \ref{ER1} is an adaption of the argument of Lebowitz and Rost \cite{Le} (see also \cite[Proposition 3.1]{GMP}) to the discrete setting. Namely, using a change of measure argument, we will show that the scaled process $\lambda X_{t/\lambda^2}$ converges (under the law $P_{\omega^\lambda}$) to a Brownian motion with drift $tD_\ell$, which yields Theorem \ref{ER1}.

For the proof of Theorem \ref{ER2}, we want to follow the strategy of Gantert, Mathieu and Piatnitski \cite{GMP}. Arguments in the proof of \cite[Proposition 5.1]{GMP}  show that if we can construct a sequence of random times $(\tau_n)_{n\in\mathbb{N}}$ (called the {\it regeneration times}) that
divides the random path into iid (under the annealed measure) pieces, then good moment estimates of the regeneration times yield Theorem \ref{ER2}.
In the construction of the regeneration times in \cite{GMP}, a heat kernel estimate \cite[Lemma 5.2]{GMP} for reversible diffusions is crucially employed.
However, due to the lack of reversibility, we don't have a good heat kernel estimate for RWRE. In this thesis, we construct the regeneration times differently, so that they divide the random path into ``almost iid" parts. Moreover, our regeneration times have good moment bounds, which lead to a proof of Theorem \ref{ER2}. The key ingredients in our construction are Kuo and Trudinger's \cite{KT} Harnack inequality for discrete harmonic functions and the ``$\epsilon$-coins" trick introduced by Comets and Zeitouni \cite{CZ1}.


\chapter{Limiting Velocity in Mixing Random Environment}
\label{LV chapter}
This chapter is devoted to the proof of Theorem~\ref{LVthm2}. The organization of the proof is as follows. In Section \ref{seccomb}, we prove a refined 
version of \cite[Lemma 3]{Ze}. With this combinatorial result, we will prove the CLLN \eqref{ICLLN} in Section \ref{seclln}, using coupling arguments. In Section \ref{sechke}, using coupling, we obtain
heat kernel estimates, which is later used in
Section \ref{secunique} to show the uniqueness of the non-zero limiting velocity. 

Throughout this chapter, we assume that \textit{the environment is uniformly elliptic with ellipticity constant $\kappa$ and satisfies $(G)$}. We use $c, C$ to denote finite positive constants that depend only on 
the dimension $d$ and the environment measure $P$ (and implicitly, on the parameters $\kappa,r$ and $\gamma$ of the environment). They may differ from line to line. 
We denote by $c_1,c_2,\ldots$ positive constants which are fixed throughout, and which depend only on $d$ and the measure $P$. Let $\{e_1,\ldots,e_d\}$ be the natural basis of $\mathbb{Z}^d$.

\section{A combinatorial lemma and its consequences}\label{seccomb}
In this section we consider the case that $\mathbb{P}(\varlimsup_{n\to\infty}X_n\cdot e_1/n>0)>0$. We will adapt the arguments in \cite{Ze} and prove that with positive probability,
the number of visits to the $i$-th level $\mathcal{H}_i=\mathcal{H}_i(X_0):=\{x:x\cdot e_1=X_0\cdot e_1+ i\}$ grows
slower than $Ci^2$.
An important ingredient of the proof is a refinement of a combinatorial lemma of Zerner \cite[Lemma 3]{Ze} about deterministic paths.

We say that a sequence $\{x_i\}_{i=0}^{k-1}\in (\mathbb{Z}^d)^{k}$, $2\le k\le\infty$, is a \textit{path} if 
$|x_i-x_{i-1}|=1$ for $i=1,\cdots, k-1$. For $i\ge 0$ and an infinite path $X_\cdot=\{X_n\}_{n=0}^\infty$ such that $\sup_n X_n\cdot e_1=\infty$, let
\[T_i=\inf\{n\ge 0: X_n\in\mathcal{H}_i\}.\]
For $0\le i<j$ and $k\ge 1$, let $T_{i,j}^1:=T_i$ and define
recursively
\[
T_{i,j}^{k+1}=\inf\{n\ge T_{i,j}^k: X_n\in\mathcal{H}_i \text{ and } n<T_j\}\in \mathbb{N}\cup  \{\infty\}.
\]
That is, $T_{i,j}^k$ is the time of the $k$-th visit to $\mathcal{H}_i$ before hitting
$\mathcal{H}_j$. Let
\[
N_{i,j}=\sup\{k: T_{i,j}^k<\infty\}
\]
be the total number of visits to $\mathcal{H}_i$ before hitting
$\mathcal{H}_j$.

As in \cite{Ze}, for $i\ge 0, l\ge 1$, let
\[
h_{i,l}=T_{i,i+l}^{N_{i,i+l}}-T_i
\]
denote the time spent between the first and the last visits to $\mathcal{H}_i$ before hitting $\mathcal{H}_{i+l}$. 
For $m,M, a\ge 0$ and $l\ge 1$, set
\[
H_{m,l}=\sum_{i=0}^{l-1}N_{m+i,m+l}/(i+1)^2
\]
and
\[
E_{M,l}(a)=\frac{\#\{0\le m\le M: h_{m,l}\le a \text{ and } H_{m,l}\le a\}}{M+1}.
\]
Note that $E_{M,l}(a)$ decreases in $l$ and increases in $a$.

The following lemma is
a minor adaptation of \cite[Lemma 3]{Ze}.
\begin{lemma}\label{LVl5}
For any path $X_\cdot$ with $\varlimsup_{n\to\infty}X_n\cdot e_1/n>0$, 
\begin{equation}\label{LVe27}
\sup_{a\ge 0}\inf_{l\ge 1}\varlimsup_{M\to\infty}E_{M,l}(a)>0.
\end{equation}
\end{lemma}
\pf
Since $\varlimsup_{n\to\infty}n/T_n=\varlimsup_{n\to\infty}X_n\cdot e_1/n>0$,
there exist an increasing sequence $(n_k)_{k=0}^\infty$ and $\delta<\infty$ such that
\[
T_{n_k}<\delta n_k \text{ for all }k.
\]
Thus for any $m$ such that $n_k/2\le m\le n_k$,
\begin{equation}\label{LV*17}
T_m\le 2\delta m.
\end{equation}
Set $M_k=\lceil n_k/2\rceil$, where $\lceil x\rceil\in\mathbb{N}$ denotes the
smallest integer which is not smaller than $x$. Then for all $k$ and
$1<l<\lfloor n_k/2 \rfloor$,
\begin{align}\label{LVe28}
\sum_{m=0}^{M_k} H_{m,l}
&=\sum_{i=0}^{l-1}\Big(\sum_{m=0}^{M_k}N_{m+i,m+l}\Big)/(i+1)^2\nonumber\\
&\le \sum_{i=0}^{l-1} T_{M_k+l}/(i+1)^2
\stackrel{(\ref{LV*17})}{\le} 4\delta (M_k+l).
\end{align}

By the same argument as in Page 193-194 of \cite{Ze}, we will show that there exist
constants $c_1, c_2>0$ such that
\begin{equation}\label{LV*18}
\inf_{l\ge 1}\varlimsup_{k\to\infty}
\frac{\#\{0\le m\le M_k: h_{m,l}\le c_1 \}}{M_k+1}>c_2.
\end{equation}
Indeed, if (\ref{LV*18}) fails, 
then for any $u>0$, 
\[
	\varlimsup_{k\to\infty}\dfrac{\#\{0\le m\le M_k, h_{m,l}\le u\}}{M_k+1}\longrightarrow 0 
\]
as $l\to\infty$ (note that the right side is decreasing in $l$). Hence,
one can find a sequence $(l_i)_{i\ge 0}$ with $l_{i+1}>l_i, l_0=0,$
such that for all $i\ge 0$,
\begin{equation}\label{LV*19}
\varlimsup_{k\to\infty}\dfrac{\#\{0\le m\le M_k, h_{m,l_{i+1}}\le 6\delta l_i\}}{M_k+1}<\frac{1}{3}.
\end{equation}
On the other hand, for $i\ge 0$
\begin{align}
&\varlimsup_{k\to\infty}\dfrac{\#\{0\le m\le M_k, h_{m,l_i}\ge 6\delta l_i\}}{M_k+1}\nonumber\\
&\le 
\varlimsup_{k\to\infty}\frac{1}{(M_k+1)6\delta l_i}\sum_{m=0}^{M_k}(T_{m+l_i}-T_m)\nonumber\\
&\le 
\varlimsup_{k\to\infty}\frac{l_i T_{M_k+l_i}}{6\delta l_i(M_k+1)}
\stackrel{(\ref{LV*17})}{\le}\frac{1}{3}. \label{LV*20}
\end{align}
By (\ref{LV*19}) and (\ref{LV*20}) 
, for any $i\ge 0$,
\begin{equation}\label{LV*21}
\varlimsup_{k\to\infty}\dfrac{\#\{0\le m\le M_k, h_{m,l_{i+1}}> h_{m,l_i}\}}{M_k+1}
\ge \frac{1}{3}.
\end{equation}
Therefore, for any $j\ge 1$, noting that 
\[
\sum_{i=0}^{j-1}1_{h_{m,l_{i+1}}>h_{m,l_i}}\le N_{m,m+l_j}\le H_{m,l_j},
\]
 we have 
\begin{align*}
\frac{j}{3}&\stackrel{(\ref{LV*21})}{\le}
\varlimsup_{k\to\infty}
\sum_{i=0}^{j-1}\dfrac{\#\{0\le m\le M_k, h_{m,l_{i+1}}> h_{m,l_i}\}}{M_k+1}\\
&\le\varlimsup_{k\to\infty}
\frac{1}{M_k+1}\sum_{m=0}^{M_k}H_{m,l_j}\stackrel{(\ref{LVe28})}{\le} 4\delta,
\end{align*}
which is a contradiction if $j$ is large. This proves (\ref{LV*18}).

It follows from (\ref{LV*18}) that, for any $l\ge 1$, there is a subsequence
$(M'_k)$ of $(M_k)$ such that
\[
\frac{\#\{0\le m\le M'_k: h_{m,l}\le c_1 \}}{M'_k+1}>c_2
\]
for all $k$. 
Letting $c_3=9\delta/c_2$, we have that when $k$ is large enough,
\[
\frac{1}{M'_k+1}\sum_{m=0}^{M'_k}1_{h_{m,l}\le c_1, H_{m,l}>c_3}
\le 
\frac{1}{c_3(M'_k+1)}\sum_{m=0}^{M'_k}H_{m,l}
\stackrel{(\ref{LVe28})}{\le}
\frac{c_2}{2}.
\]
Hence for any $l>1$ and large $k$,
\begin{align*}
E_{M_k',l}(c_1\vee c_3)
&\ge \frac{1}{M'_k+1}\sum_{m=0}^{M'_k}1_{h_{m,l}\le c_1,H_{m,l}\le c_3}\\
&= 
\frac{1}{M'_k+1}\sum_{m=0}^{M'_k}(1_{h_{m,l}\le c_1}-1_{h_{m,l}\le c_1,H_{m,l}> c_3})
\ge \frac{c_2}{2}.
\end{align*}
This shows the lemma, and what is more, with explicit constants.\qed\\

For $i\ge 0$, let $N_i=\lim_{j\to\infty} N_{i,j}$ denote the total number of visits to $\mathcal{H}_i$.
With Lemma \ref{LVl5}, one can deduce that with positive probability, $N_i\le C(i+1)^2$ for all $i\ge 0$:
\begin{theorem}\label{LVthm3}
If $\mathbb{P}(\varlimsup_{n\to\infty}X_n\cdot e_1/n>0)>0$, then there exists a constant
$c_5$ such that
\[\mathbb{P}(R=\infty)>0,\]
where $R$ is the stopping time defined by
\begin{align*}
R&=R_{e_1}(X_\cdot, c_5)\\
&:=
\inf\{n\ge 0: \sum_{i=0}^n 1_{X_i\in\mathcal{H}_j}>c_5(j+1)^2 \text{ for some }j\ge 0\}\wedge D,
\end{align*}
and $D:=\inf\{n\ge 1: X_n\cdot e_1\le X_0\cdot e_1\}$.
\end{theorem}

\begin{figure}[h]
\centering
\includegraphics[width=0.7\textwidth]{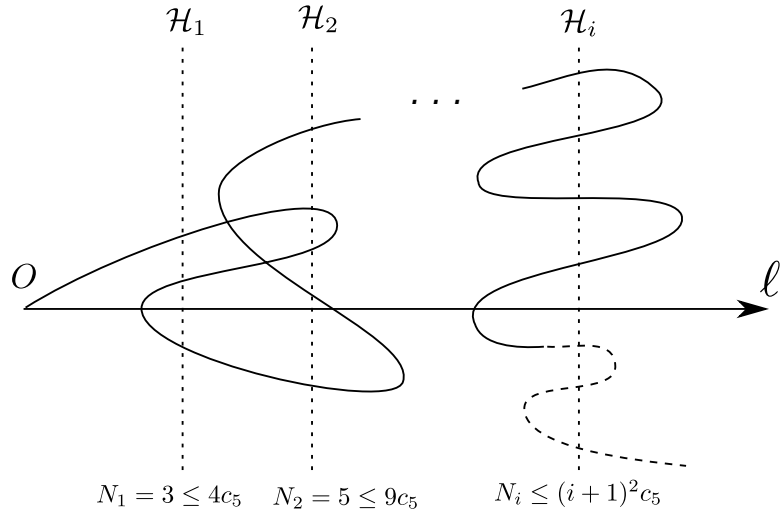}
\caption{On $\{R=\infty\}$, the path visits the $i$-th level no more than
$c_5(i+1)^2$ times.}
\end{figure}

Note that for any $L>0$ and a path $(X_i)_{i=0}^\infty$ with $X_0=o$,
\begin{align}\label{LVe32}
\sum_{\substack{y: y\cdot e_1\le -L\\0\le i\le R}}e^{-\gamma d(y,X_i)}
&\le
\sum_{j=0}^\infty (\#\text{visits to $\mathcal{H}_j$ before time $R$})e^{-\gamma (j+L)}\nonumber\\
&\le 
C\sum_{j=0}^\infty c_5(j+1)^2e^{-\gamma(j+L)}\le Ce^{-\gamma L}.
\end{align}
Hence on the event $\{R=\infty\}$, by (\ref{LVe32}) and $(G)$,
 the trajectory $(X_i)_{i=0}^\infty$ is ``almost 
independent" with the environments $\{\omega_x:x\cdot e_1\le -L\}$ when $L$ is large.
This fact will be used in our definition of the regeneration times in the Section \ref{seclln}.

To prove Theorem \ref{LVthm3}, we need the following lemma. 
Recall that $r,\gamma$ are parameters of the environment measure $P$.
Let $S$ be a countable set of finite paths.
With abuse of notation, we
also use $S$ as the synonym for the event
\begin{equation}\label{LV2e*}
\bigcup_{(x_i)_{i=0}^N\in S}\{X_i=x_i \text{ for }0\le i\le N\}.
\end{equation}
\begin{lemma}\label{LVc2}
Let $a>0$ and $A\subset\Lambda\subset\mathbb{Z}^d$.
Suppose $S\neq\emptyset$ is a countable set of finite paths 
$x_\cdot=(x_i)_{i=0}^N, N<\infty$ that satisfy $d(x_\cdot, \Lambda)\ge r$ and
\[
\sum_{y\in A, 0\le i\le N}e^{-\gamma d(y,x_i)}\le a.
\]
Then, $P$-almost surely,
\begin{equation}\label{LVe31}
	\exp(-Ca)\le\frac{E_P [P_\omega(S)|\omega_x: x\in\Lambda]}{E_P [P_\omega(S)|\omega_x: x\in\Lambda\setminus A]}\le\exp(Ca).
\end{equation}
\end{lemma}

\pf
We shall first show that for any $(x_i)_{i=0}^N\in S$, $P$-almost surely,
\begin{align}\label{LVe40}&E_P[P_\omega(X_i=x_i,0\le i\le N)|\omega_y:y\in\Lambda]\nonumber\\&\le \exp(Ca)E_P[P_\omega(X_i=x_i,0\le i\le N)|\omega_y:y\in\Lambda\setminus A].\end{align}
Note that when $\Lambda^c$ is a finite subset of $\mathbb{Z}^d$,  (\ref{LVe40})
is an easy consequence of $(G)$. For general $\Lambda$, we let
\[
\Lambda_n=\Lambda\cup\{x:|x|\ge n\}.
\]
When $n$ is sufficiently big, $(G)$ implies that 
\begin{equation*}
	\frac{E_P[P_\omega(X_i=x_i,0\le i\le N)|\omega_y:y\in\Lambda_n]}{E_P[P_\omega(X_i=x_i,0\le i\le N)|\omega_y:y\in\Lambda_n\setminus A]}\le\exp(Ca).
\end{equation*}
Since $\Lambda_n\downarrow \Lambda$ as $n\to\infty$, 
(\ref{LVe40}) follows by taking $n\to\infty$ in the above inequality.

Summing over all $(x_i)_{i=0}^N\in S$ on both sides of (\ref{LVe40}), we conclude that
$P$-almost surely,
\[
E_P[P_\omega(S)|\omega_y:y\in\Lambda]
\le \exp(Ca)E_P[P_\omega(S)|\omega_y:y\in\Lambda\setminus A]. 
\]
The upper bound of (\ref{LVe31}) is proved. The lower bound follows likewise.\qed

Now we can prove the theorem. Our proof is a modification
of the proof of Theorem 1 in \cite{Ze}:\\

\noindent\textit{Proof of Theorem \ref{LVthm3}:}
It follows by Lemma \ref{LVl5} that there exists a constant $c_4>0$ such that
\begin{equation}\label{LV*22}
\mathbb{P}(\inf_{l\ge 1}\varlimsup_{M\to\infty}E_{M,l}(c_4)>0)>0.
\end{equation}
For $l>r$, $k\ge 0$ and $z\in\mathbb{Z}^d$ with $z\cdot e_1=r$, let $B_{m,l}(z,k,c)$ denote the event
\[
\{N_{m+r,m+l}=k,X_{T_{m+r,m+l}^k}=X_{T_m}+z,H_{m+r,l-r}\le c\}.
\]

Note that on the event $\{h_{m,l}\le c_4\text{ and }H_{m,l}\le c_4\}$, we have
\begin{align*}
T_{m+r,m+l}^{N_{m+r,m+l}}-T_m
&\le h_{m,l}+\sum_{i=0}^r N_{m+i,m+l}\\
&\le c_4+\sum_{i=0}^r (i+1)^2c_4\le (1+r)^3c_4,
\shortintertext{and}
H_{m+r,l-r}
&\le \sum_{i=0}^{l-r-1}(r+1)^2N_{m+r+i,m+l}/(r+i+1)^2\\
&\le (r+1)^2c_4=:c_5.
\end{align*}
Hence
$
\{h_{m,l}\le c_4\text{ and }H_{m,l}\le c_4\}\subset \bigcup_{|z|,k\le (r+1)^3c_4}B_{m,l}(z,k,c_5),
$
and
\[
\lim_{l\to\infty}\varlimsup_{M\to\infty}E_{M,l}(c_4)
\le 
\sum_{|z|,k\le (r+1)^3c_4}
\varlimsup_{l\to\infty}\varlimsup_{M\to\infty}\frac{1}{M+1}
\sum_{m=0}^M 1_{B_{m,l}(z,k,c_5)}.
\]
Thus by (\ref{LV*22}), for some $k_0$ and $z_0$ with $z_0\cdot e_1=r$,
\begin{equation}\label{LVe29}
\mathbb{P}(\varlimsup_{l\to\infty}\varlimsup_{M\to\infty}
\frac{1}{M+1}\sum_{m=0}^M 1_{B_{m,l}(z_0,k_0,c_5)}>0)>0.
\end{equation}
In what follows, we write $B_{m,l}(z_0,k_0,c_5)$ simply as
$B_{m,l}$. 

For any $l>r$ and any fixed $i\le l-1$, let $m_j=m_j(l,i):=i+jl$, i.e. $(m_j)_{j\ge 0}$ is
the class of residues of $i(\text{mod }l)$.
Now take any $j\in \mathbb{N}$.	Observe that for any event $E=\{1_{B_{m_{j-1},l}}=\cdot,\ldots,1_{B_{m_0,l}}=\cdot\}$
and $x\in \mathcal{H}_{m_j}$,
\begin{align}\label{LV*1}
\MoveEqLeft P_\omega(\{X_{T_{m_j}}=x\}\cap E\cap B_{m_j,l})\\
&\le 
P_\omega(\{X_{T_{m_j}}=x\}\cap E)
P_\omega^{x+z_0}(D>T_{l-r},H_{0,l-r}\le c_5).\nonumber
\end{align} 
Moreover, for any $x\in \mathcal{H}_{m_j}$, there exists a countable set $S$ of finite paths $(x_i)_{i=0}^N$ that satisfy
$m_j+r\le x_i\cdot e_1\le m_j+l$ and $\#\{k\le N: x_k\in\mathcal{H}_i(x_0)\}\le c_5(i+1)^2$ for
$0\le i\le N$, such that
\begin{align*}
&\{X_0=x+z_0, D>T_{l-r},H_{0,l-r}\le c_5\}\\
&=\cup_{(x_i)_{i=0}^N\in S}\{X_i=x_i
\text{ for }0\le i\le N\}.
\end{align*}
Noting that (by the same argument as in (\ref{LVe32})) for any $(x_i)_{i=0}^N\in S$,
\[
\sum_{\substack{y:y\cdot e_1\le m_j\\i\le N}}e^{-\gamma d(y,x_i)}\le Ce^{-\gamma r},
\]
by Lemma \ref{LVc2} we have
\begin{align*}
&E_P[P_\omega^{x+z_0}(D>T_{l-r},H_{0,l-r}\le c_5)|\omega_y:y\cdot e_1\le m_j]\\
&\le 
\exp{(Ce^{-\gamma r})}\mathbb{P}(D>T_{l-r},H_{0,l-r}\le c_5).
\end{align*}
Thus for $j\ge 0$ and $l>r$,
\begin{align*}
&\mathbb{P}(E\cap B_{m_j,l})\\
&\stackrel{(\ref{LV*1})}{\le}
\sum_{x\in\mathcal{H}_{m_j}}
E_P \big[P_\omega(\{X_{T_{m_j}}=x\}\cap E)P_\omega^{x+z_0}(D>T_{l-r},H_{0,l-r}\le c_5)\big]\\
&\le 
\exp{(Ce^{-\gamma r})}
\sum_{x\in\mathcal{H}_{m_j}}
\mathbb{P}(\{X_{T_{m_j}}=x\}\cap E)\mathbb{P}(D>T_{l-r},H_{0,l-r}\le c_5)\\
&=
C\mathbb{P}(E)
\mathbb{P}(D>T_{l-r},H_{0,l-r}\le c_5).
\end{align*}
Hence, for any $j\ge 0$ and $l>r$,
\begin{equation*}
\mathbb{P}(1_{B_{m_j,l}}=1|1_{B_{m_{j-1},l}},\ldots,1_{B_{m_0,l}})
\le 
C\mathbb{P}(D>T_{l-r},H_{0,l-r}\le c_5),
\end{equation*}
which implies that $\mathbb{P}$-almost surely,
\begin{equation}\label{LVe30}
\varlimsup_{n\to\infty}
\frac{1}{n}\sum_{j=0}^{n-1} 1_{B_{m_j,l}}
\le C\mathbb{P}(D>T_{l-r},H_{0,l-r}\le c_5).
\end{equation}
Therefore, $\mathbb{P}$-almost surely,
\begin{align*}
\varlimsup_{l\to\infty}\varlimsup_{M\to\infty}
\frac{1}{M+1}\sum_{m=0}^M 1_{B_{m,l}}
&\le \varlimsup_{l\to\infty}\frac{1}{l}\sum_{i=0}^{l-1}
\varlimsup_{M\to\infty}\frac{l}{M+1}
\sum_{\substack{0\le m\le M\\m\text{ mod }l=i}} 1_{B_{m,l}}\\
&\stackrel{(\ref{LVe30})}{\le} \lim_{l\to\infty}
C\mathbb{P}(D>T_{l-r},H_{0,l-r}\le c_5)\\
&=C\mathbb{P}(D=\infty, \sum_{i=0}^\infty N_i/(i+1)^2\le c_5).
\end{align*}
This and (\ref{LVe29}) yield
$\mathbb{P}(D=\infty, \sum_{i=0}^\infty N_i/(i+1)^2\le c_5)>0$. 
The theorem follows.\qed

\section{The conditional law of large numbers}\label{seclln}
In this section we will prove the conditional law of large numbers \eqref{ICLLN}, using
regeneration times and coupling. 
Given the dependence structure of the environment, we want 
to define regeneration times in such a way that what happens after a regeneration time has little 
dependence on the past. To this end, we will use the ``$\epsilon$-coins" trick introduced in \cite{CZ1} and the stopping time $R$ to define the regeneration times. 
Intuitively, at a regeneration time, the past and the future movements have nice properties. That
is, the walker has walked straight for a while without paying attention to the environment, and his
future movements have little dependence on his past movements.

We define the $\epsilon$-coins $(\epsilon_{i,x})_{i\in\mathbb{N}, x\in \mathbb{Z}^d}=:\epsilon$
to be iid random variables with distribution
$Q$ such that 
\[Q(\epsilon_{i,x}=1)=d\kappa \text{ and }Q(\epsilon_{i,x}=0)=1-d\kappa.\]

For fixed $\omega$, $\epsilon$, $P_{\omega,\epsilon}^x$ is the law of the Markov chain $(X_n)$ such that $X_0=x$ and that for any $e\in\mathbb{Z}^d$
such that $|e|=1$,
\[
P_{\omega,\epsilon}^x(X_{n+1}=z+e|X_n=z)
=\frac{1_{\epsilon_{n,z}=1}}{2d}+\frac{1_{\epsilon_{n,z}=0}}{1-d\kappa}[\omega(z,z+e)-\frac{\kappa}{2}].
\]
Note that the law of $X_\cdot$ under $\bar{P}_\omega^x=Q\otimes P_{\omega,\epsilon}^x$ coincides with its
law under $P_\omega^x$. Sometimes we also refer to 
$P_{\omega,\epsilon}^x (\cdot)$ as a measure on the sets of paths, without indicating the specific random path. 

Denote by
$\bar{\mathbb{P}}=P\otimes Q\otimes P_{\omega,\epsilon}^o$ the law of the triple 
$(\omega, \epsilon, X_\cdot)$. 

Now we define the regeneration times in the direction $e_1$. 
Let $L$ be a fixed number which is sufficiently large.
Set $R_0=0$. Define inductively for $k\ge 0$:
\begin{align*}
&S_{k+1}=\inf\{n\ge R_k: X_{n-L}\cdot e_1>\max\{X_m\cdot e_1: m<n-L\},\\
&\qquad\qquad\qquad \epsilon_{n-i, X_{n-i}}=1, X_{n-i+1}-X_{n-i}=e_1 \text{ for all }1\le i\le L\},\\
&R_{k+1}=R\circ\theta_{S_{k+1}}+S_{k+1},
\end{align*}
where $\theta_n$ denotes the time shift of the path, i.e.,
$\theta_n X=(X_{n+i})_{i=0}^\infty$.

Let \[K=\inf\{k\ge 1: S_k<\infty,R_k=\infty\}\] and $\tau_1=\tau_1(e_1,\epsilon,X_\cdot):=S_K.$
For $k\ge 1$, the ($L$-)regeneration times are defined inductively by \[\tau_{k+1}=\tau_1\circ\theta_{\tau_k}+\tau_k .\]

By similar argument as in \cite[Lemma 2.2]{CZ1}, we can show:
\begin{lemma}\label{LVl7}
If $\mathbb{P}(\lim_{n\to\infty}X_n\cdot e_1/n=0)<1$, then
\begin{equation}\label{LVe19}
\mathbb{P}(A_{e_1}\cup A_{-e_1})=1.
\end{equation}
Moreover, on $A_{e_1}$, $\tau_i$'s are $\bar{\mathbb P}$-almost surely finite.
\end{lemma}
\pf
If $\mathbb{P}(\lim_{n\to\infty}X_n\cdot e_1/n=0)<1$, 
\[
\mathbb{P}(\varlimsup_{n\to\infty}X_n\cdot e_1/n>0)>0\quad\text{ or }\quad
\mathbb{P}(\varlimsup_{n\to\infty}X_n\cdot (-e_1)/n>0)>0.
\]
Without loss of generality, assume that 
\[\mathbb{P}(\varlimsup_{n\to\infty}X_n\cdot e_1/n>0)>0.\]
 It then follows from Theorem \ref{LVthm3} that $\mathbb{P}(R=\infty)>0$. 
 We want to show that 
$R_k=\infty$ for all but finitely many $k$'s.

For $k\ge 0$,
\begin{align*}
&\bar{\mathbb P}(R_{k+1}<\infty)\\
&= \bar{\mathbb P}(S_{k+1}<\infty, R\circ\theta_{S_{k+1}}<\infty)\\
&=\sum_{n,x}\bar{\mathbb P}
(S_{k+1}=n, X_n=x, R\circ\theta_n <\infty)\\
&=\sum_{n,x}E_{P\otimes Q}
\big[P_{\omega,\epsilon}(S_{k+1}=n,X_n=x)
P_{\omega,\theta^n\epsilon}^x(R<\infty)\big],
\end{align*}
where $\theta^n\epsilon$ denotes the time shift of the coins $\epsilon$, 
i.e. $(\theta^n\epsilon)_{i,x}=\epsilon_{n+i,x}$.
Note that $P_{\omega,\epsilon}(S_{k+1}=n,X_n=x)$ and 
$P_{\omega,\theta^n\epsilon}^x(R<\infty)$ are independent under the measure $Q$,
since the former is a function of $\epsilon$'s before time $n$, and the latter
involves $\epsilon$'s after time $n$.
It then follows by induction that
\begin{align*}
&\bar{\mathbb P}(R_{k+1}<\infty)\\
&=\sum_{n,x}E_P
\big[\bar{P}_\omega(S_{k+1}=n,X_n=x)
\bar{P}_\omega^x(R<\infty)\big]\\
&=\sum_{n,x}E_P\big[\bar{P}_\omega(S_{k+1}=n,X_n=x)
E_P[\bar{P}_\omega^x(R<\infty)|\omega_y: y\cdot e_1\le x\cdot e_1-L]
\big]\\
&\stackrel{(\ref{LVe32}), \text{Lemma }\ref{LVc2}}{\le} \bar{\mathbb P}(R_k<\infty)\exp{(e^{-cL})}\bar{\mathbb P}(R<\infty)\\
&\le [\exp{(e^{-cL})}\bar{\mathbb P}(R<\infty)]^{k+1},
\end{align*}
where we used in the second equality the fact that $\bar{P}_\omega(S_{k+1}=n,X_n=x)$
is $\sigma(\omega_y: y\cdot e_1\le x\cdot e_1-L)$-measurable.	
Hence, by taking $L$ sufficiently large and by the Borel-Cantelli Lemma, $\bar{\mathbb P}$-almost surely,
$R_k=\infty$ except for finitely many values of $k$.

 Let $\mathcal{O}_{e_1}$
denote the event that the signs of $X_n\cdot e_1$ change infinitely many often.
It is easily seen that (by the ellipticity of the environment)
\begin{align*}
	&\mathbb{P}(\mathcal{O}_{e_1}\cup A_{e_1}\cup A_{-e_1})=1
\shortintertext{and}
&\mathcal{O}_{e_1}\subset \{\sup_n X_n\cdot e_1=\infty\}.
\end{align*} 
However, on $\{\sup_n X_n\cdot e_1=\infty\}$, given that
$R_k$ is finite, $S_{k+1}$ is also finite.
Hence $\tau_1$ is $\bar{\mathbb P}$-almost surely finite on $\{\sup_n X_n\cdot e_1=\infty\}$,
and so are the regeneration times $\tau_2,\tau_3\ldots$. 
Therefore,
\[
\mathbb{P}(\mathcal{O}_{e_1})=\bar{\mathbb P}(\mathcal{O}_{e_1}\cap\{\tau_1<\infty\}).
\]
Since $\mathcal{O}_{e_1}\cap\{\tau_1<\infty\}=\emptyset$, we get $\mathbb{P}(\mathcal{O}_{e_1})=0$.	This gives (\ref{LVe19}). \qed\\

When $\mathbb{P}(R=\infty)>0$, we let 
\[
\hat{\mathbb P}(\cdot):=\bar{\mathbb P}(\cdot|R=\infty).
\]
The following proposition is a consequence of Lemma \ref{LVc2}.
\begin{proposition}
Assume $\mathbb{P}(R=\infty)>0$.
Let $l>r$ and
$\Lambda\subset\{x:x\cdot e_1<-r\}.$
Then for any $A\subset\Lambda\cap\{x: x\cdot e_1<-l\}$ and $k\in\mathbb{N}$,
\begin{equation}\label{LVprop1}
\exp(-Ce^{-\gamma l})
\le 
\dfrac{E_P\big[\bar{P}_\omega\big((X_i)_{i=0}^{\tau_k}\in\cdot, R=\infty\big)|\omega_y:y\in\Lambda\setminus A]}
{E_P\big[\bar{P}_\omega\big((X_i)_{i=0}^{\tau_k}\in\cdot, R=\infty\big)|\omega_y:y\in\Lambda]}
\le 
\exp(Ce^{-\gamma l}).
\end{equation}
Furthermore, for any $k\in\mathbb{N}$ and $n\ge 0$, $\hat{\mathbb P}$-almost surely,
\begin{equation}\label{LVprop2}
\exp(-e^{-cL})
\le 
\frac{\hat{\mathbb P}\big((X_{\tau_n+i}-X_{\tau_n})_{i=0}^{\tau_{n+k}-\tau_n}\in\cdot|X_{\tau_n}\big)}
{\hat{\mathbb P}\big((X_i)_{i=0}^{\tau_k}\in\cdot\big)}
\le 
\exp(e^{-cL}).
\end{equation}
\end{proposition}
\pf
First, we shall prove (\ref{LVprop1}).
By the definition of the regeneration times, for any finite path $x_\cdot=(x_i)_{i=0}^N, N<\infty$, there exists an event
$G_{x_\cdot}\in\sigma(\epsilon_{i,X_i},X_i: i\le N)$
 such that $G_{x_\cdot}\subset\{R>N\}$  and
\[
\{(X_i)_{i=0}^{\tau_k}=(x_i)_{i=0}^N, R=\infty\}
=G_{x_\cdot}\cap\{R\circ\theta_N=\infty\}.
\]
(For example, when $k=1$, we let 
\[
G_{x_\cdot}=\bigcup_{j=1}^\infty  
\{(X_i)_{i=0}^N=(x_i)_{i=0}^N, S_j=N, R>N\}.
\]
Then $\{(X_i)_{i=0}^{\tau_1}=(x_i)_{i=0}^N, R=\infty\}
=G_{x_\cdot}\cap\{R\circ\theta_N=\infty\}.$)

For $n\in\mathbb{N}$, we let 
\[
E_n:=G_{x_\cdot}\cap\{R\circ\theta_N\ge n\}.
\]
Note that $E_n\in\sigma(\epsilon_{i,X_i},X_i:i\le N+n)$ can be interpreted (in the sense of (\ref{LV2e*}))
as a set of paths with lengths $\le N+n$. Also note that $E_n\subset\{R>N+n\}$.
Then by Lemma \ref{LVc2} and (\ref{LVe32}), we have
\[
\exp(-Ce^{-\gamma l})
\le 
\dfrac{E_P\big[\bar{P}_\omega\big(E_n)|\omega_y:y\in\Lambda\setminus A]}
{E_P\big[\bar{P}_\omega\big(E_n\big)|\omega_y:y\in\Lambda]}
\le 
\exp(Ce^{-\gamma l}).
\]
(\ref{LVprop1}) follows by letting $n\to\infty$.

Next, we shall prove (\ref{LVprop2}).
Let $x\in\mathbb{Z}^d$ be any point that satisfies 
\[
\bar{\mathbb P}(X_{\tau_n}=x)>0.
\]
By the definition of the regeneration times, 
for any $m\in\mathbb{N}$,
there exists an event
$G_m^x\in\sigma\{\epsilon_{i,X_i},X_i: i\le m\}$
such that $\bar{P}_\omega(G_m^x)$ is 
$\sigma(\omega_y:y\cdot e_1\le x\cdot e_1-L)$-measurable, and
\[
\{\tau_n=m,X_m=x, R=\infty\}=G^x_m\cap\{R\circ\theta_m=\infty\}.
\]
Thus 
\begin{align}\label{LV2e2}
&\bar{\mathbb P}\big((X_{\tau_n+i}-X_{\tau_n})_{i=0}^{\tau_{n+k}-\tau_n}\in\cdot,X_{\tau_n}=x, R=\infty\big)\nonumber\\
&=\sum_{m}\bar{\mathbb P}\big((X_{\tau_n+i}-X_{\tau_n})_{i=0}^{\tau_{n+k}-\tau_n}\in\cdot, \tau_n=m,X_m=x,R=\infty\big)\nonumber\\
&=
\sum_{m}E_P\big[\bar{P}_\omega(G_m^x)
\bar{P}_\omega^x((X_i-x)_{i=0}^{\tau_k}\in\cdot,R=\infty)\big]\nonumber\\
&\stackrel{(\ref{LVprop1})}{\le}
\exp(Ce^{-\gamma L})\sum_{m}\bar{\mathbb P}(G_m^x)
\bar{\mathbb P}\big((X_i)_{i=0}^{\tau_k}\in\cdot, R=\infty\big).
\end{align}
On the other hand, 
\begin{align}\label{LV2e11}
\bar{\mathbb P}(X_{\tau_n}=x, R=\infty)
&=\sum_{m} E_P[\bar{P}_{\omega}(G^x_m)\bar{P}_\omega^x(R=\infty)]\nonumber\\
&\stackrel{(\ref{LVprop1})}{\ge}\exp(-Ce^{-\gamma L})
\sum_{m}
\bar{\mathbb P}(G_m^x)\bar{\mathbb P}(R=\infty).
\end{align}
By (\ref{LV2e2}) and (\ref{LV2e11}), we have (note that $L$ is sufficiently big)
\[
\hat{\mathbb P}\big((X_{\tau_n+i}-X_{\tau_n})_{i=0}^{\tau_{n+k}-\tau_n}\in\cdot|X_{\tau_n}=x\big)
\le 
\exp(e^{-cL})\hat{\mathbb P}\big((X_i)_{i=0}^{\tau_k}\in\cdot\big).
\]
The right side of (\ref{LVprop2}) is proved. The left side of (\ref{LVprop2}) follows likewise.
\qed

The next lemma describes the dependency of a regeneration on its remote past.
It is a version of Lemma 2.2 in \cite{CZ2}. (The denominator is omitted in the 
last equality in \cite[page 101]{CZ2}, which is corrected here, see the equality in (\ref{LV*8}).)

Set $\tau_0=0$.
Denote the truncated path between $\tau_{n-1}$ and $\tau_n-L$ by 
\[
P_n=(P_n^i)_{0\le i\le \tau_{n}-\tau_{n-1}-L}:=(X_{i+\tau_{n-1}}-X_{\tau_{n-1}})_{0\le i\le \tau_n-\tau_{n-1}-L}.
\]
Set
\begin{align*}
 W_n &=(\omega_{x+X_{\tau_{n-1}}})_{x\in P_n}=:\omega_{X_{\tau_{n-1}}+P_n},\\
 F_n &=X_{\tau_n}-X_{\tau_{n-1}},\\
J_n &=(P_n,W_n,F_n,\tau_n-\tau_{n-1}).
\end{align*}
For $i\ge 0$, let $h_{i+1}(\cdot|j_i,\ldots,j_1):=\hat{\mathbb{P}}(J_{i+1}\in\cdot|J_{i},\ldots,J_1)|_{J_{i}=j_i,\ldots,J_1=j_1}$ denote the transition kernel of $(J_n)$. Note that when $i=0$, $h_{i+1}(\cdot|j_i,\ldots,j_1)=h_1(\cdot|\emptyset)=\hat{\mathbb P}(J_1\in\cdot)$.
\begin{lemma}\label{LVl4} 
Assume $\mathbb{P}(R=\infty)>0$, $0\le k\le n$. Then $\hat{\mathbb P}$-almost surely,
\begin{equation}\label{LVe26}
\exp{(-e^{-c(k+1)L})}
\le 
\frac{h_{n+1}(\cdot|J_n,\ldots,J_1)}{h_{k+1}(\cdot|J_n,\ldots,J_{n-k+1})}
\le \exp{(e^{-c(k+1)L})}.
\end{equation}
\end{lemma}

\pf 
For $j_m=(p_m,w_m,f_m,t_m),m=1,\ldots n$, let
\begin{align*}
&\bar{x}_m:=f_1+\cdots+f_m,\\
&\bar{t}_m:=t_1+\cdots+t_m,\\
&B_{p_1,\ldots,p_m}:=\{R=\infty, P_i=p_i\text{ for all }i=1,\ldots,m\},\\
\text{ and }\quad
& \omega_{p_1,\ldots,p_m}:=(\omega_{\bar{x}_{i-1}+p_i})_{i=1}^m.
\end{align*}

First, we will show that for any $1\le k\le n$,
\begin{equation}\label{LV*8}
h_{k+1}(\cdot|j_k,\ldots,j_1)
=\frac{E_P\big[\bar{P}_\omega^{\bar{x}_k}(J_1\in \cdot,R=\infty)|\omega_{p_1,\ldots,p_k}\big]}
{E_P\big[\bar{P}_\omega^{\bar{x}_k}(R=\infty)|\omega_{p_1,\ldots,p_k}\big]}
\Big|_{\omega_{p_1,\ldots,p_k}=(w_i)_{i=1}^k}.
\end{equation}
By the definition of the regeneration times, there exists an event
\[
G_{p_1,\ldots,p_k}\in\sigma(X_{i+1},\epsilon_{i,X_i}, 0\le i\le \bar{t}_k-1)
\]
such that 
\begin{equation}\label{LV*7}
B_{p_1,\ldots,p_k}=G_{p_1,\ldots,p_k}\cap\{R\circ\theta_{\bar{t}_k}=\infty\}.
\end{equation}
On the one hand, for any $\sigma(J_k,\ldots, J_1)$-measurable function $g(J_k,\ldots,J_1)$,
\begin{align}\label{LV*15}
&E_{\bar{\mathbb P}}\big[h_{k+1}(\cdot|J_k,\ldots,J_1)g(J_k,\ldots, J_1)1_{B_{p_1,\ldots,p_k}}\big]\nonumber\\
&=E_{\bar{\mathbb P}}\big[g1_{B_{p_1,\ldots,p_k}}1_{J_{k+1}\in \cdot}\big]\nonumber\\
&=
E_P [g1_{B_{p_1,\ldots,p_k}}\bar{P}_\omega(J_{k+1}\in\cdot,B_{p_1,\ldots,p_k})]\nonumber\\
&\stackrel{(\ref{LV*7})}{=}
E_P \big[g1_{B_{p_1,\ldots,p_k}}\bar{P}_\omega(G_{p_1,\ldots,p_k})\bar{P}_\omega^{\bar{x}_k}(J_1\in \cdot,R=\infty)\big].
\end{align}
On the other hand, we also have
\begin{align}\label{LV*16}
& E_{\bar{\mathbb P}}\big[h_{k+1}(\cdot|J_k,\ldots,J_1)
g(J_k,\ldots, J_1)1_{B_{p_1,\ldots,p_k}}\big]\nonumber\\
&=
E_P \big[
h_{k+1}(\cdot|J_k,\ldots,J_1)g1_{B_{p_1,\ldots,p_k}}
\bar{P}_\omega(B_{p_1,\ldots,p_k})\big]\nonumber\\
&\stackrel{(\ref{LV*7})}{=}
E_P \big[
h_{k+1}(\cdot|J_k,\ldots,J_1)g1_{B_{p_1,\ldots,p_k}}
\bar{P}_\omega(G_{p_1,\ldots,p_k})\bar{P}_\omega^{\bar{x}_k}(R=\infty)\big].
\end{align}
Comparing (\ref{LV*15}) and (\ref{LV*16}) and observing that on $B_{p_1,\ldots,p_k}$, 
$\bar{P}_\omega(G_{p_1,\ldots,p_k})$ and all functions of $J_1,\ldots,J_k$
are $\sigma(\omega_y: y\in\bar{x}_{i-1}+p_i, i\le k)$-measurable
, we obtain that on $B_{p_1,\ldots,p_k}$, $P$-almost surely,
\begin{equation*}
h_{k+1}(\cdot|J_k,\ldots,J_1)
=\frac{E_P\big[\bar{P}_\omega^{\bar{x}_k}(J_1\in \cdot,R=\infty)|\omega_{\bar{x}_{i-1}+p_i},i\le k\big]}
{E_P\big[\bar{P}_\omega^{\bar{x}_k}(R=\infty)|\omega_{\bar{x}_{i-1}+p_i},i\le k\big]}.
\end{equation*}
Noting that
\[
B_{p_1,\ldots,p_k}\cap\{\omega_{p_1,\ldots,p_k}=(w_i)_{i=1}^k\}
=
\{J_i=j_i,1\le i\le k\},
\]
(\ref{LV*8}) is proved.

Next, we will prove the lower bound in \eqref{LVe26}.

When $n\ge k\ge 1$,	by formula (\ref{LV*8}) and (\ref{LVprop1}), we have
\begin{align}\label{LV3e2}
&h_{n+1}(\cdot| j_n,\ldots,j_1)\nonumber\\
&=
\frac{E_P[\bar{P}_\omega^{\bar{x}_n}(J_1\in\cdot,R=\infty)|\omega_{p_1,\ldots,p_n}]}
{E_P\big[\bar{P}_\omega^{\bar{x}_n}(R=\infty)|\omega_{p_1,\ldots,p_n}\big]}\bigg|_{\omega_{p_1,\ldots,p_n}=(w_i)_{i=0}^n}\nonumber\\
&\le
\frac{\exp(Ce^{-\gamma(k+1)L})
E_P[\bar{P}_\omega^{\bar{x}_n}(J_1\in\cdot,R=\infty)|\omega_{\bar{x}_{i-1}+p_i},n-k+1\le i\le n]}
{\exp(-Ce^{-\gamma(k+1)L})
E_P[\bar{P}_\omega^{\bar{x}_n}(R=\infty)|\omega_{\bar{x}_{i-1}+p_i},n-k+1\le i\le n]}\nonumber\\
&\qquad
\big|_{\omega_{p_1,\ldots,p_n}=(w_i)_{i=0}^n}\nonumber\\
&=
\exp(2Ce^{-\gamma(k+1)L})
\frac{E_P[\bar{P}_\omega^{\bar{x}_n-\bar{x}_{n-k}}(J_1\in\cdot,R=\infty)|\omega_{p_{n-k+1},\ldots, p_n}]}
{E_P[\bar{P}_\omega^{\bar{x}_n-\bar{x}_{n-k}}(R=\infty)|\omega_{p_{n-k+1},\ldots, p_n}]}\nonumber\\
&\qquad
\big|_{\omega_{p_{n-k+1},\ldots,p_n}=(w_i)_{i=n-k+1}^n}\nonumber\\
&\stackrel{(\ref{LV*8})}{=}
\exp(2Ce^{-\gamma(k+1)L})h_{k+1}(\cdot|j_{n},\ldots,j_{n-k+1}),
\end{align}
where we used the translation invariance of the measure $P$ in the last but one equality.

When $k=0$ and $n\ge 1$, by formula \eqref{LV*8} and \eqref{LVprop1},
\begin{align}\label{LV3e1}
	h_{n+1}(\cdot| j_n,\ldots,j_1)
	&\le 
	\frac{\exp(Ce^{-\gamma L})E_P[\bar{P}_\omega^{\bar{x}_n}(J_1\in\cdot,R=\infty)]}
	{\exp(-Ce^{-\gamma L})E_P[\bar{P}_\omega^{\bar{x}_n}(R=\infty)]}\nonumber\\
	&=\exp(2Ce^{-\gamma L})\hat{\mathbb P}(J_1\in\cdot)\nonumber\\
	&=\exp(2Ce^{-\gamma L})h_1(\cdot|\emptyset).
\end{align}

When $k=n=0$, \eqref{LVe26} is trivial.
Hence combining \eqref{LV3e2} and \eqref{LV3e1}, the lower bound in (\ref{LVe26}) follows as we take $L$ sufficiently big. The upper bound follows likewise.\qed
 
\begin{lemma}\label{LVl6}
Suppose that a sequence of non-negative random variables $(X_n)$ satisfies
\[
a\le \frac{\ud P(X_{n+1}\in\cdot|X_1,\ldots, X_n)}{\ud \mu}\le b
\]
for all $n\ge 1$, where $a\le 1\le b$ are constants and $\mu$ is a probability measure. Let $m_\mu\le \infty$ be
the mean of $\mu$. Then almost surely,
\begin{equation}\label{LVe21}
a m_\mu\le\varliminf_{n\to\infty}
\frac{1}{n}\sum_{i=1}^n X_i
\le
\varlimsup_{n\to\infty}
\frac{1}{n}\sum_{i=1}^n X_i
\le b m_\mu.
\end{equation}
\end{lemma}

Before giving the proof, let us recall the ``splitting representation" of random variables:
\begin{proposition}\cite[Page 94]{Tho}\label{LVprop}
Let $\nu$ and $\mu$ be probability measures. Let $X$ be a random variable with law $\nu$. If 
for some $a\in(0,1)$,
\[
\frac{\ud\nu}{\ud\mu}\ge a,
\]
then, enlarging the probability space if necessary, we can find independent
random variables $\Delta, \pi, Z$ such that
\begin{itemize}
\item[i)] $\Delta$ is Bernoulli with parameter $1-a$, i.e., $P(\Delta=1)=1-a$,
$P(\Delta=0)=a$;
\item[ii)] $\pi$ is of law $\mu$, and $Z$ is of law $(\nu-a\mu)/(1-a)$;
\item[iii)] $X=(1-\Delta)\pi+\Delta Z$.
\end{itemize}
\end{proposition}

\noindent{\it Proof of Lemma \ref{LVl6}}:\\
By Proposition \ref{LVprop}, enlarging the probability space if necessary, there are
random variables $\Delta_i,\pi_i,Z_i,i\ge 1$,
such that for any $i\in\mathbb{N}$,
\begin{itemize}
\item 
$\Delta_i$ is Bernoulli with parameter $(1-a)$, and $\pi_i$ is of law $\mu$;
\item 
$\Delta_i, \pi_i$ and $Z_i$ are mutually independent;
\item 
$(\Delta_i,\pi_i)$
is independent of
$\sigma(\Delta_k,\pi_k,Z_k: k<i)$;
\item
$X_i=(1-\Delta_i)\pi_i+\Delta_i Z_i$.
\end{itemize}
Note that since $X_i$'s are supported on $[0,\infty)$, $\pi_i\ge 0$ and $Z_i\ge 0$ for all $i\in\mathbb{N}$.
Thus by the law of large numbers, almost surely,
\[
\varliminf_{n\to\infty}\frac{1}{n}\sum_{i=1}^n X_i\ge
\lim_{n\to\infty}\frac{1}{n}\sum_{i=1}^n (1-\Delta_i)\pi_i=a m_\mu.
\]
This proves the first inequality of (\ref{LVe21}). 

If $m_\mu=\infty$, the last inequality
of (\ref{LVe21}) is trivial. Assume that $m_\mu<\infty$.
Let $(\tilde{\Delta}_i)_{i\ge 1}$ be an iid Bernoulli sequence with
parameter $1-b^{-1}$ such that every $\tilde{\Delta}_i$ is independent of
all the $X_n$'s.
By a similar splitting procedure, we can construct non-negative
random variables $\tilde{\pi}_i,\tilde{Z}_i, i\ge 1$,
such that $(\tilde{\pi}_i)_{i\ge 1}$ are iid with law $\mu$, and
\[
\tilde{\pi}_i=(1-\tilde{\Delta}_i)X_i+\tilde{\Delta}_i\tilde{Z}_i.
\]

Let $Y_i=(1-b^{-1}-\tilde{\Delta}_i)X_i 1_{X_i\le i}$, we will first show that
\begin{equation}\label{LVe22}
\lim_{n\to\infty}\frac{1}{n}\sum_{i=1}^n Y_i=0.
\end{equation}
By Kronecker's Lemma, it suffices to show that
\[
\sum_{i=1}^\infty \frac{Y_i}{i} \text{ converges.}
\]
Observe that $(\sum_{i=1}^n Y_i/i)_{n\in\mathbb{N}}$ is a martingale sequence.
Moreover, for all $n\in\mathbb{N}$,
\begin{align*}
E\big(\sum_{i=1}^n \frac{Y_i}{i}\big)^2
= \sum_{i=1}^n EY_i^2/i^2
&\le \sum_{i=1}^\infty EX_i^2 1_{X_i\le i}/i^2\\
&\le b \sum_{i=1}^\infty E\tilde{\pi}_i^2 1_{\tilde{\pi}_i\le i}/i^2\\
&= b\int_0^\infty x^2 (\sum_{i\ge x}\frac{1}{i^2})\ud\mu\\
&\le C\int_0^\infty x\ud\mu=Cm_\mu<\infty.
\end{align*}
By the $L^2$-martingale convergence theorem, $\sum Y_i/i$
converges a.s. and in $L^2$. This proves (\ref{LVe22}). 

Since
\[
\sum_i P(Y_i\neq (1-b^{-1}-\tilde{\Delta}_i)X_i)
\le 
\sum_i P(X_i>i)
\le 
b \sum_i P(\pi_1>i)\le b m_\mu<\infty,
\]
by the Borel-Cantelli lemma, it follows from (\ref{LVe22}) that
\[
\lim_{n\to\infty}\frac{1}{n}\sum_{i=1}^n (1-b^{-1}-\tilde{\Delta}_i)X_i=0, \text{a.s.}.
\]
Hence almost surely,
\begin{equation*}
m_\mu=\lim_{n\to\infty}\frac{1}{n}\sum_{i=1}^n \tilde{\pi}_i
\ge 
\varlimsup_{n\to\infty}\frac{1}{n}\sum_{i=1}^n (1-\tilde{\Delta}_i)X_i
=\varlimsup_{n\to\infty}\frac{1}{n}\sum_{i=1}^n b^{-1} X_i.
\end{equation*}
The last inequality of (\ref{LVe21}) is proved.\qed 

\begin{theorem}\label{LVlln}
There exist two deterministic numbers $v_{e_1},v_{-e_1}\ge 0$ such that $\mathbb{P}$-almost surely,
\begin{equation}\label{LVe25}
\lim_{n\to\infty}\frac{X_n\cdot e_1}{n}=v_{e_1} 1_{A_{e_1}}-v_{-e_1}1_{A_{-e_1}}.
\end{equation}
Moreover, if $v_{e_1}>0$, then $E_{\hat{\mathbb P}}\tau_1<\infty$ and 
$\mathbb{P}(A_{e_1}\cup A_{-e_1})=1$.
\end{theorem}
\pf
We only consider the nontrivial case that $\mathbb{P}(\lim X_n\cdot e_1/n=0)<1$,
which by Lemma \ref{LVl7} implies
$\mathbb{P}(A_{e_1}\cup A_{-e_1})=1$. 
Without loss of generality, assume $\mathbb{P}(\varlimsup_{n\to\infty}X_n\cdot e_1/n>0)>0$.
We will show that on $A_{e_1}$,
\[
\lim_{n\to\infty}X_n\cdot e_1/n=v_{e_1}>0, \text{ $\mathbb{P}$-a.s..}
\]

By (\ref{LVprop2}) and Lemma \ref{LVl6}, we obtain that
$\mathbb{P}(\cdot|A_{e_1})$-almost surely,
\begin{align}
\exp{(-e^{-cL})}E_{\hat{\mathbb P}}X_{\tau_1}\cdot e_1
&\le \varliminf_{n\to\infty}\frac{X_{\tau_n}\cdot e_1}{n}\nonumber\\
&\le \varlimsup_{n\to\infty}\frac{X_{\tau_n}\cdot e_1}{n}
\le \exp{(e^{-cL})}E_{\hat{\mathbb P}}X_{\tau_1}\cdot e_1,\label{LVe23}\\
\exp{(-e^{-cL})}E_{\hat{\mathbb P}}\tau_1
&\le \varliminf_{n\to\infty}\frac{\tau_n}{n}
\le \varlimsup_{n\to\infty}\frac{\tau_n}{n}
\le \exp{(e^{-cL})}E_{\hat{\mathbb P}}\tau_1. \label{LVe24}
\end{align}
Note that (\ref{LVe23}), (\ref{LVe24}) hold even if 
$E_{\hat{\mathbb P}}X_{\tau_1}\cdot e_1=\infty$ or $E_{\hat{\mathbb P}}\tau_1=\infty$.
But it will be shown later that under our assumption, both of them are finite.

We claim that 
\begin{equation}\label{LVe5}
E_{\hat{\mathbb P}}X_{\tau_1}\cdot e_1<\infty.
\end{equation}
To see this, let $\Theta:=\{i: X_{\tau_k}\cdot e_1=i \text{ for some }k\in\mathbb{N}\}$.
Since $\tau_i$'s are finite on $A_{e_1}$,
there exist (recall that $\tau_0=0$) a sequence $(k_n)_{n\in\mathbb{N}}$ such that
$X_{\tau_{k_n}}\cdot e_1\le n<X_{\tau_{k_n+1}}\cdot e_1$ for all $n\in\mathbb{N}$ and 
$\lim_{n\to\infty}k_n=\infty$.
Hence for $n\ge 1$,
\[
\frac{\sum_{i=1}^n 1_{i\in \Theta}}{n}\le 
\frac{k_n+1}{X_{\tau_{k_n}}\cdot e_1}, \quad\text{ $\hat{\mathbb P}$-a.s..}
\]
Then, $\hat{\mathbb P}$-a.s.,
\[
\varlimsup_{n\to\infty}
\frac{\sum_{i=1}^n 1_{i\in \Theta}}{n}
\le
\varlimsup_{n\to\infty}\frac{n}{X_{\tau_n}\cdot e_1}.
\]

Let $B_k=\{\epsilon_{k,X_k}=0, X_{k+1}-X_k=e_1, \epsilon_{k+i,X_{k+i}}=1,X_{k+i+1}-X_{k+i}=e_1
 \text{ for all }1\le i\le L\}$. Then
\[
 \bar{P}_\omega(B_k)
 \ge 
 (d\kappa)^L(1-d\kappa)(\frac{\kappa}{2})(\frac{1}{2d})^L
 \stackrel{1\ge 2d\kappa}{>}(\frac{\kappa}{2})^{L+2}.
\]
Observe that by the definition of the regeneration times, for $n> L+1$,
\begin{align*}
&\{T_{n-L-1}=k,X_k= x-(L+1)e_1, R>k\}\cap B_k\cap\{R\circ\theta_{k+L+1}=\infty\}\\
&\subset\{R=\infty, n\in\Theta, T_n=k+L+1,X_{T_n}=x\}.
\end{align*}
Hence for $n> L+1$,
\begin{align*}
& \hat{\mathbb P}(n\in \Theta)\\
&\ge  
\sum_{k\in\mathbb{N},x\in\mathcal{H}_n}
\hat{\mathbb P}(B_k\cap\{T_{n-L-1}=k,X_k= x-(L+1)e_1,R\circ\theta_{k+L+1}=\infty\})\\
&\ge  \sum_{k\in\mathbb{N},x\in\mathcal{H}_n}
E_P \big[P_\omega\big(T_{n-L-1}=k,X_k= x-(L+1)e_1,R>k\big)(\frac{\kappa}{2})^{L+2}\\
&\qquad\qquad\qquad\qquad\qquad\qquad\qquad\quad\times 
P_\omega^x(R=\infty)\big]/\mathbb{P}(R=\infty).
\end{align*}
Since by (\ref{LVprop1}) and the translation invariance of $P$,
\[
E_P\big[P_\omega^x(R=\infty)|\omega_y:y\cdot e_1\le x\cdot e_1-L-1\big]
\ge 
\exp(-e^{-cL})\mathbb{P}(R=\infty),
\]
we have for $n>L+1$,
\begin{align}\label{LV2e10}
& \hat{\mathbb P}(n\in \Theta)\nonumber\\
&\ge
(\frac{\kappa}{2})^{L+2}\exp(-e^{-cL})
\sum_{k\in\mathbb{N},x\in\mathcal{H}_n}\mathbb{P}(T_{n-L-1}=k,X_k= x-(L+1)e_1,R>k)\nonumber\\
&\ge (\frac{\kappa}{2})^{L+2} e^{-1}\mathbb{P}(R=\infty).
\end{align}
Hence
\begin{align*}
\frac{C}{E_{\hat{\mathbb P}}X_{\tau_1}\cdot e_1}
\stackrel{(\ref{LVe23})}{\ge} E_{\hat{\mathbb P}}\varlimsup_{n\to\infty}\frac{n}{X_{\tau_n}\cdot e_1}
&\ge  E_{\hat{\mathbb P}}\varlimsup_{n\to\infty}
\frac{\sum_{i=1}^n 1_{i\in \Theta}}{n}\\
&\ge \varlimsup_{n\to\infty}E_{\hat{\mathbb P}}\frac{\sum_{i=1}^n 1_{i\in \Theta}}{n}\\
&\stackrel{(\ref{LV2e10})}{\ge} (\frac{\kappa}{2})^{L+2} e^{-1}\mathbb{P}(R=\infty)>0.
\end{align*}
This gives (\ref{LVe5}).

Now we can prove the theorem.
By (\ref{LVe23}) and (\ref{LVe24}), 
\begin{align}\label{LV*9}
\exp{(-2e^{-cL})}\frac{E_{\hat{\mathbb P}}X_{\tau_1}\cdot e_1}{E_{\hat{\mathbb P}}\tau_1}
&\le \varliminf_{n\to\infty}\frac{X_{\tau_n}\cdot e_1}{\tau_{n+1}}\nonumber\\
&\le \varlimsup_{n\to\infty}\frac{X_{\tau_{n+1}}\cdot e_1}{\tau_n}
\le \exp{(2e^{-cL})}\frac{E_{\hat{\mathbb P}}X_{\tau_1}\cdot e_1}{E_{\hat{\mathbb P}}\tau_1},
\end{align}
$\mathbb{P}(\cdot|A_{e_1})$-almost surely.
Further, by the fact that $|X_i|\le i$ and the obvious inequalities
\begin{equation*}
\varliminf_{n\to\infty}\frac{X_{\tau_n}\cdot e_1}{\tau_{n+1}}
\le 
\varliminf_{n\to\infty}\frac{X_n\cdot e_1}{n}
\le 
\varlimsup_{n\to\infty}\frac{X_n\cdot e_1}{n}
\le 
\varlimsup_{n\to\infty}\frac{X_{\tau_{n+1}}\cdot e_1}{\tau_n},
\end{equation*}
we have that
\begin{equation*}
\varlimsup_{n\to\infty}
\Bigl\lvert \frac{X_n\cdot e_1}{n}-
\frac{E_{\hat{\mathbb P}}X_{\tau_1}\cdot e_1}{E_{\hat{\mathbb P}}\tau_1}\Bigr\rvert
\le \exp{(2e^{-cL})}-1, \text{ $\mathbb{P}(\cdot|A_{e_1})$-a.s.}
\end{equation*}
Therefore, $\mathbb{P}(\cdot|A_{e_1})$-almost surely,
\[
\lim_{n\to\infty} \frac{X_n\cdot e_1}{n}
=
\lim_{L\to\infty}
\frac{E_{\hat{\mathbb P}}X_{\tau_1^{(L)}}\cdot e_1}{E_{\hat{\mathbb P}}\tau_1^{(L)}}:=v_{e_1},
\]
where $\tau_1$ is written as $\tau_1^{(L)}$ to indicate that it is an $L$-regeneration time.
Moreover, our assumption $\mathbb{P}(\varlimsup_{n\to\infty}X_n\cdot e_1/n>0)>0$ implies that
$v_{e_1}>0$ and (by (\ref{LV*9}))
\[
E_{\hat{\mathbb P}}\tau_1<\infty.
\]
Our proof is complete.\qed\\

If $v_{e_1}>0$, then it follows by (\ref{LVe24}) that
\begin{equation}\label{LVetau}
E_{\hat {\mathbb P}}\tau_n\le CnE_{\hat{\mathbb P}}\tau_1<\infty.
\end{equation}

Observe that although Theorem \ref{LVlln} is stated for $e_1$, the previous arguments, if properly modified, still work if one replaces
$e_1$ with any $z\in\mathbb{R}^d\setminus\{o\}$. So Theorem \ref{LVlln} is true for the general case. That is, for any $z\neq o$, there exist two deterministic constants $v_z, v_{-z}\ge 0$ such that
\[
\lim_{n\to\infty}\frac{X_n\cdot z}{n}=v_z1_{A_z}-v_{-z}1_{A_{-z}} 
\]
and that $\mathbb{P}(A_z\cup A_{-z})=1$ if $v_z>0$.
Then, by the same argument as in \cite[page 1112]{Go}, one concludes that the
limiting velocity $\lim_{n\to\infty}X_n/n$ can take at most two antipodal values.
This proves the first part of Theorem \ref{LVthm2}.

\section{Heat kernel estimates}\label{sechke}
The following heat kernel estimates are crucial for the proof of the uniqueness of the
non-zero velocity in the next section. Although in the mixing case we don't have iid
regeneration slabs, we know that (by Lemma \ref{LVl4}) a regeneration slab has little dependence on its remote
past. This allows us to use coupling techniques to get the same heat
kernel estimates as in \cite{Be}:
\begin{theorem}\label{hke}
Assume $v_{e_1}>0$. For $x\in\mathbb{Z}^d$ and $n\in \mathbb{N}$,
we let 
\[Q(n,x):=\hat{\mathbb P}(x \text{ is visited in }[\tau_{n-1},\tau_n)).\]
Then for any $x\in\mathbb{Z}^d$ and $n\in \mathbb{N}$,
\begin{align}
&\hat{\mathbb P}(X_{\tau_n}=x)\le Cn^{-d/2},\label{LVehke}\\
&
\sum_{x\in\mathbb{Z}^d}Q(n,x)^2\le C(E_{\hat{\mathbb P}}\tau_1)^2 n^{-d/2}.\label{LVehke2}
\end{align}
\end{theorem}

By Lemma \ref{LVl4}, we have for $n\ge 2$ and $1\le k\le n-1$, $\hat{\mathbb P}$-almost surely,
\begin{align}\label{LVnew1}
\frac{h_{k+1}(\cdot|J_{n-1},\ldots,J_{n-k})}{h_k(\cdot|J_{n-1},\ldots,J_{n-k+1})}
&=
\frac{h_{k+1}(\cdot|J_{n-1},\ldots,J_{n-k})}{h_n(\cdot|J_{n-1},\ldots,J_1)}
\frac{h_n(\cdot|J_{n-1},\ldots,J_1)}{h_k(\cdot|J_{n-1},\ldots,J_{n-k+1})}\nonumber\\
&\ge\exp(-e^{-c(k+1)L}-e^{-ckL})\nonumber\\
&\ge 1-e^{-ckL}
\end{align}
for large $L$. Hence for $n\ge 2$ and $1\le k\le n-1$, we can define a (random) probability measure $\zeta_{n,k}^{J_{n-1},\ldots,J_{n-k}}$ that satisfies
\begin{align}\label{LVnew2}
\MoveEqLeft
h_{k+1}(\cdot|J_{n-1},\ldots,J_{n-k})\\
&=e^{-ckL}\zeta_{n,k}^{J_{n-1},\ldots,J_{n-k}}(\cdot)+(1-e^{-ckL})h_k(\cdot|J_{n-1},\ldots, J_{n-k+1}).\nonumber
\end{align}

To prove Theorem \ref{hke}, we will first construct a sequence of random variables $(\tilde J_i, i\in\mathbb{N})$ such that for any $n\in\mathbb{N}$, 
\begin{equation}\label{LVnew0}
(\tilde J_1,\ldots,\tilde J_n)\sim \hat{\mathbb P}(J_1\in\cdot,\ldots, J_n\in\cdot),
\end{equation}
where ``$X\sim\mu$" means ``$X$ is of law $\mu$". 
\subsection{Construction of the $\tilde J_i$'s}
Our construction consists of three steps:

\noindent{\it Step 1.}
We let $\tilde J_1, \tilde J_{2,1}, \tilde \Delta_{2,1}$ be independent random variables  such that 
\[
\tilde J_1\sim h_1(\cdot|\emptyset),\quad \tilde J_{2,1}\sim h_1(\cdot|\emptyset)
\] 
and $\tilde \Delta_{2,1}$ is Bernoulli with parameter $e^{-cL}$. Let $\tilde Z_{2,1}$ be independent of $\sigma(\tilde J_{2,1}, \tilde \Delta_{2,1})$ such that
\[
P(\tilde Z_{2,1}\in\cdot|\tilde J_1)=
\zeta_{2,1}^{\tilde J_1}(\cdot).
\]
Setting $\tilde J_{2}:=(1-\tilde \Delta_{2,1})\tilde J_{2,1}+\tilde\Delta_{2,1}\tilde Z_{2,1}$, by \eqref{LVnew2} we have
\[
(\tilde J_1, \tilde J_2)\sim
\hat{\mathbb P}(J_1\in\cdot,J_2\in\cdot).
\]
\noindent{\it Step 2.}
For $n\ge 3$, assume that we have constructed $\tilde J_1$ and $(\tilde J_{i,1}, \tilde\Delta_{i,j}, \tilde Z_{i,j}, 1\le j<i\le n-1)$ such that
\[
(\tilde J_1,\ldots, \tilde J_{n-1})
\sim
\hat{\mathbb P}(J_1\in\cdot,\ldots,J_{n-1}\in\cdot),
\]
where for $2\le j\le i<n$,
\[
\tilde J_{i,j}:=(1-\tilde\Delta_{i,j-1})\tilde J_{i,j-1}+\tilde\Delta_{i,j-1}\tilde Z_{i,j-1}
\]
and
\[
\tilde J_i:=
\tilde J_{i,i}.
\]
Then we define $\tilde J_{n,1}$ and $(\tilde\Delta_{n,k}, \tilde Z_{n,k}, 1\le k\le n-1)$ to be random variables such that conditioning on the values of $\tilde J_1$ and $(\tilde J_{i,1}, \tilde\Delta_{i,j}, \tilde Z_{i,j}, 1\le j<i\le n-1)$,
\begin{itemize}
\item $(\tilde J_{n,1}, \tilde\Delta_{n,k}, \tilde Z_{n,k}, 1\le k\le n-1)$ are conditionally independent;
\item The conditional distribution of $\tilde J_{n,1}$ is $h_1(\cdot|\emptyset)$;
\item For $1\le k\le n-1$, the conditional distributions of $\tilde Z_{n,k}$ and $\tilde\Delta_{n,k}$ are $\zeta_{n,k}^{\tilde J_{n-1},\ldots, \tilde J_{n-k}}(\cdot)$ and Bernoulli with parameter $e^{-ckL}$, respectively.
\end{itemize}
\noindent{\it Step 3.} For $2\le k\le n$, set
\begin{align*}
&\tilde J_{n,k}:=(1-\tilde\Delta_{n,k-1})\tilde J_{n,k-1}+\tilde\Delta_{n,k-1}\tilde Z_{n,k-1}\\
\text{and }&\tilde J_n:=\tilde J_{n,n}.
\end{align*}
Then (by \eqref{LVnew2}) almost surely,
\begin{equation}\label{LV3e3}
P(\tilde J_{n,k}\in\cdot|\tilde J_{n-1},\ldots,\tilde J_1)=h_k(\cdot|\tilde J_{n-1},\ldots,\tilde J_{n-k+1}).
\end{equation}
It follows immediately that 
\begin{equation*}
(\tilde J_1,\ldots,\tilde J_n)\sim \hat{\mathbb P}(J_1\in\cdot,\ldots, J_n\in\cdot).
\end{equation*}
Therefore, by induction, we have constructed $(\tilde J_i, i\in\mathbb{N})$ such that \eqref{LVnew0} holds for all $n\in\mathbb{N}$. 

In what follows,  with abuse of notation, we will identify $\tilde J_i$ with
$J_i$ and simply write $\tilde J_{i,j}, \tilde \Delta_{i,j}, \tilde Z_{i,j}$ as
$J_{i,j}, \Delta_{i,j}$ and $Z_{i,j}$, $1\le j<i$. We still use $\hat{\mathbb
P}$ to denote the law of the random variables in the enlarged probability space.
\begin{remark}
To summarize, we have introduced random variables $J_{i,j}, \Delta_{i,j},
Z_{i,j}$, $1\le j<i$ such that for any $n\ge 2$,
\begin{align*}
&J_{n,2}=(1-\Delta_{n,1})J_{n,1}+\Delta_{n,1}Z_{n,1},\\
&\ldots,\\
&J_{n,n-1}=(1-\Delta_{n,n-2})J_{n,n-2}+\Delta_{n,n-2}Z_{n,n-2},\\
&J_n=J_{n,n}=(1-\Delta_{n,n-1})J_{n,n-1}+\Delta_{n,n-1}Z_{n,n-1}.
\end{align*}
Intuitively, we flip a sequence of ``coins" $\Delta_{n,n-1},\ldots,\Delta_{n,1}$ to determine whether $J_1,\ldots,J_{n-1}$ are in the ``memory" of $J_n$. For instance, if 
\[
\Delta_{n,n-1}=\cdots=\Delta_{n,n-i}=0,
\]
 then $J_n=J_{n,n-i}$ doesn't ``remember" $J_1,\ldots, J_i$ (in the sense that
\[
\hat{\mathbb P}(J_{n,n-i}\in\cdot|J_{n-1},\ldots, J_1)
=
h_{n-i}(\cdot|J_{n-1},\ldots, J_{i+1}).
\]
See \eqref{LV3e3}).
\end{remark}
\subsection {Proof of Theorem \ref{hke}}
For $1<i\le n$, let $I_n(i)$ be the event that 
$\Delta_{i,i-1}=\ldots=\Delta_{i,1}=0$ and
$\Delta_{m,m-1}=\ldots=\Delta_{m,m-i}=0$ for all
$i<m\le n$.
Note that on $I_n(i)$, 
\begin{equation}\label{LVnew3}
J_i=J_{i,1}\text{ and }J_m=J_{m,m-i}
\text{ for all }i<m\le n.
\end{equation}
Setting
\[M_n:=\{1\le i\le n: I_n(i)\neq\emptyset\},\]
we have
\begin{lemma}\label{LVliid}
For $n\ge 2$, let $H$ be a nonempty subset of $\{2,\ldots, n\}$, and set
\[M_n:=\{1< i< n: I_n(i)\neq\emptyset\}.\]
 Conditioning on the
event $\{M_n=H\}$, the sequence $(J_i)_{i\in H}$ is iid and independent of
$(J_i)_{i\in \{1,\ldots,n\}\setminus H}$.
\end{lemma}
\noindent{\it Proof of Lemma \ref{LVliid}:}
From our construction it follows that for any $i>1$, $J_{i,1}$ is independent
of 
\[
\sigma(\Delta_{k,j}, 1\le j<k)
\vee
\sigma(J_l, 1\le l<i)
\vee
\sigma(J_{m,m-i},m>i).
\]
Hence by \eqref{LVnew3}, for any $i\in H$ and any appropriate measurable sets $(V_j)_{1\le j\le n}$,
\begin{align*}
&\hat{\mathbb P}(J_j\in V_j, 1\le j\le n|M_n=H)\\
&=\hat{\mathbb P}(J_{i,1}\in V_i)
\hat{\mathbb P}(J_j\in V_j, 1\le j\le n, j\neq i|M_n=H).
\end{align*}
By induction, we get
\begin{align*}
&\hat{\mathbb P}(J_j\in V_j, 1\le j\le n|M_n=H)\\
&=\prod_{i\in H}\hat{\mathbb P}(J_{i,1}\in V_i)
\hat{\mathbb P}(J_j\in V_j, 1\le j\le n, j\notin H|M_n=H).
\end{align*}
The lemma is proved.\qed

\noindent{\it Proof of Theorem \ref{hke}:}
By Lemma \ref{LVliid}, for $i\in H$ and all $j\in\{1,\ldots, d\}$,
\[
\hat{\mathbb P}\big(X_{\tau_i}-X_{\tau_{i-1}}=(L+1)e_1\pm e_j|M_n=H\big)
=
\hat{\mathbb P}(X_{\tau_1}=(L+1)e_1\pm e_j)>0,
\]
where the last inequality  is due to ellipticity.
Hence arguing as in \cite[pages 736, 737]{Be}, using Lemma~\ref{LVliid} and the heat kernel estimate for bounded iid random walks in $\mathbb{Z}^d$, we get that for any $x\in\mathbb{Z}^d$,
\[
\hat{\mathbb P}(\sum_{i\in H}X_{\tau_i}-X_{\tau_{i-1}}=x|M_n=H)
\le C|H|^{-d/2},
\]
where $|H|$ is the cardinality of $H$.
Hence, for any subset $H\subset\{2,\ldots,n\}$ such that $|H|\ge n/2$,
\begin{align}\label{LVe13}
& \hat{\mathbb P}(X_{\tau_n}=x|M_n=H)\nonumber\\
&=\sum_y \hat{\mathbb P}(\sum_{i\in H}X_{\tau_i}-X_{\tau_{i-1}}=x-y, 
\sum_{i\in \{1,\ldots,n\}\setminus H}X_{\tau_i}-X_{\tau_{i-1}}=y|M_n=H)\nonumber\\
&=\sum_y 
\bigg[\hat{\mathbb P}\big(\sum_{i\in H}X_{\tau_i}-X_{\tau_{i-1}}=x-y|M_n=H\big)\nonumber\\
&\qquad\qquad\qquad\qquad\qquad\times
\hat{\mathbb P}(\sum_{i\in \{1,\ldots,n\}\setminus H}X_{\tau_i}-X_{\tau_{i-1}}=y|M_n=H)\bigg]\nonumber\\
&\le  C n^{-d/2},
\end{align}
where we used Lemma~\ref{LVliid} in the second equality.

On the other hand, 
\begin{align*}
|M_n|
&\ge n-\sum_{i=2}^n \bigg(
1_{\Delta_{i,i-1}+\cdots+\Delta_{i,1}>0}+\sum_{m=i+1}^n
1_{\Delta_{m,m-1}+\cdots+\Delta_{m,m-i}>0}
\bigg)\\
&=n-\sum_ {i=2} ^n
1_{\Delta_{i,i-1}+\cdots+\Delta_{i,1}>0}-\sum_{m=2}^n\sum_{i=2}^{m-1}1_{\Delta_{
m,m-1}+\cdots+\Delta_{m,m-i}>0}\\
&\ge n-2\sum_{m=2}^n K_m,
\end{align*}
where $K_m:=\sup\{1\le j<m:\Delta_{m,j}=1\}$. Here we follow the convention
that $\sup\emptyset=0$.
Since $K_m$'s are independent, and for $m\ge 2$,
\begin{align*}
E e^{K_m} &= \sum_{j=0}^{m-1} e^j \hat{\mathbf P}(K_m=j)\\
&\le  \sum_{j=1}^{m-1} e^j \hat{\mathbf P}(\Delta_{m,j}=1)+1\\
&\le \sum_{j=1}^\infty e^j e^{-cjL}+1\to 1 \text{ as $L\to\infty$},
\end{align*}
we can take $L$ to be large enough such that $E e^{K_m}\le e^{1/8}$ for all
$m\ge 2$ and so
\begin{align}\label{LVe14}
\hat{\mathbb P}(|M_n|<n/2)
&\le \hat{\mathbb P}(K_2+\cdots+K_n>n/4)\nonumber\\
&\le e^{-n/4}E e^{K_2+\cdots +K_n}
\le e^{-n/8}.
\end{align}
By (\ref{LVe13}) and (\ref{LVe14}), inequality (\ref{LVehke}) follows
immediately. 

Furthermore, since
\begin{align*}
& Q(n,x)\\
&= \sum_y \hat{\mathbb P}(X_{\tau_{n-1}}=y)
\hat{\mathbb P}(x \text{ is visited in }[\tau_{n-1},\tau_n)|X_{\tau_{n-1}}=y)\\
&\stackrel{\text{Lemma }\ref{LVl4}}{\le} 
C\sum_y \hat{\mathbb P}(X_{\tau_{n-1}}=y)\hat{\mathbb P}((x-y)\text{ is visited during }[0,\tau_1)),
\end{align*}
by H\"{o}lder's inequality we have
\begin{align*}
& Q(n,x)^2\\
&\le  C
\big[\sum_y\hat{\mathbb P}\big((x-y)\text{ is visited during }[0,\tau_1)\big)\big]\\
&\qquad\qquad\times\big[\sum_y\hat{\mathbb P}(X_{\tau_{n-1}}=y)^2
\hat{\mathbb P}\big((x-y)\text{ is visited during }[0,\tau_1)\big)\big]\\
&\le CE_{\hat{\mathbb P}}\tau_1 \sum_y \hat{\mathbb P}(X_{\tau_{n-1}}=y)^2
\hat{\mathbb P}\big((x-y)\text{ is visited during }[0,\tau_1)\big).
\end{align*}
Hence
\begin{align*}
& \sum_x Q(n,x)^2\\
&\le  CE_{\hat{\mathbb P}}\tau_1 
\sum_y \big[\hat{\mathbb P}(X_{\tau_{n-1}}=y)^2
\sum_x\hat{\mathbb P}\big((x-y)\text{ is visited during }[0,\tau_1)\big)\big]\\
&\le C(E_{\hat{\mathbb P}}\tau_1)^2 \sum_y \hat{\mathbb P}(X_{\tau_{n-1}}=y)^2\\
&\stackrel{(\ref{LVehke})}{\le}  C(E_{\hat{\mathbb P}}\tau_1)^2 n^{-d/2}\sum_y \hat{\mathbb P}(X_{\tau_{n-1}}=y)
=C(E_{\hat{\mathbb P}}\tau_1)^2 n^{-d/2}.
\end{align*}
Theorem \ref{hke} is proved.\qed

\section{The uniqueness of the non-zero velocity}\label{secunique}
In this section we will show that in high dimensions ($d\ge 5$), there
exists at most one non-zero velocity. The idea is the following.
Consider two random walk paths: one starts at the origin, the other starts near 
the $n$-th regeneration position of the first path. By Levy's martingale convergence
theorem, the second path is ``more and more transient" as $n$ grows (Lemma \ref{LVl1}).
On the other hand, by heat kernel estimates, when $d\ge 5$, two ballistic walks in opposite directions 
will grow further and further apart from each other (see Lemma \ref{LVl3}), thus they are almost independent.
This contradicts the previous fact that starting at the $n$-th regeneration point of the first path will prevent the second path from being transient in the opposite direction.\\

Set
$\delta=\delta(d):=\frac{d-4}{8(d-1)}$. (The reason for choosing this constant
will become clear in (\ref{LV2e6}).).
For any finite path $y_\cdot=(y_i)_{i=0}^M, M<\infty$, define $A(y_\cdot, z)$ to be the
set of paths $(x_i)_{i=0}^N, N\le\infty$ that satisfy
\begin{itemize}
\item[1)] $x_0=y_0+z$;
\item[2)] $d(x_i,y_j)>(i\vee j)^\delta$ if $i\vee j>|z|/3$.
\end{itemize}

The motivation for the definition of $A(y_\cdot, z)$ is as follows.
Note that for two paths $x_\cdot=(x_i)_{i=0}^N$ and $y_\cdot=(y_i)_{i=0}^M$ with $x_0=y_0+z$, if $i\vee j\le |z|/3$, then
\[
d(x_i,y_j)
\ge d(x_0,y_0)-d(x_0,x_i)-d(y_0,y_j)
\ge |z|-i-j
\ge |z|/3.
\]
Hence, for $(x_i)_{i=0}^N\in A(y_\cdot, z)$,
\begin{align}\label{LVe10}
\sum_{i\le N,j\le M}e^{-\gamma d(x_i,y_j)}
&\le 
\sum_{0\le i,j\le |z|/3}e^{-\gamma |z|/3}+\sum_{i\vee j>|z|/3}e^{-\gamma(i^\delta+j^\delta)/2}
\nonumber\\
&\le 
(\frac{|z|}{3})^2 e^{-\gamma |z|/3}+(\sum_{i=0}^\infty e^{-\gamma i^\delta/2})^2<C.
\end{align}
This gives us (by (G)) an estimate of the interdependence between
$\sigma\big(\omega_x: x\in (x_i)_{i=0}^N\big)$ and $\sigma\big(\omega_x: x\in (y_i)_{i=0}^M\big)$.

In what follows, we use 
\[
\tau'_\cdot=\tau_\cdot(-e_1,\epsilon,X_\cdot)
\]
 to denote the regeneration times in the $-e_1$ direction.
Assume that there are two opposite nonzero limiting velocities in directions $e_1$ and $-e_1$, i.e.,
\[
v_{e_1} \cdot v_{-e_1}>0.
\]
We let $\check{\mathbb P}(\cdot):=\mathbb{P}(\cdot|R_{-e_1}=\infty)$.
\begin{lemma}\label{LVl3}
Assume that there are two nonzero limiting velocities in direction $e_1$. We sample $(\epsilon,\tilde{X}_\cdot)$ according to $\hat{\mathbb P}$ and let
$\tilde{\tau}_\cdot=\tau_\cdot (e_1,\epsilon,\tilde{X}_\cdot)$ denote its regeneration times. 
For $n\ge 1$, we let 
\[
Y_\cdot^n=(Y_i^n)_{i=0}^{\tilde{\tau}_n}:=(\tilde{X}_{\tilde{\tau}_n-i})_{i=0}^{\tilde{\tau}_n}
\]
be the reversed path of $(\tilde{X}_i)_{i=0}^{\tilde{\tau}_n}$.
If $|z|$ is large enough, $d\ge 5$ and $n\ge 1$, then
\begin{equation}\label{LVe8}
E_{\hat{\mathbb P}}\check{\mathbb P}^{\tilde{X}_{\tilde{\tau}_n}+z}
\big(X_\cdot\in A(Y_\cdot^n, z)\big)>C>0.
\end{equation}
\end{lemma}

\pf
Let
\[m_z:=\lfloor|z|^{1/2}\rfloor.\]
Then
\begin{align}
&E_{\hat{\mathbb P}}\check{\mathbb P}^{\tilde{X}_{\tilde{\tau}_n}+z}
\big(X_\cdot\notin A(Y_\cdot^n, z)\big)\nonumber\\
&\le 
E_{\hat{\mathbb P}}\check{\mathbb P}^{\tilde{X}_{\tilde{\tau}_n}+z}(\tau'_{m_z}\ge |z|/3)+
\hat{\mathbb{P}}(\tilde{\tau}_n-\tilde{\tau}_{n-m_z}\ge |z|/3)\label{LV*12}\\
&\quad +E_{\hat{\mathbb P}}\check{\mathbb P}^{\tilde{X}_{\tilde{\tau}_n}+z}
(d(X_i,Y_\cdot^n)\le i^\delta\text{ for some } i>\tau'_{m_z})\label{LV*13}\\
&\quad +E_{\hat{\mathbb P}}\check{\mathbb P}^{\tilde{X}_{\tilde{\tau}_n}+z}
(d(\tilde{X}_{\tilde{\tau}_n-j},X_\cdot)\le j^\delta\text{ for some }j>\tilde{\tau}_n-\tilde{\tau}_{n-m_z}).\label{LV*14}
\end{align}

We will first estimate (\ref{LV*12}). By the translation invariance of the environment measure, 
\[
\check{\mathbb P}^x(\tau'_{m_z}\ge |z|/3)=\check{\mathbb P}(\tau'_{m_z}\ge |z|/3)
\text{ for any }x\in\mathbb{Z}^d.
\]
Hence
\begin{equation}\label{LV2e4}
E_{\hat{\mathbb P}}\check{\mathbb P}^{\tilde{X}_{\tilde{\tau}_n}+z}(\tau'_{m_z}\ge |z|/3)
=
\check{\mathbb P}(\tau'_{m_z}\ge |z|/3) 
\le 
\frac{3E_{\check{\mathbb P}}\tau'_{m_z}}{|z|}
\stackrel{(\ref{LVetau})}{\le}
C(E_{\check{\mathbb P}}\tau'_1)|z|^{-1/2}.
\end{equation}
Similarly,
\begin{equation}\label{LV2e5}
\hat{\mathbb P}(\tilde{\tau}_n-\tilde{\tau}_{n-m_z}\ge |z|/3)
\stackrel{(\ref{LVprop2})}{\le}
\exp{(e^{-cL})}\hat{\mathbb P}(\tau_{m_z}\ge |z|/3)\le C(E_{\hat{\mathbb P}}\tau_1) |z|^{-1/2}.
\end{equation}

To estimate (\ref{LV*13}) and (\ref{LV*14}), for $i\ge 1, n\ge j\ge 1$, we let
\begin{align*}
&Q'(i,x)=\check{\mathbb P}(x\text{ is visited in }[\tau'_{i-1},\tau'_i)),\\
&\tilde{Q}(j,x)=\hat{\mathbb P}(X_{\tau_n}+x \text{ is visited in}[\tau_{n-j},\tau_{n-j+1})).
\end{align*}
Note that by arguments that are similar to the proof of Theorem \ref{hke}, 
one can also obtain the heat kernel estimate
(\ref{LVehke}) for $Q'(i,x)$ and $\tilde{Q}(j,x)$.
For $l>0$, let $B(o,l)=\{x\in\mathbb{Z}^d: d(o,x)\le l\}$.
Recall the definition of the $r$-boundary in Definition \ref{LVdef1}.
By the translation invariance of the environment measure,
\[
\check{\mathbb P}^y(X_i=y+z)
=\check{\mathbb P}(X_i=z) \text{ for any }y,z\in\mathbb{Z}^d \text{and }i\in\mathbb{N}.
\]
Hence
\begin{align*}
&E_{\hat{\mathbb P}}\check{\mathbb P}^{\tilde{X}_{\tilde{\tau}_n}+z}
(d(X_i,\tilde{X}_\cdot)\le i^\delta\text{ for some } i>\tau'_{m_z})\\
&\le 
\sum_{i\ge m_z}\sum_{y\in\partial_1 B(o,i^\delta)}\sum_x
E_{\hat{\mathbb P}}\big[\check{\mathbb P}^{\tilde{X}_{\tilde{\tau}_n}+z}
(\tilde{X}_{\tilde{\tau}_n}+z+x\text{ is visited in }[\tau'_i,\tau'_{i+1}))\\
&\qquad\qquad\qquad\qquad\qquad\times 1_{\tilde{X}_{\tilde{\tau}_n}+z+x+y\in Y_\cdot^n}\big]\\
&=
\sum_{i\ge m_z}\sum_{y\in\partial_1 B(o,i^\delta)}\sum_x
\check{\mathbb P}(x\text{ is visited in }[\tau'_i,\tau'_{i+1}))
\hat{\mathbb P}(\tilde{X}_{\tilde{\tau}_n}+z+x+y\in Y_\cdot^n)\\
&=\sum_{i\ge m_z}\sum_{y\in\partial_1 B(o,i^\delta)}\sum_{j\le n}\sum_x
Q'(i,x)\tilde{Q}(j,x+z+y).
\end{align*}
By the heat kernel estimates and H\"{o}lder's inequality, 
\begin{align*}
\sum_{j\le n}\sum_x Q'(i,x)\tilde{Q}(j,x+z+y)
&\le 
\sqrt{\sum_x Q'(i,x)^2}\sum_{j\le n}\sqrt{\sum_x \tilde{Q}(j,x+y)^2}\\
&\le
 C(E_{\check{\mathbb P}}\tau'_1)i^{-d/4}
 \sum_{j\le n}(E_{\hat{\mathbb P}}\tau_1)j^{-d/4}\\
&\stackrel{d\ge 5}{\le} 
Ci^{-d/4}E_{\check{\mathbb P}}\tau'_1 E_{\hat{\mathbb P}}\tau_1.
\end{align*}
Thus
\begin{align}\label{LV2e6}
&E_{\hat{\mathbb P}}\check{\mathbb P}^{\tilde{X}_{\tilde{\tau}_n}+z}
(d(X_i,\tilde{X}_\cdot)\le i^\delta\text{ for some } i>\tau'_{m_z})\nonumber\\
&\le 
C\sum_{i\ge m_z}\sum_{y\in\partial_1 B(o,i^\delta)}i^{-d/4}
E_{\check{\mathbb P}}\tau'_1 E_{\hat{\mathbb P}}\tau_1\nonumber\\
&\le 
C\sum_{i\ge m_z}i^{(d-1)\delta}i^{-d/4}E_{\check{\mathbb P}}\tau'_1 E_{\hat{\mathbb P}}\tau_1
\le C 
|z|^{-(d-4)/8}
E_{\check{\mathbb P}}\tau'_1 E_{\hat{\mathbb P}}\tau_1,
\end{align}
where we used $d\ge 5$ and $\delta=\frac{d-4}{8(d-1)}$ in the last inequality.
Similarly, we have
\begin{align}\label{LV2e7}
&E_{\hat{\mathbb P}}\check{\mathbb P}^{\tilde{X}_{\tilde{\tau}_n}+z}
(d(\tilde{X}_{\tilde{\tau}_n-j},X_\cdot)
\le j^\delta\text{ for some }j>\tilde{\tau}_n-\tilde{\tau}_{n-m_z})\nonumber\\
&\le 
C |z|^{-(d-4)/8}E_{\check{\mathbb P}}\tau'_1 E_{\hat{\mathbb P}}\tau_1.
\end{align}

Combining (\ref{LV2e4}), (\ref{LV2e5}), (\ref{LV2e6}) and (\ref{LV2e7}), we conclude that
\begin{equation*}
E_{\hat{\mathbb P}}\check{\mathbb P}^{\tilde{X}_{\tilde{\tau}_n}+z}
\big(X_\cdot\in A(Y_\cdot^n, z)\big)>C>0,
\end{equation*}
if $|z|$ is large enough and $d\ge 5$.\qed\\

Let
\[
T^o=\inf\{i\ge 0: X_i\cdot e_1<0\}.
\]
For every fixed $\omega\in\Omega$ and $P_{\omega,\epsilon}^o$-almost every
$X_\cdot$,
\[
P_{\omega,\theta^n\epsilon}^{X_n}(T^o=\infty)1_{T^o>n}=P_{\omega,\epsilon}^o (T^o=\infty|X_1,\ldots,X_n),
\]
and so by Levy's martingale convergence theorem,
\[
\lim_{n\to\infty}P_{\omega,\theta^n\epsilon}^{X_n}(T^o=\infty)1_{T^o> n}= 1_{T^o=\infty}, 
\quad\text{$P_{\omega,\epsilon}^o$-almost surely}.
\]
Hence, for $(\omega, \epsilon,\tilde{X}_\cdot)$ sampled according to
$\hat{\mathbb P}$,
\[
\lim_{n\to\infty}
P_{\omega,\theta^{\tilde{\tau}_n}\epsilon}^{\tilde{X}_{\tilde{\tau}_n}}(T^o=\infty)=1, \quad\text{$\hat{\mathbb P}$-almost surely}.
\]
It then follows by the dominated convergence theorem that
\begin{equation}\label{LV3e7}
\lim_{n\to\infty}
E_{\hat{\mathbb P}}
P_{\omega,\theta^{\tilde{\tau}_n}\epsilon}^{\tilde{X}_{\tilde{\tau}_n}}(T^o<\infty)
=0.
\end{equation}

\begin{lemma}\label{LVl1}
For any $z\in\mathbb{Z}^d$,
\begin{equation}\label{LVe2}
\lim_{n\to\infty}
E_{\hat{\mathbb P}}
P_{\omega,\theta^{\tilde{\tau}_n}\epsilon}^{\tilde{X}_{\tilde{\tau}_n}+z}(T^o<\infty)
=0.
\end{equation}
\end{lemma}

\pf
For $n>|z|$, obviously
\[(\tilde{X}_{\tilde{\tau}_n}+z)\cdot e_1>0.\]
This together with ellipticity yields
\[
P_{\omega,\theta^{\tilde{\tau}_n}\epsilon}^{\tilde{X}_{\tilde{\tau}_n}}(T^o<\infty)
\ge{(\frac{\kappa}{2})^{|z|}}
P_{\omega,\theta^{\tilde{\tau}_n+|z|}\epsilon}^{\tilde{X}_{\tilde{\tau}_n}+z}(T^o<\infty).
\]
Hence using \eqref{LV3e7},
\[
\lim_{n\to\infty}
E_{\hat{\mathbb P}}
P_{\omega,\theta^{\tilde{\tau}_n+|z|}\epsilon}^{\tilde{X}_{\tilde{\tau}_n}+z}(T^o<\infty)
=0.
\]
On the other hand, noting that $\{R>\tau_1\}=\{R=\infty\}$,
\begin{align*}
&E_{\hat{\mathbb P}}
P_{\omega,\theta^{\tilde{\tau}_n+|z|}\epsilon}^{\tilde{X}_{\tilde{\tau}_n}+z}(T^o<\infty)\\
&= \sum_{m,x}
E_{P\otimes Q}[P_{\omega,\theta^{m+|z|}\epsilon}^{x+z}(T^o<\infty)
P_{\omega, \epsilon}^o (R>\tau_1,\tau_n=m,X_m=x)]/\mathbb{P}(R=\infty)\\
&=\sum_{m,x}
E_{P\otimes Q}[P_{\omega,\theta^{m}\epsilon}^{x+z}(T^o<\infty)
P_{\omega, \epsilon}^o (R>\tau_1,\tau_n=m,X_m=x)]/\mathbb{P}(R=\infty)\\
&=
E_{\hat{\mathbb P}}P_{\omega,\theta^{\tilde{\tau}_n}\epsilon}^{\tilde{X}_{\tilde{\tau}_n}+z}(T^o<\infty),
\end{align*}
where we used the independence (under $Q$) of $P_{\omega,
\theta^m\epsilon}^{x+z}(T^o<\infty)$ and
$P_{\omega,\epsilon}^o(R>\tau_1,\tau_n=m, X_m=x)$ in the second to last
equality. The conclusion follows.
\qed\\

\noindent\textit{Proof of the uniqueness of the non-zero velocity when $d\ge 5$, as stated in Theorem \ref{LVthm2}:} 
If the two antipodal velocities are both non-zero, we assume that 
\[v_{e_1}\cdot v_{-e_1}>0.\]

Sample $(\omega,\epsilon_\cdot,\tilde{X}_\cdot)$ according
to $\hat{\mathbb P}$.
Henceforth, we take $z=z_0$ such that (\ref{LVe8}) holds
and 
\[z_0\cdot e_1<-L.\] 

We will prove Theorem \ref{LVthm2} by showing that
\begin{equation}\label{LVcontradiction}
E_{\hat{\mathbb P}} 
P_{\omega,\theta^{\tilde{\tau}_n}\epsilon}^{\tilde{X}_{\tilde{\tau}_n}+z_0}(T^o<\infty)>C
\end{equation}
for all $n>|z_0|$, which contradicts with (\ref{LVe2}).

First, let $\mathcal{G}$ denote the set of finite paths $y_\cdot=(y_i)_{i=0}^M$ that satisfy $y_M=0, M<\infty$.
Then
\begin{align}\label{LVe11}
& E_{\hat{\mathbb P}} 
P_{\omega,\theta^{\tilde{\tau}_n}\epsilon}^{\tilde{X}_{\tilde{\tau}_n}+z_0}(T^o<\infty)\\
&\ge 
E_{\hat{\mathbb P}} 
P_{\omega,\theta^{\tilde{\tau}_n}\epsilon}^{\tilde{X}_{\tilde{\tau}_n}+z_0}
\big((X_i)_{i=0}^{T^o}\in A(Y_\cdot^n, z_0),T^o<\infty\big)\nonumber\\
&=
\sum_{y_\cdot=(y_i)_{i=0}^M\in \mathcal{G}}
E_{\hat{\mathbb P}} 
[P_{\omega,\theta^M\epsilon}^{y_0+z_0}
\big((X_i)_{i=0}^{T^o}\in A(y_\cdot, z_0),T^o<\infty\big)
1_{Y_\cdot^n=y_\cdot}]\nonumber\\
&=
\frac{1}{\mathbb{P}(R=\infty)}
\sum_{y_\cdot\in \mathcal{G}}
\sum_{\substack{N<\infty\\(x_i)_{i=0}^N\in A(y_\cdot, z_0)}}
E_{P\otimes Q}
[P_{\omega,\theta^M\epsilon}^{y_0+z_0}
\big((X_i)_{i=0}^{T^o}=x_\cdot\big)
P_{\omega,\epsilon}(Y_\cdot^n=y_\cdot)].\nonumber
\end{align}
By the definition of the regeneration times, for any finite path $y_\cdot=(y_i)_{i=0}^M$,
there exists an event $G_{y_\cdot}$ such that $P_{\omega,\epsilon}(G_{y_\cdot})$ is
$\sigma(\epsilon_{i,y_i}, \omega_{y_j}:0\le i\le M, 0\le j\le M-L)$-measurable and
\[
\{Y_\cdot^n=y_\cdot\}=\{(\tilde{X}_i)_{i=0}^{\tilde{\tau}_n}=(y_{M-j})_{j=0}^M\}
=G_{y_\cdot}\cap\{R\circ\theta_M=\infty\}.
\]
Hence, for and any $y_\cdot=(y_i)_{i=0}^M\in \mathcal{G}$ and
$x_\cdot=(x_i)_{i=0}^N\in A(y_\cdot,z_0)$, $N<\infty$,
\begin{align}\label{LV2e8}
& E_{P\otimes Q}
[P_{\omega,\theta^M\epsilon}^{y_0+z_0}\big((X_i)_{i=0}^{T^o}=x_\cdot, R_{-e_1}>N\big)
P_{\omega,\epsilon}(Y_\cdot^n=y_\cdot)]\nonumber\\
&=E_P[
\bar{P}_\omega^{y_0+z_0}\big((X_i)_{i=0}^{T^o}=x_\cdot, R_{-e_1}>N\big)
\bar{P}_\omega(G_{y_\cdot})
\bar{P}_\omega^{y_0}(R=\infty)]\nonumber\\
&\stackrel{(\ref{LVprop1})}{\ge}
CE_P[
\bar{P}_\omega^{y_0+z_0}\big((X_i)_{i=0}^{T^o}=x_\cdot, R_{-e_1}>N\big)
\bar{P}_\omega(G_{y_\cdot})]
\bar{\mathbb P}(R=\infty).
\end{align}
where we used in the equality that $(\epsilon_{i,x})_{i\ge 0, x\in\mathbb{Z}^d}$
are iid and in the inequality the fact that
\[
\bar{P}_\omega^{y_0+z_0}\big((X_i)_{i=0}^{T^o}=x_\cdot, R_{-e_1}>N\big)
\bar{P}_\omega(G_{y_\cdot})
\]
is $\sigma(\omega_v: v\cdot e_1\le y_0\cdot e_1-L)$-measurable (note that $z_0\cdot e_1<-L$).
Further, by Lemma \ref{LVc2} and (\ref{LVe10}), we have
\begin{align}\label{LV2e9}
&E_P[
\bar{P}_\omega^{y_0+z_0}\big((X_i)_{i=0}^{T^o}=x_\cdot, R_{-e_1}>N\big)
\bar{P}_\omega(G_{y_\cdot})]\nonumber\\
&\ge 
C\bar{\mathbb P}^{y_0+z_0}\big((X_i)_{i=0}^{T^o}=x_\cdot, R_{-e_1}>N\big)
\bar{\mathbb P}(G_{y_\cdot}).
\end{align}
Note that
\begin{align}\label{LVe12}
\bar{\mathbb P}(G_{y_\cdot})\bar{\mathbb P}(R=\infty)
&\stackrel{(\ref{LVprop1})}{\ge}
CE_P
[\bar{P}_\omega (G_{y_\cdot})\bar{P}_{\omega}^{y_0}(R=\infty)]
\nonumber\\
&=C\bar{\mathbb P}(Y_\cdot^n=y_\cdot)
\ge
C\hat{\mathbb P}(Y_\cdot^n=y_\cdot).
\end{align}
Therefore, by (\ref{LVe11}), (\ref{LV2e8}) and (\ref{LV2e9}),
\begin{align*}
& E_{\hat{\mathbb P}} 
P_{\omega,\theta^{\tilde{\tau}_n}\epsilon}^{\tilde{X}_{\tilde{\tau}_n}+z_0}(T^o<\infty)\\
&\ge
C\sum_{y_\cdot\in \mathcal{G}}
\sum_{\substack{N<\infty\\(x_i)_{i=0}^N\in A(y_\cdot, z_0)}}
\bar{\mathbb P}^{y_0+z_0}\big((X_i)_{i=0}^{T^o}=x_\cdot, R_{-e_1}>N\big)
\bar{\mathbb P}(G_{y_\cdot})\bar{\mathbb P}(R=\infty)\\
&\stackrel{(\ref{LVe12})}{\ge}
C\sum_{y_\cdot\in \mathcal{G}}
\sum_{\substack{N<\infty\\(x_i)_{i=0}^N\in A(y_\cdot, z_0)}}
\bar{\mathbb P}^{y_0+z_0}\big((X_i)_{i=0}^{T^o}=x_\cdot, R_{-e_1}>N\big)
\hat{\mathbb P}(Y_\cdot^n=y_\cdot)\\
&\ge 
CE_{\hat{\mathbb P}}
\check{\mathbb P}^{\tilde{X}_{\tilde{\tau}_n}+z_0}
(X_\cdot\in A(Y_\cdot^n, z_0))\stackrel{\text{Lemma }\ref{LVl3}}{>}C.
\end{align*}
(\ref{LVcontradiction}) is proved.\qed

\chapter{Invariance Principle for Random Walks in Balanced Random Environment}
\label{CLT chapter}


This chapter is devoted to the proofs of Theorem~\ref{CLT1} and Theorem~\ref{CLT2}.

The structure of this chapter is as follows. In Section \ref{SePeEn} we construct the ``periodized environments" as in \cite{Sz1, ZO}, and show that the
proof of $Q\sim P$ can be reduced to the proof of the inequality (\ref{CPhi0}).
Using the  maximum principle (Theorem \ref{CMP}), we then prove (\ref{CPhi0})
in Section \ref{SeMP} under the assumptions of Theorem \ref{CLT1}(i). In Section
\ref{SePercE}, which is devoted to the iid setup,
we prove Theorem \ref{CLT2}(i) using percolation tools. Section \ref{SeTran} is
devoted
to the proof of the transience of the RWRE for $d\geq 3$, thus providing
a proof of Theorem \ref{CLT1}(ii). In Section \ref{SeTriid}, we will show
a modified maximum principle for balanced difference operators, and use it to prove
Theorem \ref{CLT2}(ii).

Throughout, $C$ denotes a generic positive constant that may depend on
dimension only, and whose value may change from line to line.
\section {The periodized environments}\label{SePeEn}
As in \cite{Sz1, ZO}, the following periodic structure of the environment
is introduced.

Let $ \Delta_N (x_0)=\{x\in \mathbb{Z}^{d}: |x-x_0|_{\infty}\le N\} $ be the cube centered at $x_0$ of length $2N$. Let $\Delta_N=\Delta_N(o)$. For any $ x\in\mathbb{Z}^{d} $, set
$$ \hat{x}:=x+(2N+1)\mathbb{Z}^{d}\in \mathbb{Z}^{d}/(2N+1)\mathbb{Z}^{d}. $$

For any fixed $ \omega \in \Omega $, we define $\omega^{N}$ by setting
$ \omega^{N}(x)=\omega(x) $ for $ x \in \Delta_{N} $
and $\omega^N (y)=\omega^N (x)$ for $y \in \mathbb{Z}^{d} $ whenever $\hat{y}=\hat{x}$. Let $ \Omega^{N}=\{\omega^{N}: \omega\in\Omega\}$. Let $ \{X_{n,N}\} $ denote the random walk on $ \mathbb{Z}^{d} $ in the environment $ \omega^{N} $. Then $ \{\hat{X}_{n,N}\} $ is an irreducible finite-state Markov chain,
hence it possesses a unique invariant probability measure, which can
always be written in
the form 
\[ \dfrac{1}{(2N+1)^d}\sum_{x\in\Delta_{N}}\Phi_{N}(x)\delta_{\hat{x}} . \]
Here $\Phi_N$ is some function on $\Delta_N$ and $(2N+1)^{-d}\Phi_{N}(\cdot)$ 
sums to $1$, so that $\Phi_N$ can be interpreted as a density with
respect to the uniform measure on $\Delta_N$.

Define 
\[ Q_{N}=Q_{N,\omega}=\dfrac{1}{(2N+1)^{d}}\sum_{x\in\Delta_{N}}\Phi_N (x)\delta_{\theta^{x}\omega^{N}} \]
as a probability measure on $\Omega^N$.
Then, for any $x\in \Delta_N$, 
\begin{align*}
\sum_{y\in\Delta_N} Q_N(\theta^y \omega^N)M(\theta^y \omega^N, \theta^x\omega^N)
&=\sum_{y\in\Delta_N} \frac{\Phi_N (y)}{(2N+1)^d}\omega^N(y,x)\\
&= \frac{\Phi_N (x)}{(2N+1)^d}=Q_N (\theta^x \omega^N).
\end{align*}
This implies that
$Q_N$ is the invariance probability measure
(with respect to the kernel $M$) for the Markov chain 
$\{\bar{\omega}^{N}(n)\}$ on $\Omega^{N}$.

We will show that $Q_N$ converges weakly to some measure $ Q $ with good
properties. To do this, we first introduce a sequence of measures
\[
P_{N}=P_{N,\omega}=\dfrac{1}{(2N+1)^{d}}
\sum_{x\in\Delta_{N}}\delta_{\theta^{x}\omega^{N}},
\]
which, by the multidimensional ergodic theorem (see Theorem (14.A8) in \cite{Ge}
and also Theorem 1.7.5 in \cite{Kr}), converges weakly to $P$, $P$-a.s.

Let $ \{\omega_{\gamma}^{N}\}_{\gamma=1}^{k} $ denote the set of distinct states in $ \{ \theta^{x}\omega^{N}\}_{x\in \Delta_{N}} $ and $ C_{N}(\gamma):=\{x\in \Delta_{N}: \theta^{x}\omega^{N}=\omega_{\gamma}^{N}\} $.
Set,
for any finite subset $E\subset\mathbb{Z}^d$,
\[
\lVert f\rVert_{E,j}:=(|E|^{-1}\sum_{x \in E} |f(x)|^{j})^{\frac{1}{j}}.
\]
Since
$ \ud Q_N/\ud
P_N=\sum_{\gamma=1}^{k}\delta_{\omega_{\gamma}^{N}}|C_{N}(\gamma)|^{-1}\sum_{
x\in C_{N}(\gamma)}\Phi_{N}(x):=f_{N}$, we have that, for any measurable
function $g$ on $\Omega$,
\begin{align}\label{Q_n0}
|Q_N g|
 &\le (\int  f_{N}^{\alpha} \ud P_N)^{\frac{1}{\alpha}}(\int  |g|^{\alpha'} \ud P_N)^{\frac{1}{\alpha'}} \nonumber\\
 &\le  \big( \frac{1}{|\Delta_N |}\sum_{\gamma =1}^{k}\sum_{x \in C_N (\gamma)}\Phi_N (x)^{\alpha} \big)^{\frac{1}{\alpha}}(\int  |g|^{\alpha'} \ud P_N)^{\frac{1}{\alpha'}}
 \nonumber\\
 &=  \lVert         \Phi _N \rVert _{\Delta_N, \alpha}(P_{N}|g|^{\alpha'})^{\frac{1}{\alpha'}},
\end{align}
where
$\alpha'$ is the H\"older conjugate of $\alpha$, $1/\alpha+1/\alpha'=1$,
and we used H\"older's inequality in the first and the second inequalities.
Since $\Omega$ is compact with respect to the product topology,
along some subsequence $N_k\to\infty$,
$\{Q_{N_k}\}$ converges weakly to a limit, denoted $Q$.
Assume for the moment that
\begin{equation}\label{CPhi0}
\varlimsup_{N\to\infty}\lVert         \Phi _N \rVert _{\Delta_N, \alpha}
\le C, \quad P\mbox{- }a.s.
\end{equation}
We then show that, for a.e. $\omega\in\Omega$,
\begin{equation}\label{CLTf2}
Q\ll P   .
\end{equation}
Indeed, let $A\subset \Omega$ be measurable.
Let $ \rho $ denote a metric on the Polish space $ \Omega $.
For any closed subset
$ F \subset A $, $ \delta >0 $,
introduce the function $f(\omega)=[1-\rho (\omega, F)/\delta ]^+$
which is supported on
$ F_{\delta}=\{\omega \in \Omega: \rho (\omega, F) <\delta\}$.
Then, by (\ref{Q_n0}), (\ref{CPhi0}),
\begin{equation*}
Q F\le \varlimsup_{N\to\infty}Q_N f \le C (P f^{\alpha'})^{\frac{1}{\alpha'}}\le  C (P F_{\delta})^{\frac{1}{\alpha'}} .
\end{equation*}
Letting  $ \delta\downarrow 0 $, we get $ Q F \le C (P F)^{\frac{1}{\alpha'}}$.
Taking supremums over all closed subset $ F \subset A $, one concludes that
$Q A \le C\cdot (P A)^{\frac{1}{\alpha'}}$, which proves (\ref{CLTf2}).

Once we have (\ref{CLTf2}), it is standard to check, using ellipticity,
that $ \bar{\omega}(n) $ is ergodic with respect to $ Q $ and $Q\sim P$
(see \cite{Sz1, ZO}). (Thus,
by the ergodic theorem, $Q$ is uniquely determined by $Qg=\lim_{n\to\infty} E\sum_{j=0}^{n-1}g(\bar{\omega}_j)/n$ for every bounded measurable $g$. Hence $Q$ is \textit{the} weak limit of $Q_N$.) Therefore, to prove the invariance principle it suffices to prove (\ref{CPhi0}).
Sections
\ref{SeMP} and Section \ref{SePercE} are devoted to the proof of (\ref{CPhi0}),
under the assumptions of Theorems \ref{CLT1} and \ref{CLT2}.

\section {Maximum principle and proof of Theorem \ref{CLT1}(i)}\label{SeMP}
Throughout this section, we fix an $\omega\in \Omega$.
For any bounded set $ E \subset \mathbb{Z}^{d} $, let $\partial E =\{y \in E^{c}: \exists x\in E, |x-y| _{\infty}=1\}$, $ \bar{E}=E \bigcup \partial E $ and $ \diam(E)=\max\{|x-y| _{\infty}: x, y\in E\} $. For any function $f$ defined on $ \bar{E}$ ,
let $ L_{\omega} $ denote the operator
\begin{equation}\label{operator}
 (L_{\omega}f )(x)=\sum_{i=1}^{d}\omega(x, e_i)[f(x+e_i)+f(x-e_i)-2f(x)], \quad
x\in E .
\end{equation}

The following discrete maximum principle is an adaption of Theorem 2.1 of
\cite{KT}.
\begin{theorem}[Maximum Principle]\label{CMP}
Let $E\subset \mathbb{Z}^d$ be bounded,
and  let $u$ be a function on $ \bar{E} $. For all
$x\in E$, assume
$ \varepsilon(x)>0 $ and define
$$ I_{u}(x):=\{s \in \mathbb{R}^{d}:  u(x)-s\cdot x \ge u(z)-s\cdot z,  \forall
z \in\bar{E}\}  .$$
If $ L_{\omega} u(x) \ge -g(x)$ for all $x \in E$ such that $I_u (x)\ne
\emptyset$,
then
\begin{equation}\label{Cmpremark}
 \max_{E} u \le
  C\diam{\bar{E}}
 \bigg(\sum_{\substack{
                           x\in E\\
                           I_u (x)\ne \emptyset}
                      }|\frac{g}{\varepsilon}|^d\bigg)^{\frac{1}{d}}
 +\max_{\partial E}u  .
\end{equation}
In particular,  \[ \max_{E} u \le
C\diam{\bar{E}}\cdot |E|^{\frac{1}{d}} \lVert        
\frac{g}{\varepsilon}\rVert _{E,d}+\max_{\partial E}u   .\]
\end{theorem}
\pf See the proof of Theorem 2.1 in \cite{KT}.\qed\\

Define the stopping times $ \tau_0=0 $, $ \tau_1 =\tau :=\min \{j \ge 1: |X_{j,N}-X_{0,N}|_{\infty}> N\} $ and
$ \tau_{j+1}=\min \{n>\tau_j : |X_{n,N} -X_{\tau_j, N}|_{\infty}> N\} $.

\begin{lemma}\label{CLTE}
Let $\omega^{N}$, $\{X_{n,N}\}$ be as in Section 1 and $\tau$ as defined above, then
there exists a constant $c$ such that, for all $N$ large,
 $$E_{\theta^{x}\omega^{N}}^{o} (1-\frac{c}{N^2})^{\tau}\leq C <1   .$$
\end{lemma}
\pf
Since P is balanced, $ X_{n,N} $ is a martingale
and it follows from Doob's inequality that for any $K \ge 1  $,
\begin{align*}
P_{\theta^{x}\omega^{N}}^{o} \{\tau \le K\}
& \le 
 2 \sum_{i=1}^{d} P_{\theta^{x}\omega^{N}}^{o}\{\sup_{n \le K} X_{n, N}(i) \ge N+1\}\\
& \le 
 \frac{2}{N+1} \sum_{i=1}^{d}E_{\theta^{x}\omega^{N}}^{o} X_{K, N} (i)^{+}
\le \frac{2d}{N+1} \sqrt{K},
\end{align*}
where $ X_{n,N} (i) $ is the $i$-th coordinate of $ X_{n,N} $. Hence
\begin{equation*}
E_{\theta^{x}\omega^{N}}^{o} (1-\frac{c}{N^2})^{\tau}
 \le  (1-\frac{c}{N^2})^{K} +\frac{2d}{N+1}\sqrt{K}   .
 \end{equation*}
Taking $c=16 d^2$ and $K=N^2/16d^2$, we get
$E_{\theta^{x}\omega^{N}}^{o} (1-\frac{c}{N^2})^{\tau}\le e^{-1}+2^{-1} .$ \qed

\begin{theorem}\label{CPhi}
\begin{equation}\label{CPhi1}
\lVert         \Phi_N \varepsilon\rVert _{\Delta_N , \beta}\le C,
\end{equation}
where $\beta=d'=d/(d-1)$.
\end{theorem}
\pf
 Let $c$ be the same
constant as in the previous lemma. For any function $ h\ge 0 $ on $\Delta_N $,
\begin{align*}
&\lVert         \Phi_N \cdot h\rVert _{\Delta_N , 1}\\
&=
 \frac{c}{N^2}\sum_{x \in \Delta_N}\frac{\Phi_N (x)}{|\Delta_N |}
\sum_{m\ge 0}E_{\omega^N}^{x}\sum_{\tau_m\le\ j <\tau_{m+1}}(1-\frac{c}{N^2})^{j}h(\hat{X}_{j,N})\\
&\le 
 \frac{c}{N^2}\sum_{x \in \Delta_N}\frac{\Phi_N (x)}{|\Delta_N |}\sum_{m\ge 0}E_{\omega^N}^{x}(1-\frac{c}{N^2})^{\tau_m}
E_{\omega^N}^{\hat{X}_{\tau_m , N}}\sum_{j=0}^{\tau -1}h(\hat{X}_{j,N})\\
&\le 
 \frac{c}{N^2}\sum_{x \in \Delta_N}\frac{\Phi_N (x)}{|\Delta_N |}\sum_{m\ge 0}\big[\sup_{y\in \Delta_N}E_{\omega^N}^{y}(1-\frac{c}{N^2})^{\tau}\big]^{m}\cdot \sup_{y\in \Delta_N}E_{\omega^N}^{y}
\sum_{j=0}^{\tau -1}h(\hat{X}_{j,N}).
 \end{align*}
Since the function $f(x)=E_{\omega^N}^{x} \sum_{j=0}^{\tau -1}h(\hat{X}_{j,N})$ satisfies
\begin{equation}
\left\{
\begin{array}{rl}
L_{\omega^N}f(x)=h(x), & \text{if } x\in\Delta_N\\
f(x)=0, & \text{if } x\in\partial \Delta_N,
\end{array}
\right.
\end{equation}
we can apply the maximum principle (Theorem \ref{CMP}) and get
$$\sup_{y\in \Delta_N}E_{\omega^N}^{y}\sum_{j=0}^{\tau -1}h(\hat{X}_{j,N}) \le C N^2
\lVert
\frac{h}{\varepsilon}
\rVert _{\Delta_N ,d}.$$
This, together with Lemma \ref{CLTE} and $\sum_{x \in \Delta_N}\Phi_N
(x)/|\Delta_N |=1$, yields
$$\lVert         \Phi_N \cdot h\rVert _{\Delta_N , 1}
\le C \lVert         \frac{h}{\varepsilon}\rVert _{\Delta_N ,d}.$$
Hence by the duality of norms,
\begin{equation*}
\lVert         \Phi_N \varepsilon\rVert _{\Delta_N , \beta}=\sup_{\norm{h/\varepsilon}_{\Delta_N, d}=1}\norm{
\Phi_N h} _{\Delta_N, 1}\le C .  \quad  \mbox{\qed}
\end{equation*}

\noindent{\it Proof of (\ref{CPhi0}) under the assumption of Theorem \ref{CLT1}(i) :}\\
Assume that
\begin{equation}\label{asm}
\mathrm{E}\varepsilon(o)^{-p}< \infty
\mbox{ for some } p>d.
\end{equation}
Take $\alpha=(1-1/d+1/p)^{-1}$. We use H\"older's inequality
and Theorem \ref{CPhi}
to get
\[
\lVert \Phi_N\rVert _{\Delta_N,\alpha}\le \lVert  \Phi_N \varepsilon\rVert _{\Delta_N , \beta}\lVert         \varepsilon^{-1}\rVert _{\Delta_N, p}
\le C \lVert \varepsilon^{-1}\rVert _{\Delta_N, p}  .
\]
By the multidimensional ergodic theorem,
\begin{equation*}
\lim_{N \to \infty}\lVert         \varepsilon^{-1}\rVert _{\Delta_N, p}=(E \varepsilon(o)^{-p})^{\frac{1}{p}}<\infty,
\quad  P\mbox{- }a.s. \qed
\end{equation*}

\begin{remark}
Without the assumption (\ref{asm}),  the conclusion
\eqref{CPhi0} may fail.
To see the difficulty,
let $$A=A(\omega, \varepsilon_0)=\{x:\min_i \omega(x, e_i)<\varepsilon_0\}.$$
By (\ref{CPhi1}) we have
$$\lVert         \Phi_N 1_{A^c}\rVert _{\Delta_N, \beta}\le \lVert         \Phi_N \frac{\varepsilon}{\varepsilon_0}\rVert _{\Delta_N, \beta} \le \frac{C}{\varepsilon_0}.$$
In order to proceed as before, we need to show that
$\varlimsup_{N\to\infty}\lVert \Phi_N 1_{A}\rVert _{\Delta_N, \alpha} \le C$
for some $1 <\alpha \le \beta$ .
As Bouchaud's trap model
\cite{Bou,BAC} shows,
this is not always the case. However, if $P\{\max_{|e|=1}\omega(o,e)\ge\xi_0\}=1$,
then for $x\in A$, we have, using that the environment is balanced,
some control of $\Phi_N (x)$
by $\Phi_N|_{A^c}$ (see Lemma \ref{CPhicontrol}).
Further,
in the iid case, $A$ corresponds to  a `site percolation' model,
whose cluster sizes can be estimated.
We will show in the next section that these properties lead
to a proof of (\ref{CPhi0}) in the iid setup, without moment assumptions.
\end{remark}

\section{A percolation estimate and proof of Theorem \ref{CLT2}(i)}\label{SePercE}
In this section we consider the RWRE in the iid setting where
$\max_{|e|=1}\omega(x,e)\ge \xi_0$  for all $x\in\mathbb{Z}^d$ and all $\omega\in\Omega$.
We begin by introducing some terminology.\\

The \textit{$l^1$-distance} (graph distance) from $x$ to $y$ is defined as
$$d(x,y)=|x-y|_1=\sum_{i=1}^{d}|x_i-y_i|.$$
Note that $|x|_{\infty}\le |x|_1 \le d|x|_{\infty}$.

In an environment $\omega$, we say that a site $x$  is \textit{open}(\textit{closed}) if $\min_i \omega(x, e_i)<\varepsilon_0 (\ge\varepsilon_0, resp.)$ and that an edge of $\mathbb{Z}^d$ is open if its endpoints are open. 
Here $\varepsilon_0>0$ is a constant whose value is to be determined.
An edge is called closed if it is not open.
Let $A=A(\omega)$ denote the subgraph of $\mathbb{Z}^d$ obtained
by deleting all  closed edges and
closed sites.
We call $A(\omega)$ a \textit{site percolation}
with parameter $p=p(\varepsilon_0)=P\{\min_i \omega(x, e_i)< \varepsilon_0\}$.

A \textit{percolation cluster} is a connected component of  $A$.
(Although here a percolation cluster is defined as a graph, we also use it as
a
synonym for its set of vertices.)
The $l^1$ diameter of a percolation cluster $B$ is defined as $l(B)=\sup_{x\in B, y\in \partial B}d(x,y)$.
For $x\in A$, let $A_x$ denote the percolation cluster that contains $x$ and let $l_x$ denote its diameter. Set $A_x=\emptyset$ and $l_x=0$ if $x\notin A$.
We let $\varepsilon_0$ be small enough such that $l_x<\infty$ for all $x\in \mathbb{Z}^d$.\\

We call a sequence of sites $(x^1, \cdots, x^n)$ a \textit{path}  from $x$ to $y$ if $x^1=x$, $x^n=y$ and
$|x^j-x^{j+1}|=1$
 for $j=1,\cdots, n-1$. Let
$$\square=\{(\kappa_1,\cdots, \kappa_d)\in \mathbb{Z}^d: \kappa_i=\pm 1\}.$$ 
We say that a path $\{x^1, \cdots, x^n\}$ is a \textit{$\kappa$-path}, $\kappa\in \square$,  if
$$\omega(x^j, x^{j+1}-x^j)\ge\xi_0$$ and
$\kappa_i(x^{j+1}-x^j)_i\ge 0$ for all $i=1,\cdots, d$ and $j=1,\cdots, n-1$.
Observing that for each site there exist at least two neighbors (in opposite directions) to whom the transition probabilities are $\ge \xi_0$, we have the following property concerning the structure of the balanced environment:
\begin{itemize}
\item  For any $x\in A$ and any $\kappa\in \square$, there exists a $\kappa$-path from $x$ to some $y\in\partial A_x$, and this path is contained in $\bar{A}_x$.
\end{itemize}

This property gives us a useful inequality.
\begin{lemma}\label{CPhicontrol}
For $x\in A\cap \Delta_N$, if $l_x \le N$, then
\begin{equation}\label{CLTf7}
\Phi_N (x) \le \xi_0 ^{-l_x} \sum_{y \in \partial A_x \cap \Delta_N} \Phi_N (y).
\end{equation}
\end{lemma}
\pf Suppose that $A_x\neq\emptyset$ (otherwise the proof is trivial). Since $l_x \le N$,  $\bar{A}_x\subset \Delta_N (x)$. Note that at least one of the $2^d$ 
corners of $\Delta_N (x)$ is contained in $\Delta_N$. Without loss of generality, suppose that
$v=x+(N,\cdots,N)\subset \Delta_N$. Then there is a $(1,\cdots, 1)$-path in $\bar{A}_x$ from $x$ to
some $y\in \partial A_x\cap \Delta_N$,
as illustrated in the following figure:
\begin{center}
\begin{tikzpicture}[scale=.2]{centered}
\fill[gray!30!white][rotate=45] (10, 1.5) ellipse (4.5 and 4.3);
\draw[step=1, gray, very thin] (1.5,3.5) grid (10.5, 12.5);
\draw[dash pattern=on 2pt off 3pt on 4pt off 4pt] (-5,-5) rectangle (15, 15) node[right=1pt]{$v$};
\draw (2,4) rectangle (22, 24);
\draw[thick] (5,5)--(5,6)--(7,6)--(7,7)--(8,7)--(8,10)--(9,10)--(9,11)--(10,11)node[right=1pt]{$y$};
\draw[left=1pt] (5,5) node{$x$};
\draw (3,10) node[right=1pt]{$A_x$};
\draw (23, 14) node[right=1pt]{$\Delta_N$};
\draw (-5, 5) node[left=1pt]{$\Delta_N(x)$};
\end{tikzpicture}
\end{center}
Recalling that $\Phi_N$ is the invariant measure for $\{\hat{X}_{n,N}\}$ defined in Section 1, we have
\begin{align*}
\Phi_N (y) 
&=
 \sum_{z\in \Delta_N} \Phi_N (z) P_{\omega^N}^{d(x,y)} (\hat{z}, \hat{y})\\
&\ge 
 \Phi_N (x) P_{\omega^N}^{d(x,y)}(\hat{x}, \hat{y})
\ge  \Phi_N (x) \xi_0^{l_x}.
\end{align*}
Here $P^m_{\omega^N} (\hat{z}, \hat{y})$ denotes the $m$-step transition probability
of $\{\hat{X}_{n,N}\}$ from $\hat{z}$ to $\hat{y}$.
\qed

Let $S_n=\{x: |x|_{\infty}=n\}$ denote the boundary of $\Delta_n$.
Let $x\to y$ be the event that $y\in \bar{A}_x$
and $o\to S_n$ be the event that $o\to x$ for some $x\in S_n$.
The following theorem, which is the site percolation version of the combination of
Theorems 6.10 and 6.14 in \cite{GG}, gives an exponential bound
on the diameter of the cluster containing the origin,  when $p$ is small.
\begin{theorem}\label{perc}
There exists a function $\phi(p)$ of $p=p(\varepsilon_0)$ such that
$$P\{o\to S_n\}\le C n^{d-1} e^{-n\phi(p)}$$
and $\lim_{p\to 0}\phi(p)=\infty$.
\end{theorem}

Let $A_x(n)$ denote the connected component of $A_x\cap \Delta_n(x)$ that contains $x$ and set 
\[q_n=P\{o\to S_n\}.\]
 The proof of Theorem
\ref{perc} will proceed by showing some (approximate) subadditivity
properties of $q_n$. We thus recall the following subadditivity lemma:
\begin{lemma}
	\label{Clem-subadd}
If a sequence of finite numbers $\{b_k: k\ge 1\}$ is subadditive, that is,
$b_{m+n}\le b_m +b_n \mbox{   for all m,n}$,
then $\lim_{k\to\infty}b_k/k=\inf_{k\in\mathbb{N}} b_k/k$.
\end{lemma}
\noindent{\it Proof of Theorem \ref{perc}:}
We follow the proof given by Grimmett in \cite{GG} in the bond percolation case. 
By the BK inequality (\cite{GG}, pg. 38),
$$q_{m+n}\le \sum_{x\in S_m} P\{o\to x\} P\{x\to x+S_n\}.$$
But $P\{o\to x\}\le q_m$ for $x\in S_m$
and $P\{x\to x+S_n\}=q_n$ by translation invariance. Hence we get
\begin{equation}\label{Csub1}
q_{m+n}\le |S_m|q_m q_n.
\end{equation}
By exchanging $m$ and $n$ in (\ref{Csub1}),
\begin{equation}\label{Csub11}
q_{m+n}\le |S_{m\wedge n}|q_m q_n.
\end{equation}

On the other hand, let $U_x$ be the event that 
$x\in \overline{A_o(m)}$
 and let $V_x$ be the event that 
$\overline{A_x(n)}\cap S_{m+n}\neq\emptyset$. 
We use the FKG inequality (\cite{GG}, pg. 34) to find that
$$q_{m+n}\ge P\{U_x\}P\{V_x\}\quad \mbox{ for any $x\in S_m$}.$$
However, $\sum_{x\in S_m}P\{U_x\}\ge q_m$, which implies that
$$\max_{x\in S_m}P\{U_x\}\ge \frac{q_m}{|S_m|} .$$
Let $\gamma_n=P\{\overline{A_o(n)}\cap\{x:x_1=n\}\neq\emptyset\}$, then $P\{V_x\}\ge\gamma_n$.
 Moreover, $\gamma_n\le q_n\le 2d\gamma_n$.
Hence
\begin{equation*}
q_{m+n}\ge \frac{q_m q_n}{2d |S_m|},
\end{equation*}
and then
\begin{equation}\label{Csub2}
q_{m+n}\ge \frac{q_m q_n}{2d |S_{m\wedge n}|}.
\end{equation}

Note that $|S_m|\le C_d m^{d-1}$. Letting
$$b_k =\log q_k +\log C_d +(d-1)\log (2k),$$
one checks
using (\ref{Csub11})
that the sequence $\{b_k\}$ is subadditive.
Similarly by (\ref{Csub2}),  $\{-\log q_k +\log (2d C_d) +(d-1)\log (2k)\}$
is subadditive.
Thus, using Lemma \ref{Clem-subadd},
$$\phi (p):=-\lim_{k\to\infty}\frac{1}{k}\log q_k$$
exists and
\begin{equation}\label{Csub3}
\log q_k +\log C_d +(d-1)\log (2k)\ge -k\phi (p),
\end{equation}
\begin{equation}\label{Csub4}
-\log q_k +\log (2d C_d) +(d-1)\log (2k)\ge k\phi(p).
\end{equation}
The first part of the theorem follows simply from (\ref{Csub4}), and the second
by noting that with
$p\downarrow 0$ in (\ref{Csub3}) we have
$q_k \downarrow 0$ and then  $\phi (p)\to \infty$. \qed
\begin{remark}
	It follows from Theorem \ref{perc}
	that
\begin{equation}\label{Clo}
P\{l_o \ge n\} \le P\{o\to S_{\lfloor n/2d\rfloor}\}\le C e^{\phi(p)} n^{d-1} e^{-n\phi(p)/2d}.
\end{equation}
With (\ref{Clo}) and the Borel-Cantelli lemma, one concludes
that, P-almost surely, $l_x\le N$ is true for all
$x\in \Delta_N$ when $N$ is sufficiently large and $p$ is such that $\phi(p)>0$.
Hence the inequality (\ref{CLTf7}) holds for all $x\in\Delta_N$ when $N$
is large.
\end{remark}

\noindent{\it Proof of (\ref{CPhi0}) under the assumption of Theorem \ref{CLT2}(i):}
By H\"older's inequality,
$$\frac{1}{|\Delta_N|}\sum_{y\in \partial A_x \cap \Delta_N} \Phi(y)\le
\lVert         \Phi_N 1_{\partial A_x}\rVert _{\Delta_N, \beta}
\big(\frac{|\partial A_x|}{|\Delta_N|}\big)^{1-1/\beta} ,$$
so when $N$ is large enough we have by Lemma \ref{CPhicontrol} that for any $x\in A\cap \Delta_N$,
\begin{equation}\label{CLTf8}
\Phi_N (x)\le \xi_0^{-l_x} |\partial A_x|^{1-1/\beta	} |\Delta_N |^{1/\beta}
\lVert         \Phi_N 1_{\partial A_x}\rVert _{\Delta_N, \beta}.
\end{equation}

Hence for any $\alpha \in (1, \beta)$,
\begin{align*}
&\lVert   \Phi_N 1_A\rVert _{\Delta_N,\alpha}^\alpha\\
&\le 
 \frac{1}{|\Delta_N|}\sum_{x\in A\cap\Delta_N}\big(\xi_0^{-l_x} |\partial A_x|^{1-1/\beta} 
 |\Delta_N |^{1/\beta}\lVert   \Phi_N 1_{\partial A_x}\rVert _{\Delta_N, \beta} \big)^\alpha\\
&\le 
\left[\frac{1}{|\Delta_N|}\sum_{x\in A\cap\Delta_N}
\big(\xi_0^{-l_x} |\partial A_x|^{1-1/\beta}|A_x|^{1/\beta})^{\alpha (\beta/\alpha)'}\right]
^{1-\alpha/\beta}\\
&\qquad\times \left[\frac{1}{|\Delta_N|}\sum_{x\in A\cap\Delta_N}
\big(\frac{|\Delta_N|^{1/\beta}\lVert   \Phi_N 1_{\partial A_x}\rVert _{\Delta_N, \beta}}
{|A_x|^{1/\beta}}\big)^{\beta}
\right]^{\alpha/\beta}\\
&=
 \left[\frac{1}{|\Delta_N|}\sum_{x\in A\cap\Delta_N}
\big(\xi_0^{-l_x} |\partial A_x|^{1-1/\beta}|A_x|^{1/\beta}\big)^{\alpha\beta/(\beta -\alpha)}\right]
^{1-\alpha/\beta}\\
&\qquad\times
\left(\sum_{x\in A\cap\Delta_N}
\frac{\lVert  \Phi_N 1_{\partial A_x}\rVert _{\Delta_N, \beta}^{\beta}}{|A_x|}
\right)^{\alpha/\beta},
\end{align*}

where we used (\ref{CLTf8}) in the first inequality and H\"older's inequality in the second.

Observe that
\begin{equation}\label{CPhiaverage}
\sum_{x\in A\cap\Delta_N} \dfrac{\lVert         \Phi_N 1_{\partial A_x}\rVert _{\Delta_N, \beta}^{\beta}}{|A_x|}
\le \sum_{i=1}^n \lVert         \Phi_N 1_{\partial A_i}\rVert _{\Delta_N, \beta}^{\beta}
\le 2d \lVert         \Phi_N 1_{\partial A}\rVert _{\Delta_N, \beta}^{\beta}
\le C\varepsilon_0^{-\beta},
\end{equation}
where $A_1, \cdots, A_n$ are different clusters that intersect with $\Delta_N$.
On the other hand, the multidimensional ergodic theorem gives
\begin{align}\label{percergodic}
&\lim_{N\to\infty}\frac{1}{|\Delta_N|}\sum_{x\in A\cap\Delta_N}
\big(\xi_0^{-l_x} |\partial A_x|^{1-1/\beta}|A_x|^{1/\beta}\big)^{\alpha\beta/(\beta -\alpha)}\nonumber\\
&= E \big(\xi_0^{-l_o} |\partial A_o|^{1-1/\beta}|A_o|^{1/\beta}\big)^{\alpha\beta/(\beta -\alpha)}
\le
C E  \big(\xi_0^{-l_o} l_o^d\big)^{\alpha\beta/(\beta -\alpha)} \quad \mbox{P-a.s.,}
\end{align}
which by (\ref{Clo}) is finite when $\varepsilon_0$ is small.
\qed

\section{Transience in general ergodic environments}\label{SeTran}
In this section we will prove (ii) of Theorem \ref{CLT1} by an argument similar
to that in
\cite{ZO}. The main differences in our method are that we use a stronger control of the hitting time
(Lemma \ref{Ctau}), and that we apply a mean value
inequality (Theorem \ref{Cmvi}) instead of the discrete Harnack
inequality used in \cite{ZO}.

\begin{lemma}\label{Ctau}
Let $\{X_n\}$ be a
random walk in a balanced environment
$\omega$ such that $\omega(x,o)=0$ for all $x$. For
any $r>0$, define $\tau=\tau(r)=
\inf\{n: |X_n|>r\}$. Then $E_\omega^o \tau\le (r+1)^2$.
\end{lemma}
\pf Observe that $\{|X_n|^2-n\}$ is a (quenched) martingale with respect to
$\{\mathcal{F}_n=\sigma(X_1,\cdots,X_n)\}$.
Thus by optional stopping, $0= E_\omega^o[|X_\tau|^2-\tau]\le (r+1)^2-E_\omega^o\tau$. \qed

To prove Theorem \ref{CLT1}(ii), we shall make use of the following mean-value inequality,
which is a modification of Theorem 3.1 in \cite{KT}. Let $B_{r}(z)=\{x\in\mathbb{Z}^d: |x-z| <r\}$. 
We shall also write $B_r (o)$  as $B_r$; recall the definition of $L_\omega$ in
\eqref{operator}.
\begin{theorem}\label{Cmvi}
For any function $u$ on $\bar B_R (x_0)$ such that
$$L_{\omega} u =0 , \quad x \in B_R (x_0)$$
and any $\sigma\in (0,1)$, $0<p\le d$,  we have
\[
\max_{B_{\sigma R}(x_0)}u\le C \lVert
\frac{u^+}{\varepsilon^{d/p}}\rVert _{B_R (x_0), p},
\]
where $C$ depends on $\sigma$, $p$ and $d$.
\end{theorem}
We postpone the proof of Theorem \ref{Cmvi} to 
the next section, and now demonstrate Theorem \ref{CLT1}.

\noindent{\it Proof of Theorem \ref{CLT1}(ii):}
Note that  the transience of the random walk would not change
if we considered the walk restricted
to its jump times. That is, the transience or recurrence of the random walk
 in an environment $\omega$ is the same as in an environment
$\tilde \omega$,
where $\tilde \omega$ is defined by
$\tilde{\omega}(x,e)=\omega (x,e)/(1-\omega(x,o))$. Therefore, in the sequel we assume $\omega(x,o)=0$
for all $x$ and almost all $\omega$.

Let $K$ be any constant that is at least 3.
We denote $B_{K^i}(x)$ by $B^i(x)$ and define
$\tau_i:=\inf \{n: |X_n|> K^i\}$.
Our approach is to bound the (annealed) expected number of visits to the origin by
the walk; this requires some a-priori bounds on the moments of
$\varepsilon(o)^{-1}$.\\

For any $z\in\partial B^i$ , $y\in B^{i-1}$, noting that $E_\omega^x (\mbox{\# visits at $y$ before $\tau_{i+2}$}):=v(x)$ satisfies $L_\omega v (x)=0$ for
$x\in B^{i+2}\setminus\{y\}$, we have that, for $p\in(0,d]$,
\begin{align}\label{CLTf11}
& E_{\theta^{y}\omega}^{z} (\mbox{ \# visits at $o$ before $\tau_{i+1}$})
\nonumber\\
&\le 
 E_{\omega}^{z+y}(\mbox{\# visits at $y$ before $\tau_{i+2}$})\nonumber\\
&\le 
 \max_{x\in B^{i-1}(z)}E_{\omega}^{x}(\mbox{\# visits at $y$ before $\tau_{i+2}$})\nonumber\\
&\le 
 C\Bnorm{\frac{E_{\omega}^{x}(
               \mbox{\# visits at $y$ before $\tau_{i+2}$})}{\varepsilon_\omega (x)^{d/p}}
              } _{B_{2K^{i-1}}(z), p}\nonumber\\
&\le 
 C\Bnorm{ \frac{E_{\omega}^{x}
              (\mbox{\# visits at $y$ before $\tau_{i+2}$})}{\varepsilon_\omega (x)^{d/p}}
              } _{B^{i+2}, p},
\end{align}
where we used Theorem \ref{Cmvi} in the third inequality. Take $p=d/q$ (without loss of
generality, we always assume that $q< d$).
Then, by (\ref{CLTf11}) and Lemma \ref{Ctau},
\begin{align}\label{CLTf12}
& \sum_{y\in B^{i-1}}E_{\theta^{y}\omega}^o
(\mbox{ \# visits at $o$ in $[\tau_i,\tau_{i+1})$})\nonumber\\
&\le 
 C\sum_{y\in B^{i-1}}  \left[
\frac{1}{|B^{i+2}|}\sum_{x\in B^{i+2}}
\frac{
E_{\omega}^{x}
(\mbox{\# visits at $y$ before $\tau_{i+2}$})^{d/q}}
{\varepsilon_\omega (x)^d}
\right]^{q/d}\nonumber\\
&\le 
 C K^{-iq}\sum_{y\in B^{i-1}}\sum_{x\in B^{i+2}}
    \frac{E_{\omega}^{x}(\mbox{\# visits at $y$ before $\tau_{i+2}$})}
    {\varepsilon_\omega (x)^q}
        \nonumber\\
&= 
 C K^{-iq}\sum_{x\in B^{i+2}}
    \frac{E_{\omega}^{x}(\mbox{\# visits at $B^{i-1}$ before $\tau_{i+2}$})}
    {\varepsilon_\omega (x)^q}\nonumber\\
&\le 
 C K^{-iq}\sum_{x\in B^{i+2}}\frac{E_{\omega}^{x}\tau_{i+2}}{\varepsilon_\omega (x)^q}\nonumber\\
&\le 
 C K^{(2-q)i} \sum_{x\in B^{i+2}}       \varepsilon_\omega(x)^{-q}.
\end{align}
Taking expectations and using translation invariance, we have
\begin{equation*}
 \mathbb{E}^o(\mbox{\# visits at $o$ in $[\tau_i, \tau_{i+1})$})
\le C K^{(2-q)i} E \varepsilon^{-q}.
\end{equation*}
Therefore, if $E \varepsilon^{-q}<\infty$ for some $q>2$ , then
\begin{equation*}
\mathbb{E}^o (\mbox{\#
visits at $o$})\le C E \varepsilon^{-q}
                     \sum_{i=1}^{\infty}K^{(2-q)i}<\infty .
\end{equation*}

This proves Theorem \ref{CLT1}(ii) for $\{\Omega, P\}$ such that $\omega (x,o)=0$ for all $x$
and almost all $\omega$.
As mentioned earlier,
the general case follows by
replacing
$\varepsilon$ with $\varepsilon/(1-\omega(o,o))$.\qed\\

\begin{remark}
It is natural to expect
 that arguments similar to the proof of the invariance principle also work for
proving the transience in the iid case. Namely, one may hope to control
$P_\omega^x\{\mbox{visit $o$ in $[\tau_i,\tau_{i+1})$}\}$ using some mean value inequality (like Theorem \ref{Cmvi}), and to use percolation arguments to handle 
``bad sites''
where the ellipticity constant $\varepsilon$ is small.

This suggests considering walks that jump from bad sites to 
good sites. In \cite{KT2}, Kuo and Trudinger
proved a maximum principle and mean value inequality for balanced operators 
in general meshes, which may be applied to balanced walks with possibly
big jumps. However, their estimates, 
in the presence of a small ellipticity
constant, are not strong enough. To overcome this issue,
we will prove a modified maximum principle that involves only 
big exit probabilities, and then use it to prove the transience 
in the i.i.d case with no moment assumptions.
\end{remark}

\section{Transience in iid environments}\label{SeTriid}
In this section we prove a modified maximum principle for 
balanced environments.
We then prove Theorem 2(ii) using the corresponding mean 
value inequality (Theorem \ref{Cmvi2})
and percolation arguments.

\subsection{Difference operators}\label{SeTriid1}
Following \cite{KT2}, we introduce general difference operators.
Let $a$ be a nonnegative function on $\mathbb{Z}^d\times \mathbb{Z}^d$ 
such that for any $x$,
$a(x,y)> 0$ 
for only finitely many $y$. 
Define the linear operator $L_a$ acting on the set of functions
on $\mathbb{Z}^d$ by
\[L_a f(x)=\sum_y a(x,y)(f(y)-f(x)).\]
We say that $L_a$ is \textit{balanced} if
\begin{equation}\label{CLTe7}
\sum_y a(x,y)(y-x)=0.
\end{equation}
Throughout this section we assume that $L_a$ is a probability operator, that is,
\[\sum_y a(x,y)=1.\]
For any finite subset $E\subset\mathbb{Z}^d$, define its boundary
\[E^b=E^b(a)=\{y\notin E: a(x,y)>0 \text{ for some } x\in E
      \}, 
\] 
and set 
\begin{equation}
	\label{CLTeq-ofernew}
	\tilde{E}=E\cup E^b.
\end{equation}
Define the upper contact set of $u$ at $x\in E$ as
$$I_u(x)=I_u(x,E,a)=\{s\in\mathbb{R}^d: u(x)-s\cdot x\ge u(z)-s\cdot z \text{ for all }z\in\tilde{E}\}.$$
Set 
\begin{align*}
&h_x=h_x(a)=\max_{y: a(x,y)>0}\abs{x-y},\\
&b(x)=\sum_y a(x,y)(y-x), \mbox{ and }b_0=\sup|b|.
\end{align*}
Note that $b_0=0$ when $L_a$ is balanced.

The following lemma is useful in the proofs of various mean 
value inequalities. It is similar to Theorem 2.2 in \cite{KT2}, 
except that the proof in \cite{KT2} contains several unclear
 passages, e.g., in the inequality above (2.23) in \cite{KT2},
and so we provide a complete proof.
Throughout, we set $u^+=u\vee 0$. 
\begin{lemma}\label{Cmvilemma}
Fix $R>0$. Let $\eta(x)=\eta_R(x):=(1-\abs{x}^2/R^2)^\beta 1_{|x|<R}$ be a
function on $\mathbb{R}^d$.
For any function $u$ on $B_R$ such that 
$L_a u\ge 0$ in $B_R$ and any $\beta\ge 2$, 
we let $v=\eta u^+$.
Then, for any $x\in B_R$ with $I_v (x)=I_v(x, B_R, a)\neq\emptyset$,
\[L_a v(x)\ge -C(\beta, b_0R) \eta^{1-2/\beta}R^{-2} h_x^2 u^+,\]
where $C(\beta,b_0R)$ is a constant that depends only on $\beta$ and $b_0R$.
\end{lemma}
\pf 
We only need to consider the nontrivial case where $v\not\equiv 0$.
For $s=s(x)\in I_v(x)\neq \emptyset$,  recalling the definition of $I_v$ one has
$$|s|\le 2v(x)/(R-|x|).$$ 
Note that $I_v(x)\neq\emptyset$ implies $u(x)> 0$.
If further $R^2-|x|^2\ge 4R \abs{x-y}$ , 
computations as in \cite[pg. 426]{KT2}
reveal that
\begin{align}
 2^{-\beta}
 &\le
 \frac{\eta(y)}{\eta(x)}\le 2^\beta,\label{Cfirst}\\
 \abs{\eta(x)-\eta(y)}
 &\le
  \beta 2^\beta R^{-1} \eta(x)^{1-1/\beta}|x-y|,\label{Csecond}\\
 \abs{\eta(x)-\eta(y)-\nabla\eta(x)(x-y)}
 &\le
 \beta (\beta-1)2^{\beta} R^{-2}\eta(x)^{1-2/\beta}|x-y|^2,\label{Cthird}\\
 |s|
 &\le 
 4 \eta^{1-1/\beta}R^{-1}u,\label{Cfourth}
 \end{align}
where 
\begin{equation}\label{Afifth}
\nabla\eta=-2\beta x R^{-2}\eta^{1-1/\beta}
\end{equation}
 is the gradient of $\eta$.
 Following \cite{KT2}, we set $w(z)=v(z)-s\cdot (z-x)$. 
By the definition of $s$, we have
 $w(x)\ge w(z)$ for all $z\in \tilde E$.
Then 
\begin{align}\label{A*6}
v(x)-v(y)=\frac{\eta(x)}{\eta(y)}(v(x)-v(y))
  &+\frac{\eta(y)-\eta(x)}{\eta(y)}s(x-y)\\
  &+\frac{\eta(y)-\eta(x)}{\eta(y)}(w(x)-w(y)).\nonumber
\end{align}
Consider first $x$ such that $R^2-|x|^2\ge 4Rh_x$.
By (\ref{Csecond}), for any $y$ such that $a(x,y)>0$,
\begin{equation}\label{A*7}
\frac{\eta(y)-\eta(x)}{\eta(y)}(w(x)-w(y))
\le \beta 2^\beta R^{-1}h_x\eta(x)^{-1/\beta}\frac{\eta(x)}{\eta(y)}(w(x)-w(y)).
\end{equation}
Since 
\begin{align*}
&\sum_y a(x,y)\frac{\eta(x)}{\eta(y)}(w(x)-w(y))\\
&=\sum_y a(x,y) \Big[\frac{\eta(x)}{\eta(y)}\big(v(x)-v(y)\big)
                  +\frac{\eta(x)-\eta(y)}{\eta(y)}s (y-x)\Big]+s\cdot b(x),
\end{align*}
by \eqref{A*6}, \eqref{A*7} and noting $R^{-1}\eta^{-1/\beta}h_x\le 1/4$, we obtain
\begin{align}\label{A*8}
&\sum_y a(x,y)(v(x)-v(y))\nonumber\\
&\le 
(1+\beta 2^\beta R^{-1}h_x\eta(x)^{-1/\beta})
\sum_y a(x,y)\Big[\frac{\eta(x)}{\eta(y)}\big(v(x)-v(y)\big)+\frac{\eta(x)-\eta(y)}{\eta(y)}s (y-x)+s\cdot b(x)\Big]
\nonumber\\
&\le 
\beta 2^{\beta-1}\big[\sum_y a(x,y) \frac{\eta(x)}{\eta(y)}\big(v(x)-v(y)\big)+
4(\beta 2^{2\beta}+b_0R)\eta^{1-2/\beta}R^{-2}h_x^2u\big],              
\end{align}
where we used \eqref{Cfirst}, \eqref{Csecond} \eqref{Cfourth} in the last inequality.
Moreover,  recalling that $u(x)>0$ (because $I_v(x)\neq\emptyset$),
\begin{align}\label{A*9}
& \sum_y a(x,y) \frac{\eta(x)}{\eta(y)}\big(v(x)-v(y)\big)\nonumber\\
&= 
\sum_y a(x,y)\left [\eta(x)\big(u(x)-u^+(y)\big)+\big(\eta(x)-\eta(y)\big)u(x)+\frac{(\eta(x)-\eta(y))^2}{\eta(y)}u(x)\right]\nonumber\\
&\stackrel{a\geq 0}{\le}  
-\eta(x)L_a u(x)+\sum_y a(x,y) \left[\big(\eta(x)-\eta(y)\big)u(x)+\frac{(\eta(x)-\eta(y))^2}{\eta(y)}u(x)\right]\nonumber\\
&\stackrel{L_au\ge 0}{\le}
 \sum_y a(x,y) \left[\big(\eta(x)-\eta(y)-\nabla\eta(x)(x-y)\big)u(x)+\frac{(\eta(x)-\eta(y))^2}{\eta(y)}u(x)\right]\nonumber\\
 &\qquad\qquad-\nabla\eta(x)b(x)u(x)\nonumber\\
&\le
 (2^{3\beta+1}+b_0R)\beta^2\eta^{1-2/\beta}h_x^2R^{-2}u,
\end{align}
where we used (\ref{Cfirst}), (\ref{Csecond}), (\ref{Cthird}) and (\ref{Afifth}) in the last inequality.
Hence, by (\ref{A*8}) and (\ref{A*9}), we conclude that
\begin{equation*}
-L_a v
\le
(2^{3\beta+1}+b_0R)\beta^32^\beta \eta^{1-2/\beta}R^{-2}h_x^2 u
\end{equation*}
holds in $\{x: R^2-|x|^2\ge 4Rh_x, I_v(x)\neq\emptyset\}$.

On the other hand, if $R^2-|x|^2<4Rh_x$,  then $\eta^{1/\beta}\le 4h_x/R$. Thus by the fact that $u(x)>0$, 
we have $-L_a v\le v(x)\le 16\eta^{1-2/\beta} R^{-2} h_x^2 u$. \qed\\

\noindent{\it Proof of Theorem \ref{Cmvi}:}
Since $L_\omega$ is a balanced operator ($b_0=0$) and $h_x=1$ in this case,
by the above lemma, 
\[L_\omega v\ge -C(\beta) \eta^{1-2/\beta}R^{-2} u\] for $x\in B_R$ such that $I_u(x)\neq\emptyset$,
where $C(\beta)$ depends only on $\beta$.
Applying Theorem \ref{CMP} to $v$ and taking $\beta=2d/p\ge 2$, 
we obtain
\begin{align*}
\max_{B_R} v 
&\le 
C \Bnorm{ \eta^{1-2/\beta} \frac{u^+}{\varepsilon}} _{B_R, d}
=C\Bnorm{    v^{1-p/d}\frac{(u^+)^{p/d}}{\varepsilon}}_{B_R, d}\\
&\le 
C(\max_{B_R} v)^{1-p/d}\Bnorm{\frac{u^+}{\varepsilon^{d/p}}} _{B_R, p}^{p/d}.
\end{align*}
Hence 
\[
\max_{B_R} v \le C \Bnorm{\frac{u^+}{\varepsilon^{d/p}}} _{B_R, p},
\]
and then
\[
\max_{B_{\sigma R}}u\le (1-\sigma^2)^{-2d/p}\max_{B_{\sigma R}}v
\le C(\sigma, p, d) \Bnorm{\frac{u^+}{\varepsilon^{d/p}}} _{B_R, p}. \text{\qed}
\]

\subsection{A new maximum principle and proof of Theorem \ref{CLT2}(ii)}

For any fixed environment $\omega\in\Omega$, let $\varepsilon_0>0$ be a constant to be determined,
and define site percolation as in Section \ref{SePercE}. Recall that for $x\in\mathbb{Z}^d$, $A_x$
is the percolation cluster that contains $x$ and $l_x$ is its $l^1$-diameter.
As mentioned in the introduction, the transience would not change if we 
considered the walk restricted to its
jump times. Without loss of generality, we assume that $\omega(x,o)=0$ for all $x$, $P$-almost surely.

Recall the definition of $\square$ and $\kappa$-path for $\kappa\in\square$ in Section \ref{SePercE}. Note
that under our assumption, $\max_i \omega(x, e_i)\ge 1/2d$, so we take $\xi_0=1/2d$ in the definition of 
$\kappa$-paths.

For each $\kappa\in\square$, we pick a site $y_\kappa=y(x, \kappa)\in \partial A_x$ such that
\[d(x, y_\kappa)=\max_{\substack{y:
                  \exists \text{ $\kappa$-path in $\bar{A}_x$ }\\\text{ from $x$ to $y$}
                           }
                 } d(x,y)\] 
and let $\Lambda_x\subset \bar{A}_x$ be the union of (the points of the) $\kappa$-paths from $x$ to $y_\kappa$ over all $\kappa\in\square$.
From the definition of $y_\kappa$ one can conclude that
\begin{itemize}
\item For any $q\in\mathbb{R}^d$, we pick a $\kappa=\kappa_q\in \square$ such that
\[q_j \kappa_j\le 0 \text{ for all }j=1,\cdots, d.\]
Then $(y_\kappa-x)_j q_j\le 0$ for all $j=1,\cdots, d$. Moreover, for $i\in\{1,\cdots, d\}$, $q_i>0$ implies $y_\kappa-e_i\notin \Lambda_x$, and $q_i<0$
implies $y_\kappa+e_i\notin \Lambda_x$.
\end{itemize}

In the sequel we let $\tau_{\Lambda_x}=\inf\{n>0: X_n\notin \Lambda_x\}$ and 
\[a(x,y)=P_\omega^x \{
                      X_{\tau_{\Lambda_x}}=y
                     \}.
\]
By the fact that $X_n$ is a (quenched) 
martingale, it follows that $L_a$ is a balanced operator.

For the statement of the next theorem,
recall the definition of $\tilde{E}$ in \eqref{CLTeq-ofernew}.
\begin{theorem}\label{Cmp2}
Let $E\subset\mathbb{Z}^d$ be bounded. 
Let $u$ be a function on $\tilde{E}$. If $L_a u(x)\ge -g(x)$ for all $x\in E$ such that $I_u(x)=I_u(x,E,a)\neq \emptyset$ ,
then
$$\max_E u\le 
  \frac{d\diam \tilde{E}}{\varepsilon_0}
  \bigg(\sum_{\substack{
                           x\in E\\I_u (x)\ne \emptyset
                       }
             }
                      \abs{g(x)(2d)^{l_x}}^d
  \bigg)^{\frac{1}{d}}
 +\max_{E^b}u  .$$
\end{theorem}
\pf Without loss of generality, assume $g\ge 0$ and
$$\max_E u=u(x_0)>\max_{E^b}u$$
for some $x_0\in E$. Otherwise, there is nothing to prove.

For $ s\in \mathbb{R}^{d} $ such that $ |s|_{\infty} \le [u(x_0)-\max_{E^b}u]/(d\diam \tilde{E}) $,
we have \[ u(x_0)-u(x) \ge s \cdot (x_0-x) \] for all $ x \in E^b $,
 which implies that
$\max_{z\in\tilde{E}}u(z)-s\cdot z$
is achieved in $E$. Hence
  $ s \in \bigcup_{x \in E} I_{u}(x) $ and

\begin{equation}\label{CLTe1}
 \left[-\dfrac{u(x_0)-\max_{E^b }u}{d\diam \tilde{E}}, \dfrac{u(x_0)-\max_{E^b}u}{d\diam \tilde{E}}\right]^{d} \subset \bigcup_{x\in E} I_{u}(x)  .
 \end{equation}

Further, if $s\in I_u(x)$, we set
$$w(z)=u(z)-s(z-x).$$
Then $w(z)\le w(x)$ for all $z\in \tilde{E}$ and 
\begin{equation}\label{CLTe2}
I_u(x)=I_w(x)+s.
\end{equation}
Since for any $q\in I_w(x)$,  there is $\kappa=\kappa_q\in\square$ such that
$$q_j(x-y_\kappa)_j\ge 0 \text{ for } j=1,\cdots, d,$$ we have
$$ w(x)-w(y_\kappa\pm e_i)\ge q(x-y_\kappa\mp e_i)\ge \mp q_i.$$
Moreover, for any $i\in\{1,\cdots, d\}$, if $q_i>0$, then $y_\kappa-e_i\notin \Lambda_x$ and we have
$w(x)-w(y_\kappa-e_i)\ge |q_i|$. Similarly, if $q_i<0$, then $y_\kappa+e_i\notin \Lambda_x$ and 
$w(x)-w(y_\kappa+e_i)\ge |q_i|$. We conclude that
$$|q_i|\le \frac{\sum_{y} a(x, y)(w(x)-w(y))}{\min\limits_{\pm}\{a(x, y_\kappa\pm e_i)\}}.$$
On the other hand, from the construction of $\Lambda_x$ we obtain (noting that
$y_\kappa\in\partial A_x$)
$$a(x, y_\kappa\pm e_i)\ge (\frac{1}{2d})^{l_x} \varepsilon_0.$$
Hence, since $L_a$ is balanced,
$$|q_i|\le 
\frac{(2d)^{l_x}}{\varepsilon_0}\sum_y a(x,y)(w(x)-w(y))
=
\frac{(2d)^{l_x}}{\varepsilon_0}(-L_a u)
\le \frac{(2d)^{l_x}}{\varepsilon_0} g$$
for all $i$. Therefore
\begin{equation}\label{CLTe3}
I_w(x)\subset [-(2d)^{l_x}\varepsilon_0^{-1} g, (2d)^{l_x}\varepsilon_0^{-1} g]^d.            
\end{equation}
Combining (\ref{CLTe1}), (\ref{CLTe2}) and (\ref{CLTe3}), we conclude that
\begin{equation*}
\left(\dfrac{u(x_0)-\max_{E^b}u}{d\diam \tilde{E}}\right)^d\le 
\sum_{\substack{
                           x\in E\\I_u (x)\ne \emptyset
                       }
             }
                      \abs{g(x)(2d)^{l_x}\varepsilon_0^{-1}}^d. 
\quad\quad\quad\quad\quad
\quad\quad\quad\quad\quad
 \mbox{\qed}
\end{equation*}

As with Theorem \ref{Cmvi}, we have a corresponding mean value inequality.
\begin{theorem}\label{Cmvi2}
For any function $u$ on $B_R$ such that
$$L_a u=0, \quad x\in B_R$$ 
and any $\sigma\in (0,1)$, $0<p\le d$, we have
\[
\max_{B_{\sigma R}} u\le C\big(\frac{\diam \tilde{B}_R}{\varepsilon_0 R}\big)^{d/p}
                           \norm{[l_x^2 (2d)^{l_x}]^{d/p}u^+}_{B_R, p},
\]
where $C$ depends on $\sigma, p$ and $d$.
\end{theorem}
\pf By the same argument as in the proof of Theorem \ref{Cmvi}, Lemma \ref{Cmvilemma} and 
Theorem \ref{Cmp2} implies Theorem \ref{Cmvi2}. \qed\\

Having established
Theorem \ref{Cmvi2}, we can now 
prove the transience of the random walks in balanced
iid environment with $d\geq 3$.

\noindent{\it Proof of Theorem \ref{CLT2}(ii)}:
Let $K$ be any constant $\ge 4$ and define $B^i, \tau_i$ as in Section 5.
Let $\Omega_i=\{\omega\in\Omega: l_x\le K^{i-1} \mbox{ for all $x\in B^{i+2}$}\}$. 
For any $\omega\in\Omega_i$, $z\in\partial B^i$, $y\in B^{i-1}$, noting that
$P_\omega^x \{\mbox{visit $y$ before $\tau_{i+2}$}\}:=u(x)$ satisfies
\[L_a u(x)=0\] 
for $x\in B_{2K^{i-1}}(z)$, by similar argument as in (\ref{CLTf11}) we have
\begin{align*}
& P_{\theta^y\omega}^z \{\mbox{visit $o$ before $\tau_{i+1}$}\}1_{\omega\in\Omega_i}\\
&\le 
\max_{x\in B^{i-1}(z)}P_\omega^x \{\mbox{visit $y$ before $\tau_{i+2}$}\}1_{\omega\in\Omega_i}\\
&\le 
C \varepsilon_0^{-d}\norm{
                                [l_x^2 (2d)^{l_x}]^d P_\omega^x \{\mbox{visit $y$ before $\tau_{i+2}$}\}                              
                               }_{B_{2K^{i-1}}(z), 1}\\
&\le 
C \varepsilon_0^{-d}\abs{B^{i+2}}^{-1}
         \sum_{x\in B^{i+2}}l_x^{2d} (2d)^{dl_x} P_\omega^x \{\mbox{visit $y$ before $\tau_{i+2}$}\} ,                            
\end{align*}
where in the second inequality, we applied Theorem \ref{Cmvi2} with
$p=1$ and used the fact that $\diam \tilde{B}_{2K^{i-1}}\le 3K^{i-1}$ when $\omega\in\Omega_i$. Hence
\begin{align}\label{CLTe5}
&\sum_{y\in B^{i-1}}P_{\theta^y\omega}^o 
                     \{\mbox{visit $o$ in $[\tau_i, \tau_{i+1})$}\}
                     1_{\omega\in\Omega_i}\nonumber\\
&\le 
C \varepsilon_0^{-d}\abs{B^{i+2}}^{-1}
         \sum_{x\in B^{i+2}}l_x^{2d} (2d)^{dl_x} E_\omega^x 
         (\mbox{\# visits at $B^{i-1}$ before $\tau_{i+2}$})\nonumber\\
&\stackrel{\text{Lemma }\ref{Ctau}}{\le} 
C \varepsilon_0^{-d}K^{(2-d)i}\sum_{x\in B^{i+2}}l_x^{2d} (2d)^{dl_x}.                              
\end{align}

Since 
\begin{align}\label{CLTe6}
&\sum_{y\in B^{i-1}}P_{\theta^y\omega}^o 
                     \{\mbox{visit $o$ in $[\tau_i, \tau_{i+1})$}\}\nonumber\\
&\le \sum_{y\in B^{i-1}}P_{\theta^y\omega}^o 
                     \{\mbox{visit $o$ in $[\tau_i, \tau_{i+1})$}\}1_{\omega\in\Omega_i}
                     +\abs{B^{i-1}}1_{\omega\notin\Omega_i},                     
\end{align}
taking $P$-expectations on both sides of (\ref{CLTe6}) and using (\ref{CLTe5}) we get
\[
\mathbb{P}^o \{\mbox{visit $o$ in $[\tau_i,\tau_{i+1})$}\}
\le
C \varepsilon_0^{-d}K^{(2-d)i}El_o^{2d}(2d)^{dl_o}+P\{\omega\notin\Omega_i\}.
\]
By (\ref{Clo}), we can take $\varepsilon_0$ to be small enough such that
$El_o^{2d}(2d)^{dl_o}<\infty$ and $\sum_{i=1}^\infty P\{\omega\notin\Omega_i\}<\infty$.
Therefore when $d\ge 3$,
\[\sum_{i=1}^\infty \mathbb{P}^o \{\mbox{visit $o$ in 
$[\tau_i,\tau_{i+1})$}\}<\infty. \quad\quad\quad\quad\quad\quad
\mbox{\qed}\]

\section{Concluding remarks}
 While Bouchaud's trap model (see \cite{Bou,BAC}) provides an example of an
(iid)
environment where local traps can destroy the invariance
principle, it is interesting to note that
a counter-example to Theorem
		\ref{CLT2} in the ergodic setup also can be written.
Namely, let $d\ge 2$, write for $x\in\mathbb{Z}^d$, $z(x)=(x_2,\cdots, x_d)\in\mathbb{Z}^{d-1}$.
Let $\{\varepsilon_z\}_{z\in\mathbb{Z}^{d-1}}$ be i.i.d random variables with support in $(0, 1/2)$ and set
\begin{equation}
\omega(x, e)=\left\{
\begin{array}{rl}
\varepsilon_{z(x)}, & \text{if } e=\pm e_1\\
(1-2\varepsilon_{z(x)})/2(d-1), & \text{else }.
\end{array}
\right.
\end{equation}
It is easy to verify that $\{X_t^n\}_{t\ge 0}$ satisfies the quenched invariance principle, but that
the limiting covariance may degenerate if the tail of $\varepsilon_z$ is heavy.

\chapter{Einstein Relation for Random Walks in Balanced Random Environment}
\label{ER chapter}
%
%
 In this chapter we will give
the proof of the Einstein relation \eqref{Einstein relation} in the context of random walks in a balanced uniformly elliptic iid random environment. As mentioned in Section~\ref{IER}, our proof consists of proving Theorem~\ref{ER1} and Theorem~\ref{ER2}.

We will prove Theorem \ref{ER1} in Section \ref{ERsec1}.
 In Section \ref{ERsecreg}, we will present our new construction of the regeneration times. Furthermore, we will show in Section \ref{ERsecmo} that these regeneration times have good moment properties. Section \ref{ERsecpro} is devoted to the proof of Theorem \ref{ER2}, using the regeneration times and arguments similar to \cite[pages 219-222]{GMP}.

Throughout this chapter, we assume 
\textit{the environment $P$ is iid, balanced, and uniformly 
elliptic with ellipticity constant $\kappa>0$.}  Recall that we have obtained in Section~\ref{SePeEn} an ergodic measure $Q$ for the process $\bar\omega(n)$. By the ergodic theorem, we get
\[
\bm{D}=\big(2E_Q\omega(o,e_i)\delta_{ij}\big)_{1\le i,j\le d},
\]
where $\bm D$ is the covariance matrix defined at the beginning of Section~\ref{IER}.
\section{Proof of Theorem \ref{ER1}}\label{ERsec1}
\begin{lemma}\label{ERl1}
For any $t>0$ and any bounded continuous functional $F$ on $C([0,t],\mathbb{R}^d)$,
\[
\lim_{\lambda\to 0}E_{\omega^\lambda}
F(\lambda X_{s/\lambda^2};0\le s\le t)
=
EF(N_s+D_\ell s;0\le s\le t),
\]
where $(N_s)_{s\ge 0}$ is a $d$-dimensional Brownian motion with covariance matrix $\bm D$.
\end{lemma}
\pf
We first consider the Radon-Nikodym derivative of the measure $P_{\omega^\lambda}$ with respect to
$P_\omega$. Put
\[
G(t, \lambda)=G(t,\lambda;X_\cdot):=\log\prod_{j=1}^{\lceil t\rceil}[1+\lambda\ell\cdot(X_j-X_{j-1})].
\]
Then 
\[
E_{\omega^\lambda} F(X_s: 0\le s\le t)
=
E_\omega F(X_s: 0\le s\le t)e^{G(t,\lambda)}.
\]
In particular, taking $F\equiv 1$, we have
\begin{equation}\label{ERe0}
E_\omega e^{G(t,\lambda)}=1
\end{equation}
for any $\lambda\in (0,1)$ and $t>0$.
Moreover, by the inequality
$a-\frac{a^2}{2}\le \log (1+a)\le a-\frac{a^2}{2}+\frac{a^3}{3}$ for $a>0$, we get
\begin{align}\label{ERe1}
G(t,\lambda)
&=
\sum_{j=1}^{\lceil t \rceil}\log(1+\lambda\ell\cdot(X_j-X_{j-1}))\nonumber\\
&=\sum_{j=1}^{\lceil t \rceil}
\left[\lambda\ell\cdot(X_j-X_{j-1})-\frac{\lambda^2\big(\ell\cdot(X_j-X_{j-1})\big)^2}{2}\right]
+\lambda^2\lceil t\rceil H(\lambda)\nonumber\\
&=\lambda X_{\lceil t \rceil}\cdot\ell
-\frac{\lambda^2}{2}\sum_{j=1}^{\lceil t \rceil}\big(\ell\cdot(X_j-X_{j-1})\big)^2
+\lambda^2\lceil t\rceil H(\lambda),
\end{align}
where the random variable $H(\lambda)=H(\lambda;X_\cdot)$ satisfies $0\le H\le \lambda/3$.
Setting
$h(\omega)=\sum_{i=1}^d\omega(o,e_i)\ell_i^2$,
\[
\left(
\sum_{j=1}^n \big(\ell\cdot(X_j-X_{j-1})\big)^2-2h(\omega_{X_{j-1}})
\right)_{n\ge 0}
\]
is a martingale sequence with bounded increments. Thus $P_\omega$-almost surely,
\[
\lim_{n\to\infty}\frac{1}{n}\sum_{j=1}^n
[\big(\ell\cdot(X_j-X_{j-1})\big)^2-2h(\theta^{X_{j-1}}\omega)]=0.
\]
Further, 
by the ergodic theorem, $P\otimes P_\omega$-almost surely,
\begin{equation}\label{ERe2}
\lim_{\lambda\to 0}
\lambda^2\sum_{j=1}^{\lceil t/\lambda^2 \rceil}
\big(\ell\cdot(X_j-X_{j-1})\big)^2 
=
\lim_{\lambda\to 0}
\lambda^2\sum_{j=1}^{\lceil t/\lambda^2 \rceil} 2h(\theta^{X_{j-1}}\omega)
=
2tE_Q h.
\end{equation}
We deduce from (\ref{ERe1}) and (\ref{ERe2}) that
\[
e^{G(t/\lambda^2,\lambda)}
=\exp[\lambda X_{t/\lambda^2}\cdot\ell-tE_Q h+O_{\lambda,X_\cdot}(1)],
\]
where $O_{\lambda,X_\cdot}(1)$ denotes a quantity that depends on $\lambda$ and
$X_\cdot$, and $O_{\lambda,X_\cdot}(1)\to 0$ $P_\omega$-almost surely as
$\lambda\to 0$. 
By Theorem~\ref{LaThm},
$(\lambda X_{s/\lambda^2})_{s\ge 0}$ converges weakly (under $P_\omega$) to $(N_s)_{s\ge 0}$.
Hence for $P$-almost all $\omega$,
\begin{equation}\label{ERe3}
F(\lambda X_{s/\lambda^2};0\le s\le t)e^{G(t/\lambda^2,\lambda)}
\end{equation}
converges weakly (under $P_\omega$) to
\[
F(N_s:0\le s\le t)\exp(N_t\cdot\ell-tE_Qh).
\]

Next, we will prove that for $P$-almost every $\omega$, this convergence is also in $L^1(P_\omega)$.
It suffices to show that the class
$(e^{G(t/\lambda^2,\lambda)})_{\lambda\in (0,1)}$ is uniformly integrable under $P_\omega$, $P$-a.s.. 
Indeed, for any $\gamma>1$, it follows from (\ref{ERe1}) and the estimate on $H(\lambda)$ that
\begin{align*}
&\gamma G(t/\lambda^2,\lambda)\\
&\le G(t/\lambda^2, \gamma\lambda)
+
\frac{(\gamma^2-\gamma)\lambda^2}{2}
\sum_{j=1}^{\lceil t/\lambda^2 \rceil}
\big(\ell\cdot(X_j-X_{j-1})\big)^2
+\gamma\lambda^2\lceil t/\lambda^2\rceil H(\lambda)\\
&<
G(t/\lambda^2, \gamma\lambda)
+\gamma^2 (t+1).
\end{align*}
Hence for $\gamma>1$ and all $\lambda\in (0,1)$,
\[
E_\omega \exp(\gamma G(t/\lambda^2,\lambda))
\le 
e^{\gamma^2 (t+1)}E_\omega \exp(G(t/\lambda^2, \gamma\lambda))
\stackrel{by (\ref{ERe0})}{=}e^{\gamma^2 (t+1)},
\]
which implies the uniform integrability of $(e^{G(t/\lambda^2,\lambda)})_{\lambda\in (0,1)}$.
 So the $L^1(P_\omega)$-convergence of (\ref{ERe3}) is proved and (for $P$-almost every $\omega$)
we have 
\begin{align*}
&\lim_{\lambda\to 0}
E_{\omega^\lambda}F(\lambda X_{s/\lambda^2};0\le s\le t)\\
&=\lim_{\lambda\to 0}
E_{\omega}F(\lambda X_{s/\lambda^2};0\le s\le t)e^{G(t/\lambda^2,\lambda)}\\
&=E \big[F(N_s:0\le s\le t)\exp(N_t\cdot\ell-t E_Qh)\big].
\end{align*}

The lemma follows by noting that $tE_Qh=E(N_t\cdot\ell)^2/2$ and that, by
Girsanov's formula,
\[
E \big[F(N_s:0\le s\le t)\exp(N_t\cdot\ell-E(N_t\cdot\ell)^2/2)\big]
=
E F(N_s+D_\ell s;0\le s\le t).
\]\qed

\begin{lemma}\label{ERl2}
For any $\lambda\in (0,1), t\ge 1/\lambda^2, p\ge 1$ and any balanced environment $\omega$,
\[
E_{\omega^\lambda}\max_{0\le s\le t}|X_s|^p\le C_{p,d}(\lambda t)^p.
\]
Here we use $C_{p,d}$ to denote constants which depend only on $p$ and the dimension $d$, and
which may differ from line to line.
\end{lemma}
\pf
Since the drift of $\omega^\lambda$ at $X_n, n\in\mathbb{N}$, is
\begin{align*}
E_{\omega^\lambda}(X_{n+1}-X_n|X_n)
&=\sum_{|e|=1}\omega(X_n,e)(1+\lambda e\cdot\ell)e\\
&=\lambda\sum_{i=1}^d 2\omega(X_n,e_i)\ell_i e_i:=\lambda d_\omega(X_{n}),
\end{align*}
we get that
\begin{equation}\label{ERe4}
Y_n:=\lambda\sum_{i=1}^{n}d_\omega(X_{i-1})-X_n
\end{equation}
is a $P_{\omega^\lambda}$-martingale
with bounded increments. By the Azuma-Hoeffding inequality, we get that for any $p\ge 1$,
\[
E_{\omega^\lambda}\max_{1\le i\le n}|Y_i|^p
\le 
C_{p,d}n^{p/2}.
\]
Hence
\[
E_{\omega^\lambda}\max_{1\le i\le n}|X_i|^p
\le 
2^p(E_{\omega^\lambda}\max_{1\le i\le n}|Y_i|^p+\lambda^p n^p)
\le 
C_{p,d}\lambda^p n^p
\]
for any $n\ge 1/\lambda^2$.
The same inequality is true (with different $C_{p,d}$) if we replace $n\in\mathbb{N}$ with
any $t\in\mathbb{R}$ such that $t\ge 1/\lambda^2$.
\qed\\

\textit{Proof of Theorem \ref{ER1}:}
Note that Lemma \ref{ERl1} implies that
$\lambda X_{t/\lambda^2}$ (under the law $P_{\omega^\lambda}$) converges weakly to
$N_t+D_\ell t$ as $\lambda\to 0$. 
When $t\ge 1$, the uniform integrability of $(\lambda X_{t/\lambda^2})_{\lambda\in (0,1)}$
under the corresponding measures $P_{\omega^\lambda}$, as shown in Lemma
\ref{ERl2},
then yields that this convergence is also in $L^1$.\qed

\section{Regenerations}\label{ERsecreg}
\subsection{Auxiliary estimates}
For the rest of this section, we assume that $\ell_1=\ell\cdot e_1>0$.
Let 
\[
\lambda_1:=\big(\lceil(2\lambda\ell_1)^{-1}\rceil\big)^{-1}/2,
\]
so that $0.5/\lambda_1$ is an integer. Note that
\[
\frac{1}{2\lambda\ell_1}\le \frac{1}{2\lambda_1}<\frac{1}{2\lambda\ell_1}+1.
\]

For any $n\in\mathbb{Z}, x\in\mathbb{Z}^d$, call 
\[
\mathcal{H}_n^x=\mathcal{H}_n^x(\lambda,\ell):=
\{y\in\mathbb{Z}^d: (y-x)\cdot e_1=n/\lambda_1\}
\]
\textit{the $n$-th level} (with respect to $x$).
Denote the hitting time of the $n$-th level by
\[
T_n=T_n(X_\cdot):=\inf\{t\ge 0: (X_t-X_0)\cdot e_1=n/\lambda_1\}, n\in\mathbb{Z}.
\]
Also set
\[
T_{\pm 0.5}:=\inf\{t\ge 0: (X_t-X_0)\cdot e_1=\pm 0.5/\lambda_1\}.
\]
Since $\ell_1>0$, the random walk is transient in the $e_1$ direction.
Thus $(T_n)_{n\ge 0}$ are finite $P_{\omega^\lambda}$-almost surely.
\begin{proposition}\label{ERprop5}
For any $n,m\in\mathbb{Z}^+$ and any balanced environment $\omega$,
\[
P_{\omega^\lambda}(T_n<T_{-m})
=
\dfrac{1-q_\lambda^m}{1-q_\lambda^{m+n}},
\]
where $q_\lambda:=(\frac{1-\lambda\ell_1}{1+\lambda\ell_1})^{1/\lambda_1}$.
\end{proposition}
\pf
Observe that the jumps of $(X_n\cdot e_1)_{n\ge 0}$ are lazy random walks on $\mathbb{Z}$, with the ratio of the probabilities of left-jump to right-jump equals $(1-\lambda\ell_1)/(1+\lambda\ell_1)$.
Hence
for $i,j\in\mathbb{Z}^+$,
\[
P_{\omega^\lambda}(\tilde{T}_i<\tilde{T}_{-j})
=\frac{1-(\frac{1-\lambda\ell_1}{1+\lambda\ell_1})^j}
{1-(\frac{1-\lambda\ell_1}{1+\lambda\ell_1})^{i+j}},
\]
where $\tilde{T}_k:=\inf\{n\ge 0: (X_n-X_0)\cdot e_1=k\}, k\in\mathbb{Z}$.
The proposition follows by noting that $T_n=\tilde{T}_{n/\lambda_1}$. \qed

\begin{lemma}\label{ERl3}
For all $\lambda\in(0,1), t>0, m\in\mathbb{N}$ and any balanced environment $\omega$ with ellipticity constant 
$\kappa\in(0,1/(2d))$,
\[
P_{\omega^\lambda}(T_m\ge t/\lambda_1^2)
\le 
2 e^{-t\kappa^2/(2m)}.
\]
\end{lemma}
\pf
First, note that if $Z$ is a real-valued random variable with zero mean and supported on $[-c,c]$, then for $\theta>0$,
$Ee^{\theta Z}\le \exp{(\frac{1}{2}\theta^2c^2)}$. (By Jensen's inequality,
$e^{\theta Z}\le \frac{c-Z}{2c}e^{\theta c}+\frac{c+Z}{2c}e^{-\theta c}$. Taking expectations on both sides gives the inequality.)
Recall the definition of $Y_n$ in (\ref{ERe4}). Since $Y_n\cdot e_1$ is a $P_{\omega^\lambda}$-martingale with increments bounded by $2$, for $\theta>0$,
\begin{align*}
&E_{\omega^\lambda}(e^{\theta Y_{n+1}\cdot e_1}|X_i,i\le n)\\
&=
e^{\theta Y_n\cdot e_1}E_{\omega^\lambda}[e^{\theta(Y_{n+1}-Y_n)\cdot e_1}||X_i,i\le n]
\le 
e^{\theta Y_n\cdot e_1+2\theta^2}.
\end{align*}
Hence
\[
\exp{(\theta Y_n\cdot e_1-2n\theta^2)}
\]
is a $P_{\omega^\lambda}$-supermartingale.
By the optional stopping theorem and ellipticity,
\begin{align*}
1
& \ge 
E_{\omega^\lambda} \exp[\theta Y_{T_m}\cdot e_1-2 T_m\theta^2]\\
& \ge
E_{\omega^\lambda}
\exp[\theta(2\lambda\ell_1\kappa T_m-X_{T_m}\cdot e_1)-2T_m\theta^2].
\end{align*}
Letting $\theta=\kappa\lambda\ell_1/2$
 in the above inequality and noting that 
$X_{T_m}\cdot e_1=m/\lambda_1$, we obtain
\begin{align*}
1 &\ge 
E_{\omega^\lambda}\exp\big((\kappa\lambda\ell_1)^2T_m/2-\kappa\lambda\ell_1 m/(2\lambda_1)\big)\\
&\ge
E_{\omega^\lambda}\exp(\kappa^2\lambda_1^2T_m/2-\kappa m),
\end{align*}
where we used $\lambda_1\le \lambda\ell_1\le 2\lambda_1$ in the second
inequality.
Hence by H\"{o}lder's inequality,
\[
E_{\omega^\lambda}\exp(\kappa^2\lambda_1^2T_m/(2m)-\kappa)\le 1.
\]
Therefore,
\[
P_{\omega^\lambda}(T_m\ge t/\lambda_1^2)
\le  
e^{\kappa-\kappa^2 t/(2m)}
<2 e^{-\kappa^2 t/(2m)}. \qed
\]

\begin{proposition}\label{ERprop2}
There exists a constant $C_0=C_0(\kappa,d)>0$ such that
\[
P_{\omega^\lambda}(\max_{0\le s\le T_1}|X_s|\ge C_0/\lambda_1)<0.5.
\] 
\end{proposition}
\pf
By Lemma \ref{ERl2} and Lemma \ref{ERl3}, for any $m\ge 1$,
\begin{align*}
&P_{\omega^\lambda}(\max_{0\le s\le T_1}|X_s|\ge m/\lambda_1)\\
&\le 
P_{\omega^\lambda}(T_1\ge \sqrt{m}/\lambda_1^2)
+P_{\omega^\lambda}(\max_{0\le s\le\sqrt{m}/\lambda_1^2}|X_s|\ge m/\lambda_1)\\
&\le 
2e^{-\sqrt{m}\kappa^2/2}+C/\sqrt{m},
\end{align*}
which is less than $0.5$ if $m$ is large enough.\qed

\begin{lemma}\label{ERl4}
There exists a constant $c_1\in (0,1]$ such that for any $\lambda\in (0,1)$, $x\in\mathbb{Z}^d$
and balanced environment $\omega$,
\begin{equation}\label{ERe5}
P_{\omega^\lambda}^{x}(X_{T_1}=\cdot)
\ge 
c_1
P_{\omega^\lambda}^{x+0.5e_1/\lambda_1}(X_{T_{0.5}}=\cdot|T_{0.5}<T_{-0.5}).
\end{equation}
\end{lemma}

\pf
For any $x\in\mathbb{Z}^d$, let 
\[
\mathcal{H}_{0.5}^x
:=\{y\in\mathbb{Z}^d: (y-x)\cdot e_1=0.5/\lambda_1\}.
\]
Fix $w\in\mathcal{H}_1^x$. Then the function
\[
f(z):=
P_{\omega^\lambda}^z(X_\cdot \text{ visits $\mathcal{H}_1^x$ for the first time at }w)
\]
satisfies
\[
L_{\omega^\lambda}f(z)=0
\]
for all $z\in\{y: (y-x)\cdot e_1<1/\lambda_1\}$.
By the Harnack inequality for discrete harmonic functions (See Theorem
\ref{ERharnack} in the Appendix. In this case $a=\omega^\lambda,
R=0.5/\lambda_1$ and $b_0\le\lambda$), there exists a constant $C_2$ such that, 
for any $y, z\in \mathcal{H}_{0.5}^x$ with $|z-y|<0.5/\lambda_1$, 
\[
f(z)\ge C_2 f(y).
\]
Hence, for any $z\in \mathcal{H}_{0.5}^x$ such that
$|z-(x+0.5e_1/\lambda_1)|<C_0/\lambda_1$, we have
\begin{equation}\label{ERe27}
f(z)\ge C_2^{2C_0} f(x+0.5e_1/\lambda_1).
\end{equation}
Therefore,
\begin{align*}
P_{\omega^\lambda}^x(X_{T_1}=w)
&\ge 
\sum_{|y-x|<C_0/\lambda}
P_{\omega^\lambda}^x(X_{T_{0.5}}=y)P_{\omega^\lambda}^y(X_{T_{0.5}}=w)\\
&\stackrel{\eqref{ERe27}}{\ge}
C P_{\omega^\lambda}^x(|X_{T_{0.5}}-x|<C_0/\lambda_1)P_{\omega^\lambda}^{x+0.5e_1/\lambda_1}
(X_{T_{0.5}}=w)\\
&\ge 
c_1 P_{\omega^\lambda}^{x+0.5e_1/\lambda_1}
(X_{T_{0.5}}=w|T_{0.5}<T_{-0.5})
\end{align*}
where in the last inequality we used the facts that (by Proposition \ref{ERprop2})
\[
P_{\omega^\lambda}^x(|X_{T_{0.5}}-x|<C_0/\lambda_1)>\frac{1}{2}
\]
and
\[
P_{\omega^\lambda}^{x+0.5e_1/\lambda_1}(T_{0.5}<T_{-0.5})>\frac{1}{2}.\qed
\]

\subsection{Construction of the regeneration times}
Let
\[
\mu_{\omega^\lambda,1}^x(\cdot)=P_{\omega^\lambda}^{x+0.5e_1/\lambda_1}
(X_{T_{0.5}}=\cdot|T_{0.5}<T_{-0.5}).
\]
Recall that $c_1$ is the constant in Lemma \ref{ERl4}. For any $\beta\in (0,c_1)$, we set
\[
\mu_{\omega^\lambda,0}^x(\cdot)=\mu_{\omega^\lambda,0}^{x,\beta}(\cdot)
:=
\big[P_{\omega^\lambda}^x(X_{T_1}=\cdot)-\beta\mu_{\omega^\lambda,1}^x(\cdot)\big]/(1-\beta).
\]
Then by (\ref{ERe5}), both $\mu_{\omega^\lambda,1}^x$ and $\mu_{\omega^\lambda,0}^x$ are probability measures on $\mathcal{H}^x_{0.5}$ and
\[
P_{\omega^\lambda}^x(X_{T_1}=u)=\beta\mu_{\omega^\lambda,1}^x(u)+(1-\beta)\mu_{\omega^\lambda,0}^x(u).
\]

For any $\mathcal{O}\in\sigma(X_1,X_2,\ldots, X_{T_1}), x\in\mathbb{Z}^d$ and $i\in\{0,1\}$, put
\begin{align}\label{ERe6}
\nu_{\omega^\lambda,i}^x(\mathcal{O})
&=
\nu_{\omega^\lambda,i}^{x,\beta}(\mathcal{O})\nonumber\\
&:=
\sum_y 
\big[i\mu_{\omega^\lambda,1}^x(y)+(1-i)\mu_{\omega^\lambda,0}^x(y)\big]
P_{\omega^\lambda}^x(\mathcal{O}|X_{T_1}=y).
\end{align}
Notice that under the environment measure $P$, 
\[
\nu_{\omega^\lambda,1}^x(X_{T_1}\in\cdot)=\mu_{\omega^\lambda,1}^x(\cdot)
\]
 is independent of 
$\sigma(\omega_y:y\cdot e_1\le x\cdot e_1)$.\\

We will now define the regeneration times.

\begin{figure}[h]
\centering
\includegraphics[width=0.9\textwidth]{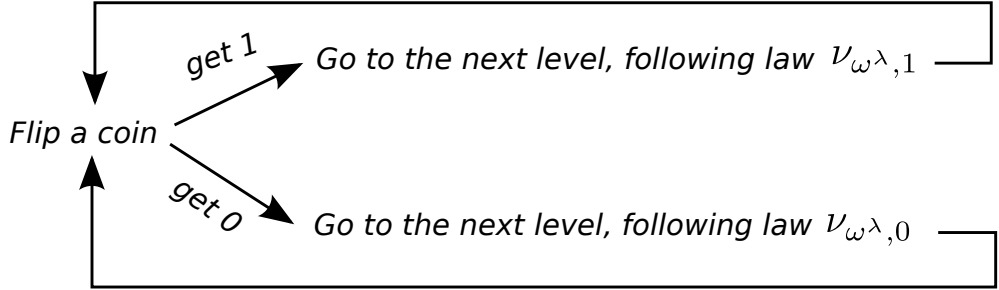}
\caption{The law $\bar P_{\omega^\lambda,\epsilon}$ for the walks.} 
\label{ERfig0}
\end{figure}

We first sample a sequence $(\epsilon_i)_{i=1}^\infty\in\{0,1\}^\mathbb{N}$ of iid Bernoulli random variables according to the law $Q_\beta$ defined by
\[
Q_\beta(\epsilon_i=1)=\beta \text{ and } Q_\beta(\epsilon_i=0)=1-\beta.
\]
Then, fixing $\epsilon:=(\epsilon_i)_{i=1}^\infty$, we will define a new law $P_{\omega^\lambda,\epsilon}$
on the paths as follows (see Figure \ref{ERfig0}). For $x\in\mathbb{Z}^d$, set
\[P_{\omega^\lambda,\epsilon}^x(X_0=x)=1.\]
Assume that the $P_{\omega^\lambda,\epsilon}^x$-law for finite paths of
length$\le n$ is defined. For any path $(x_i)_{i=0}^{n+1}$ with $x_0=x$, define
\begin{align*}
&P_{\omega^\lambda,\epsilon}^{x}
(X_{n+1}=x_{n+1},\ldots, X_{0}=x_0)\\
&:=
P_{\omega,\epsilon}^{x}(X_I=x_I,\ldots, X_0=x_0)
\nu_{\omega^\lambda,\epsilon_J}^{x_I}(X_{n+1-I}=x_{n+1},\ldots, X_1=x_{I+1}),
\end{align*}
where
\[
J=J(x_0,\ldots,x_n):=\max\{j\ge 0: \mathcal{H}_{j}^{x_0}\cap\{x_i, 0\le i\le n\}\neq\emptyset\}
\] 
is the highest level visited by $(x_i)_{i=0}^{n}$ and
\[
I=I(x_0,\ldots,x_n):=\min\{0\le i\le n: x_i\in\mathcal{H}_J^{x_0}\}
\]
is the hitting time to the $J$-th level. By induction, the law $P_{\omega^\lambda,\epsilon}^x$
is well-defined for paths of all lengths. 

Note that a path sampled by $P_{\omega^\lambda,\epsilon}^x$ is not a Markov chain, but
the law of $X_\cdot$ under 
\[
\bar{P}_{\omega^\lambda}^x=\bar{P}_{\omega^\lambda,\beta}^{x}:=Q_\beta\otimes P_{\omega^\lambda,\epsilon}^x
\]
coincides with $P_{\omega^\lambda}^x$. That is, 
\[
\bar{P}_{\omega^\lambda}^x(X_\cdot\in\cdot)
=
P_{\omega^\lambda}^x(X_\cdot\in\cdot).
\] 
Denote by 
$\bar{\mathbb P}_\lambda=\bar{\mathbb P}_{\lambda,\beta}:=P\otimes\bar{P}_{\omega^\lambda,\beta}$ the law of the triple
$(\omega,\epsilon, X_\cdot)$. Expectations with respect to $\bar{P}_{\omega^\lambda}^x$ and $\bar{\mathbb P}_\lambda$ are denoted by 
$\bar{E}_{\omega^\lambda}^x$ and $\bar{\mathbb E}_\lambda (=\bar{\mathbb E}_{\lambda,\beta})$, respectively.

\begin{figure}[h]
\centering
\includegraphics[width=0.9\textwidth]{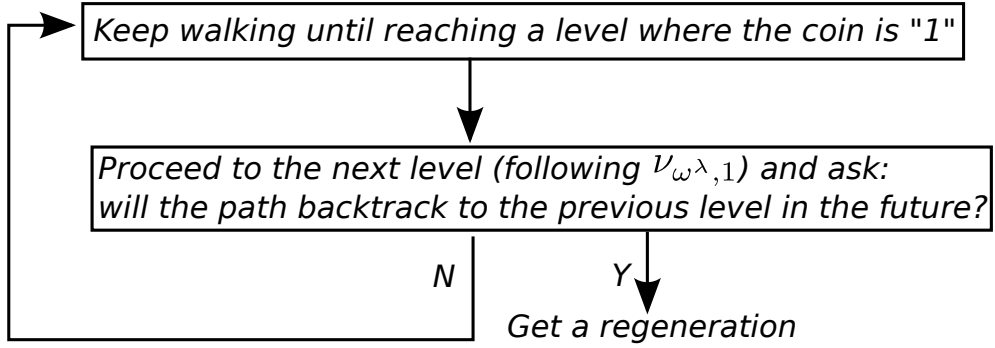}
\caption{The definition of a regeneration time.}
\label{ERfig2}
\end{figure}

Next, for a path $(X_n)_{n\ge 0}$ sampled according to $P_{\omega^\lambda,\epsilon}^o$, we will define the regeneration times. See Figure \ref{ERfig2} for an illustration.

\begin{figure}[h]
\centering
\includegraphics[width=0.8\textwidth]{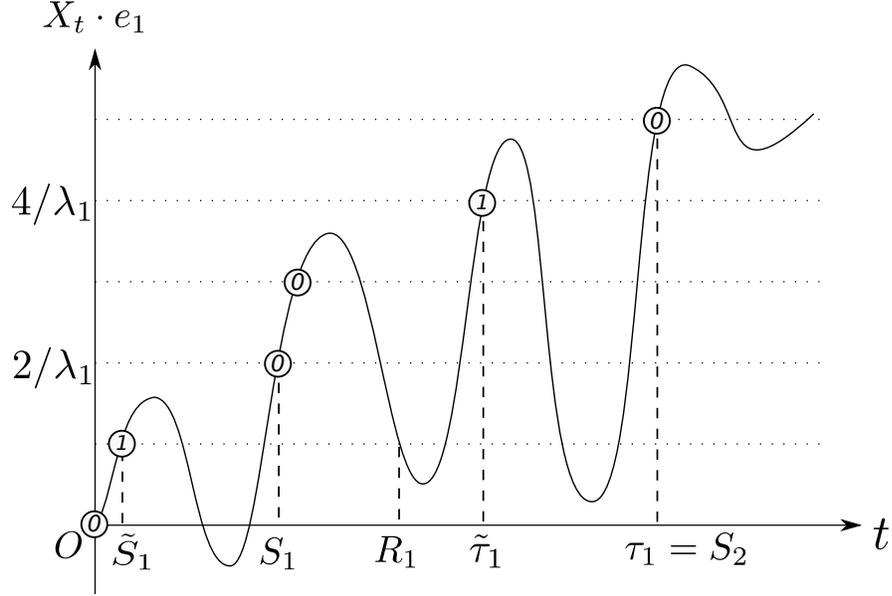}
\caption{In this picture, $K=2, X_{\tau_1}=5/\lambda_1, M_1=4/\lambda_1$.}
\label{ERfig1}
\end{figure}

To be specific, put $S_0=0, M_0=0$,
and define inductively
\begin{align*}
&S_{k+1}=\inf\{T_{n+1}: n/\lambda_1\ge M_k \text{ and }\epsilon_n=1\},\\
&R_{k+1}=S_{k+1}+T_{-1}\circ\theta_{S_{k+1}},\\
&M_{k+1}=X_{S_{k+1}}\cdot e_1+N\circ \theta_{S_{k+1}}, \qquad k\ge 0.
\end{align*}
Here $\theta_n$ denotes the time shift of the path, i.e, $\theta_n X_\cdot=(X_{n+i})_{i=0}^\infty$, 
and
\[
N:=\inf\{n/\lambda_1: n/\lambda_1>(X_i-X_0)\cdot e_1 \text{ for all }i\le T_{-1}\}.
\]
Set
\begin{align*}
&K:=\inf\{k\ge 1: S_k<\infty, R_k=\infty\},\\
&\tau_1:=S_K,\\
&\tau_{k+1}=\tau_k+\tau_1\circ\theta_{\tau_k}.
\end{align*}
We call $(\tau_k)_{k\ge 1}$ the ($\beta$-)\textit{regeneration times}.
Intuitively, under $\bar{P}_{\omega^\lambda}^x$,  whenever the walker visits a new level $\mathcal{H}_i, i\ge 0$, 
he flips a coin $\epsilon_i$. 
If $\epsilon_i=0$ (or $1$), he then walks
following the law $\nu_{\omega^\lambda,0}$ (or $\nu_{\omega^\lambda,1}$) until he hits the $(i+1)$-th level. The regeneration time $\tau_1$ is defined to be the first time of visiting a new level 
$\mathcal{H}_k$ such that the outcome $\epsilon_{k-1}$ of the previous 
coin-tossing is ``$1$" and the path will never backtrack to the level 
$\mathcal{H}_{k-1}$ in the 
future. See Figure \ref{ERfig1}.
\subsection{The renewal property of the regenerations}

The regeneration times possess good renewal properties in the following sense:
\begin{enumerate}
\item 
Since the ratio of the probabilities of left-jump and right-jump of the lazy random walks $(X_n\cdot e_1)_{n\ge 0}$ (in $\mathbb{Z}$) is $(1-\lambda\ell_1)/(1+\lambda\ell_1)$,
the law of $(X_{\tau_n}\cdot e_1)_{n\ge 1}$ does not depend on the
environment $\omega$. (Indeed, if we only observe
the chain $(X_n\cdot e_1)_{n\ge 0}$ at the times when it \textit{moves} and forget about its laziness, we get a random walk on $\mathbb{Z}$ with probabilities $(1-\lambda\ell_1)/2$ and $(1+\lambda\ell_1)/2$ of jumping to the left and to the right, respectively.) Furthermore, under $\bar{P}_{\omega^\lambda}$, the inter-regeneration distances
$(e_1\cdot X_{\tau_1}\circ\theta_{\tau_n})_{n=1}^\infty$ in the direction $e_1$ are iid random variables which are independent of $X_{\tau_1}\cdot e_1$, and
\[
\bar{P}_{\omega^\lambda}(e_1\cdot X_{\tau_1}\circ\theta_{\tau_n}\in\cdot)
=
\bar{P}_{\omega^\lambda}(X_{\tau_1}\cdot e_1\in\cdot|T_{-1}=\infty), n\ge 1.
\]
\item 
For $k\ge 0$, define
\begin{align*}
&\tilde{S}_{k+1}:=\inf\{T_n: n/\lambda\ge M_k \text{ and }\epsilon_n=1\},\\
&\tilde{\tau}_1:=\tilde{S}_K,\\
&\tilde{\tau}_{k+1}:=\tau_k+\tilde{\tau}_1\circ\theta_{\tau_k}.
\end{align*}
Note that for $k\ge 1$,
\begin{align*}
&S_{k}=\tilde S_{k}+T_1\circ\theta_{\tilde S_k},\\
&X_{\tau_k}\cdot e_1
=X_{\tilde{\tau}_k}\cdot e_1+1/\lambda_1.
\end{align*}
Conditioning on $X_{\tilde{\tau}_k}=x$, the law of $X_{\tau_k}$ is
$\mu_{\omega^\lambda,1}^x$, which is independent (under the environment measure $P$) of 
$\sigma(\omega_y:y\cdot e_1\le x\cdot e_1)$. Moreover, after time $\tau_k$,
the path will never visit $\{y:y\cdot e_1\le x\cdot e_1\}$. Thus the movement of the path
after time $\tau_k$ is independent (under $\bar{\mathbb P}_\lambda$) 
of $(X_n)_{n\le \tilde{\tau}_k}$, and therefore, we expect 
\[
(\tilde{\tau}_1\circ\theta_{\tau_k})_{k\ge 1}
\]
to be iid random variables under $\bar{\mathbb P}_\lambda$.
See Proposition \ref{ERprop1} for a rigorous proof.

\item 
Although the inter-regeneration distances $(X_{\tau_1}\circ\theta_{\tau_k})_{k\ge 1}$
and $(\tilde{\tau}_1\circ\theta_{\tau_k})_{k\ge 1}$  are both iid sequences,
the inter-regeneration times $(\tau_1\circ\theta_{\tau_k})_{k\ge 1}$ are not even independent.
However, letting
\[
\Delta_k:=T_1\circ\theta_{X_{\tilde{\tau}_k}}=\tau_k-\tilde{\tau}_k \text{ for }k\ge 1,
\]
we can show that for every $k\ge 1$, $\lambda_1^2 \Delta_k$ is bounded by a constant plus an
exponential random variable. 
So $\Delta_k$ is much less than $\tau_1\circ\theta_{\tau_k}$, which is roughly $C/(\beta\lambda_1^2)$
(as will be shown in Proposition \ref{ERprop3}).
In this sense, the inter-regeneration times $\tau_1\circ\theta_{\tau_k}$ are 
\textit{almost} iid if $\beta$ is sufficiently small.
\end{enumerate}

The rest of this subsection is devoted to the proof that 
$(\tilde{\tau}_1\circ\theta_{\tau_k})_{k\ge 1}$ are iid 
(Proposition \ref{ERprop1}) and
that $\Delta_k$'s are dominated by iid random variables of sizes $1/\lambda^2$ 
(Proposition \ref{ERprop4}).

We introduce the $\sigma$-field
\[
\mathcal{G}_k
:=
\sigma\big(
\tilde{\tau}_k, (X_i)_{i\le\tilde{\tau}_k},(\omega_y)_{y\cdot e_1\le X_{\tilde{\tau}_k}\cdot e_1}
\big).
\]

\begin{lemma}\label{ERl5}
For any appropriate measurable sets $B_1, B_2$
and any event 
\[
B:=\{(X_i)_{i\ge 0}\in B_1, (\omega_y)_{y\cdot e_1>-1/\lambda_1}\in B_2\},
\]
we have, for $k\ge 1$,
\[
\bar{\mathbb P}_\lambda(B\circ\bar{\theta}_{\tau_k}|\mathcal{G}_k)
=
\frac{E_P\big[
\sum_y
\mu_{\omega^\lambda,1}(y)
\bar{P}_{\omega^\lambda}^y(B\cap\{T_{-1}=\infty\})\big]}
{
E_P\big[
\sum_y
\mu_{\omega^\lambda,1}(y)
\bar{P}_{\omega^\lambda}^y(T_{-1}=\infty)\big]
}.
\]
Here  $\bar{\theta}_n$ is the shift defined by
\[
B\circ\bar{\theta}_n
=
\{(X_i)_{i\ge n}\in B_1, (\omega_y)_{(y-X_n)\cdot e_1>-1/\lambda_1}\in B_2\}.
\]
\end{lemma}
\pf
For simplicity, let us consider the case $k=1$. 
We use $\theta^n$ to denote the shift of the $\epsilon$-coins, i.e.,
$\theta^n \epsilon_\cdot=(\epsilon_i)_{i\ge n}$.
For any $A\in\mathcal{G}_1$,
\begin{align*}
&\bar{\mathbb P}_\lambda(B\circ\bar{\theta}_{\tau_1}\cap A)\\
&=
E_{P\otimes Q_\beta}\big[
\sum_{k\ge 1,x}P_{\omega^\lambda,\epsilon}
(A\cap\{\tilde{S}_k<\infty,R_k=\infty, X_{\tilde{S}_k}=x\}\cap B\circ\bar{\theta}_{S_k})
\big]\\
&=
E_{P\otimes Q_\beta}\big[
\sum_{k\ge 1,x,y}P_{\omega^\lambda,\epsilon}
(A\cap\{\tilde{S}_k<\infty,X_{\tilde{S}_k}=x\})
\nu_{\omega^\lambda,1}^x(X_{T_1}=x+y)\\
&\qquad\qquad\qquad\qquad\qquad\times
P_{\omega^\lambda,\theta^{k+1}\epsilon}^{x+y}(B\cap\{T_{-1}=\infty\})\big].
\end{align*}

Note that in the last equality,
\[P_{\omega^\lambda,\epsilon}
(A\cap\{\tilde{S}_k<\infty,X_{\tilde{S}_k}=x\})\]
is 
$\sigma\big((\epsilon_i)_{i\le k},(\omega_z)_{(z-x)\cdot e_1\le 0}\big)$-
measurable,
whereas
\[\nu_{\omega^\lambda,1}^x(X_{T_1}=x+y)P_{\omega^\lambda,\theta^{k+1}\epsilon}^{x+y}(B\cap\{T_{-1}=\infty\})\]
is
$\sigma\big((\epsilon_i)_{i\ge k+1}, (\omega_z)_{(z-x)\cdot e_1>0}\big)$-
measurable for $y\in\mathcal{H}_1^x$.
Hence they are independent under $P\otimes Q_\beta$ and we have
\begin{align}\label{ERe7}
&\bar{\mathbb P}_\lambda(B\circ\bar{\theta}_{\tau_1}\cap A)\\
&=
\sum_{k\ge 1}\bar{\mathbb P}_\lambda
(A\cap\{\tilde{S}_k<\infty\})
E_P\big[
\sum_y
\nu_{\omega^\lambda,1}(X_{T_1}=y)
\bar{P}_{\omega^\lambda}^y(B\cap\{T_{-1}=\infty\})\big].\nonumber
\end{align}
Substituting $B$ with the set of all events, we get
\begin{equation}\label{ERe8}
\bar{\mathbb P}_\lambda(A)=
\sum_{k\ge 1}\bar{\mathbb P}_\lambda
(A\cap\{\tilde{S}_k<\infty\})
E_P\big[
\sum_y
\mu_{\omega^\lambda,1}(y)
\bar{P}_{\omega^\lambda}^y(T_{-1}=\infty)\big].
\end{equation}
(\ref{ERe7}) and (\ref{ERe8}) yield that
\[
\bar{\mathbb P}_\lambda(B\circ\bar{\theta}_{\tau_1}|A)
=
\frac{E_P\big[
\sum_y
\mu_{\omega^\lambda,1}(y)
\bar{P}_{\omega^\lambda}^y(B\cap\{T_{-1}=\infty\})\big]}
{
E_P\big[
\sum_y
\mu_{\omega^\lambda,1}(y)
\bar{P}_{\omega^\lambda}^y(T_{-1}=\infty)\big]
}
.\]
The lemma is proved for the case $k=1$. The general case $k>1$ follows by
induction. (The reasoning for the induction step is the same, although the
notation becomes more cumbersome.)\qed

The following proposition is an immediate consequence of the lemma.
\begin{proposition}\label{ERprop1}
Under $\bar{\mathbb P}_\lambda$, $\tilde{\tau}_1,\tilde{\tau}_1\circ\theta_{\tau_1},\ldots,
\tilde{\tau}_1\circ\theta_{\tau_k},\ldots$ are independent random variables. Furthermore, 
$(\tilde{\tau}_1\circ\theta_{\tau_k})_{k\ge 1}$ are iid with law
\[
\bar{\mathbb P}_\lambda(\tilde{\tau}_1\circ\theta_{\tau_k}\in\cdot)
=
\frac{
E_P\big[
\sum_y
\mu_{\omega^\lambda,1}(y)
\bar{P}_{\omega^\lambda}^y(\tilde{\tau}_1\in\cdot,T_{-1}=\infty)\big]
}
{
E_P\big[
\sum_y
\mu_{\omega^\lambda,1}(y)
\bar{P}_{\omega^\lambda}^y(T_{-1}=\infty)\big]
}.
\]
\end{proposition}
 Note that the inter-regeneration times $(\tau_1\circ\theta_{\tau_k})_{k\ge 1}$ are not independent.
However, the  differences between $\tau_1\circ\theta_{\tau_k}$ and
$\tilde{\tau}_1\circ\theta_{\tau_k}, k\ge 1$ are controlled by iid exponential random variables.

For any $x\in\mathbb{Z}^d, t\ge 0$ ,
\begin{align*}
&\nu_{\omega^\lambda,1}^x(\lambda_1^2T_1\ge t)\\
&\stackrel{\eqref{ERe6}}{=}
\sum_y \mu_{\omega^\lambda,1}^x(y)P_{\omega^\lambda}^x(\lambda_1^2T_1\ge t|X_{T_1}=y)\\
&\stackrel{\eqref{ERe5}}{\le} c_1^{-1}
\sum_y P_{\omega^\lambda}^x(X_{T_1}=y)P_{\omega^\lambda}^x(\lambda_1^2T_1\ge t|X_{T_1}=y)\\
&=
c_1^{-1}P_{\omega^\lambda}^x(\lambda_1^2T_1\ge t)
\stackrel{\text{Lemma }\ref{ERl3}}{\le} 2c_1^{-1}e^{-t\kappa^2/2}.
\end{align*}
Hence for $k\ge 1$,
\[
P_{\omega^\lambda,\epsilon}(\lambda_1^2\Delta_k\ge t|X_i,i\le \tilde{\tau}_k)=
\nu_{\omega^\lambda,1}^{X_{\tilde{\tau}_k}}(\lambda_1^2T_1\ge t)
\le 
2c_1^{-1}e^{-t\kappa^2/2},
\]
which implies that $\lambda_1^2\Delta_k$ is stochastically dominated by an exponential random
variable (with rate $\kappa^2/2$) plus a constant $c_2:=2\kappa^{-2}\log(2/c_1)$.
Thus we conclude:
\begin{proposition}\label{ERprop4}
Enlarging the probability space if necessary, one can couple $(\Delta_k)_{k\ge 1}$ 
with an iid sequence $(\xi_k)_{k\ge 1}$ such that each $\xi_k$ is the sum of $c_2$ and an exponential random
variable with rate $\kappa^2/2$, and that
\[
\lambda_1^2\Delta_k\le \xi_k, \text{ for all }k\ge 1.
\]
Therefore, for any $n\ge 1$,
\begin{equation}\label{ERe9}
\tilde{\tau}_1+\sum_{i=1}^{n-1}\tilde{\tau}_1\circ\theta_{\tau_i}
\le 
\tau_n
\le 
\tilde{\tau}_1+\sum_{i=1}^{n-1}\tilde{\tau}_1\circ\theta_{\tau_i}+\sum_{i=1}^n\xi_i/\lambda_1^2.
\end{equation}
\end{proposition}

\section{Moment estimates}\label{ERsecmo}
Throughout this section, we assume that 
\[\ell\cdot e_1>0.\]
Set $\tau_0=0$. We will show that the typical values of $e_1\cdot X_1\circ\theta_{\tau_k}$
and $\tau\circ\theta_{\tau_k}$, $k\ge 0$ are $C/(\beta\lambda)$ and $C/(\beta\lambda^2)$, respectively.

\begin{theorem}\label{ERregdist}
Let $\omega$ be an elliptic and balanced environment. If $\lambda>0$
 and $\beta>0$ are small enough,
then
\[
\bar{E}_{\omega^\lambda}\exp(\beta\lambda_1 X_{\tau_1}\cdot e_1/2)<12.
\] 
\end{theorem}
\pf
For $0\le k\le K-1$, set
\[
L_{k+1}=\inf\{n\ge \lambda_1 M_k: \epsilon_n=1\}-\lambda_1 M_k+1.
\]
Then $L_1$ is the number of coins tossed to get the first `$1$' and
\[
X_{S_1}\cdot e_1=L_1/\lambda_1.
\]
 Moreover, for $1\le k\le K-1$, let
 \[
N_k=N\circ\theta_{S_k}.
 \]
 Then
\[
(X_{S_{k+1}}-X_{S_k})\cdot e_1=N_{k}+L_{k+1}/\lambda_1, \quad k\ge 1.
\]
So
\begin{equation}\label{ERe10}
X_{\tau_1}\cdot e_1=\sum_{i=1}^K L_i/\lambda_1+\sum_{i=1}^{K-1}N_i.
\end{equation}

First, we will compute the exponential moment of $L_i, i\le K$.
Since $(L_i)_{i\ge 1}$ depends only on the coins $(\epsilon_i)_{i\ge 0}$,
it is easily seen that they are iid geometric random variables with
parameter $\beta$. Hence for $i\ge 1$
(noting that $(1-\beta)e^{\beta/2}< e^{-\beta/2}<1$),
\[
\bar{E}_{\omega^\lambda}
[e^{\beta L_i/2}]
=
\sum_{n=0}^\infty e^{\beta(n+1)/2}(1-\beta)^n\beta
=
\dfrac{\beta e^{\beta/2}}{1-(1-\beta)e^{\beta/2}}.
\]
If $\beta>0$ is small enough, we have 
\begin{equation}\label{ERe11}
\bar{E}_{\omega^\lambda}
[e^{\beta L_i/2}]<3.
\end{equation}

Next, we will compute the exponential moment of $N_i, i\le K-1$.
By Proposition \ref{ERprop5}, putting
\[
p_\lambda:=\bar{P}_{\omega^\lambda}(T_{-1}=\infty)
=1-q_\lambda,
\]
we have
\begin{align*}
&\bar{P}_{\omega^\lambda}(N=(n+1)/\lambda_1)\\
&=
\bar{P}_{\omega^\lambda}(T_n<T_{-1}<T_{n+1})\\
&=
\bar{P}_{\omega^\lambda}(T_n<T_{-1})-
\bar{P}_{\omega^\lambda}(T_{n+1}<T_{-1})\\
&=
\frac{p_\lambda}{1-q_\lambda^{n+1}}-\frac{p_\lambda}{1-q_\lambda^{n+2}}
=\dfrac{q_\lambda^{n+1}p_\lambda^2}{(1-q_\lambda^{n+1})(1-q_\lambda^{n+2})}, 
\quad n\ge 0.
\end{align*}
Observe that conditioning on $K$, 
$(N_i)_{1\le i<K}$ 
are iid under $\bar{P}_{\omega^\lambda}$.
Hence
\begin{align*}
\bar{P}_{\omega^\lambda}(N_i=(n+1)/\lambda_1|K>i)
&=\bar{P}_{\omega^\lambda}(N=(n+1)/\lambda_1|T_{-1}<\infty)\\
&=
\dfrac{q_\lambda^{n}p_\lambda^2}{(1-q_\lambda^{n+1})(1-q_\lambda^{n+2})}\le q_\lambda^n,
\end{align*}
and
\[
\bar{E}_{\omega^\lambda}
[e^{\beta\lambda_1 N_i/2}|K>i]
\le 
\dfrac{e^{\beta/2}}{1-e^{\beta/2}q_\lambda}.
\]
Noting that $\lim_{\lambda\to 0}q_\lambda=e^{-2}$, we can take both $\lambda$ and $\beta$ to be small enough such that
\begin{equation}\label{ERe12}
\bar{E}_{\omega^\lambda}
[e^{\beta\lambda_1 N_i/2}|K>i]
<\frac{1}{4q_\lambda}.
\end{equation}

Finally, note that, under $\bar{P}_{\omega^\lambda}=Q_\beta\otimes
P_{\omega^\lambda,\epsilon}^o$, 
$K$ is a geometric random variable with success parameter $p_\lambda$, and
$(L_i)_{1\le i\le K}$ and $(N_i)_{1\le i\le K}$ are iid sequences when
conditioned on $K$. Therefore, by (\ref{ERe10}), (\ref{ERe11}) and
(\ref{ERe12}),
\[
\bar{E}_{\omega^\lambda}\exp(\beta\lambda_1 X_{\tau_1}\cdot e_1/2)
\le 
\bar{E}_{\omega^\lambda} \frac{3^K}{(4q_\lambda)^{K-1}}
=\sum_{n=0}^\infty \frac{3^{n+1}}{(4q_\lambda)^n}q_\lambda^np_\lambda
<12
\]
if both $\beta,\lambda>0$ are small enough. \qed

\begin{corollary}\label{ERcor1}
For $t\ge 1$ and small enough $\lambda, \beta>0$,
\[
\bar{P}_{\omega^\lambda}(\beta\lambda_1^2\tau_1\ge t)
\le 14\exp(-\kappa^2\sqrt{t}/4).
\]
\end{corollary}
\pf
By Lemma \ref{ERl3} and Theorem \ref{ERregdist},
\begin{align*}
&\bar{P}_{\omega^\lambda}(\beta\lambda_1^2\tau_1\ge t)\\
&\le 
\bar{P}_{\omega^\lambda}(\beta\lambda_1^2T_{\lceil\sqrt{t}/\beta\rceil}\ge t)
+\bar{P}_{\omega^\lambda}(T_{\lceil\sqrt{t}/\beta\rceil}<\tau_1)\\
&\le 
2\exp(-\frac{\kappa^2t/\beta}{2(\sqrt{t}/\beta+1)})
+\bar{P}_{\omega^\lambda}(\lceil\sqrt{t}/\beta\rceil/\lambda_1<X_{\tau_1}\cdot e_1)\\
&\le 
2e^{-\kappa^2\sqrt{t}/4}+12e^{-\sqrt{t}/2}
\le 
14 e^{-\kappa^2\sqrt{t}/4}. \qed
\end{align*}

It follows from Corollary \ref{ERcor1} and Lemma \ref{ERl5} 
(and noting that $P_{\omega^\lambda}(T_{-1}=\infty)=p_\lambda>1/2$) that, 
for $k\ge 1$,
\begin{equation}\label{ERe13}
\bar{\mathbb P}_\lambda
(\beta\lambda_1^2\tau_1\circ\theta_{\tau_k}\ge t)
\le 
28\exp(-\kappa^2\sqrt{t}/4).
\end{equation}
Hence by Theorem \ref{ERregdist}, Corollary \ref{ERcor1} and \eqref{ERe13},
we conclude that, for any $p\ge 1, k\ge 0$, there exists a constant
$C(p)<\infty$ such that
\begin{align}
&\bar{\mathbb E}_\lambda 
(\beta\lambda_1^2\tau_1\circ\theta_{\tau_k})^p
<C(p),\label{ERe14}\\
&\bar{\mathbb E}_\lambda 
(\beta\lambda_1 X_{\tau_1}\circ\theta_{\tau_k})^p
<C(p).\label{ERe15}
\end{align}
Moreover, since $\bar{\mathbb P}_\lambda$-almost surely,
\[
v_\lambda\cdot e_1=\lim_{n\to\infty}\frac{X_{\tau_n}\cdot e_1}{\tau_n},
\]
by \eqref{ERe9}
and the law of large numbers, we
have
\begin{equation}\label{ERe16}
L^{\beta,\lambda}:=\dfrac{\bar{\mathbb E}_\lambda[e_1\cdot X_{\tau_1}\circ\theta_{\tau_1}]}
{\bar{\mathbb E}_\lambda[\tilde{\tau}_1\circ\theta_{\tau_1}]+E\xi_1/\lambda_1^2}
\le 
v_\lambda\cdot e_1
\le 
\dfrac{\bar{\mathbb E}_\lambda[e_1\cdot X_{\tau_1}\circ\theta_{\tau_1}]}
{\bar{\mathbb E}_\lambda[\tilde{\tau}_1\circ\theta_{\tau_1}]}=:R^{\beta,\lambda}.
\end{equation}

\begin{proposition}\label{ERprop3}
When $\lambda, \beta>0$ are small enough,
\begin{equation}\label{ERe17}
\bar{\mathbb E}_\lambda[\tilde{\tau}_1\circ\theta_{\tau_1}]\ge \frac{C}{\beta\lambda_1^2}.
\end{equation}
\end{proposition}
\pf
By the definition of $L_i, i\ge 1$, we get
\begin{equation}\label{ERe18}
\bar{\mathbb E}_\lambda [e_1\cdot X_{\tau_1}\circ\theta_{\tau_1}]
\ge 
\bar{\mathbb E}_\lambda L_1/\lambda_1
\ge 
\frac{1}{\beta\lambda_1}.
\end{equation}
On the other hand, Lemma \ref{ERl2} implies that
\[
|v_\lambda|\le C\lambda \text{ for all }\lambda\in(0,1).
\]
This, together with \eqref{ERe16} and \eqref{ERe18}, yields
\[
\bar{\mathbb E}_\lambda[\tilde{\tau}_1\circ\theta_{\tau_1}]+E\xi_1/\lambda_1^2
\ge \dfrac{C}{\beta\lambda_1^2}.
\]
Recalling (see Proposition \ref{ERprop4}) that $E\xi_1$ is an exponential random variable with rate $\kappa^2/2$, (\ref{ERe17}) then follows by taking $\beta$ sufficiently small.\qed

Note that, by \eqref{ERe16} and \ref{ERe17},
\begin{equation}\label{ERe19}
R^{\beta,\lambda}\le (1+C\beta)L^{\beta,\lambda}\le C\lambda.
\end{equation}

\section{Proof of the Einstein relation}\label{ERsecpro}
\begin{lemma}\label{ERl6}
Assume $\ell\cdot e_1>0$. Then when $\beta>0$ and $\lambda>0$ are small enough,
there exists a constant $C$ such that
\[
\left|
\dfrac{\bar{\mathbb E}_\lambda X_{\tau_n}\cdot e_1}{\lambda\bar{\mathbb E}_\lambda\tau_n}
-\frac{v_\lambda\cdot e_1}{\lambda}
\right|\le C\beta+\frac{C}{n} \quad \text{ for all }n\ge 2.
\]
\end{lemma}
\pf
For $n\ge 2$, since
\[
\dfrac{\bar{\mathbb E}_\lambda X_{\tau_n}\cdot e_1}{\bar{\mathbb E}_\lambda\tau_n}
\ge 
\dfrac{(n-1)\bar{\mathbb E}_\lambda[e_1\cdot X_{\tau_1}\circ\theta_{\tau_1}]}
{\bar{\mathbb E}_\lambda\tau_1+(n-1)\big(\bar{\mathbb E}_\lambda[\tilde{\tau}_1\circ\theta_{\tau_1}]+E\xi_1/\lambda_1^2\big)},
\]
and
\[
\dfrac{\bar{\mathbb E}_\lambda X_{\tau_n}\cdot e_1}{\bar{\mathbb E}_\lambda\tau_n}
\le 
\dfrac{\bar{\mathbb E}_\lambda X_{\tau_1}\cdot e_1+(n-1)\bar{\mathbb E}_\lambda[e_1\cdot X_{\tau_1}\circ\theta_{\tau_1}]}
{(n-1)\bar{\mathbb E}_\lambda[\tilde{\tau}_1\circ\theta_{\tau_1}]},
\]
by the moment bounds \eqref{ERe14}, \eqref{ERe15},  \eqref{ERe17} and \eqref{ERe18}, we have
(for small $\beta$ and $\lambda$)
\[
\frac{L^{\beta,\lambda}/\lambda}{C/(n-1)+1}
\le 
\dfrac{\bar{\mathbb E}_\lambda X_{\tau_n}\cdot e_1}{\lambda\bar{\mathbb E}_\lambda\tau_n}
\le 
\frac{C}{n-1}+\frac{R^{\beta,\lambda}}{\lambda}.
\]
Hence when $\beta>0$ and $\lambda>0$ are small enough and $n\ge 2$, by \eqref{ERe16},
\[
\left|\dfrac{\bar{\mathbb E}_\lambda X_{\tau_n}\cdot e_1}{\lambda\bar{\mathbb E}_\lambda\tau_n}
-\frac{v_\lambda\cdot e_1}{\lambda}\right|
\le 
\frac{C}{n-1}+\frac{R^{\beta,\lambda}}{\lambda}-\frac{L^{\beta,\lambda}/\lambda}{C/(n-1)+1}
\stackrel{(\ref{ERe19})}{\le}
C\beta+\frac{C}{n-1}.
\]
The lemma is proved.\qed

\begin{lemma}\label{ERthm3}
Assume $\ell\cdot e_1>0$.
Let $\alpha_n=\alpha_n(\beta,\lambda):=\bar{\mathbb E}_\lambda\tau_n$. Then when $\beta>0$ and $\lambda>0$ are small enough,
\[
\left|
\dfrac{\bar{\mathbb E}_\lambda X_{\tau_n}\cdot e_1}{\lambda\alpha_n}
-\dfrac{\bar{\mathbb E}_\lambda X_{\alpha_n}\cdot e_1}{\lambda\alpha_n}
\right|\le \frac{C}{n^{1/4}}\quad\text{ for all }n\in \mathbb{N}.
\]
\end{lemma}

Note that, by \eqref{ERe14} and \eqref{ERe17},
\begin{equation}\label{ERe20}
\frac{Cn}{\beta\lambda^2}
\le \alpha_n \le 
\frac{C(1)n}{\beta\lambda^2}.
\end{equation}
\pf
Assume that both $\lambda$ and $\beta$ are sufficiently small.

First, for any $\rho\in(0,1)$,
\begin{align}\label{ERe21}
&\bar{\mathbb E}_\lambda
[|X_{\alpha_n}-X_{\tau_n}|1_{|\tau_n-\alpha_n|\le \rho\alpha_n}]\\
&\le 
\bar{\mathbb E}_\lambda
\big[\max_{(1-\rho)\alpha_n\le s\le (1+\rho)\alpha_n}|X_s-X_{\alpha_n}|\big]\stackrel{\text{Lemma \ref{ERl2}}}{\le}
C\rho\lambda\alpha_n.\nonumber
\end{align}

Second,
\begin{align}\label{ERe22}
&\bar{\mathbb E}_\lambda
[|(X_{\alpha_n}-X_{\tau_n})\cdot e_1|1_{|\tau_n-\alpha_n|> \rho\alpha_n}]\nonumber\\
&\le 
\sqrt{\bar{\mathbb E}_\lambda[|(X_{\alpha_n}-X_{\tau_n})\cdot e_1|^2]
\bar{\mathbb P}_\lambda(|\tau_n-\alpha_n|> \rho\alpha_n)}\nonumber\\
&\stackrel{\text{Lemma \ref{ERl2}, \eqref{ERe15}}}{\le }
Cn(\beta\lambda)^{-1}\sqrt{\bar{\mathbb P}_\lambda(|\tau_n-\alpha_n|> \rho\alpha_n)}.
\end{align}

Furthermore, we can show that
\begin{equation}\label{ERe23}
\bar{\mathbb P}_\lambda(|\tau_n-\alpha_n|> \rho\alpha_n)\le C/(n\rho^2).
\end{equation}
Indeed, put
\[
A_n:=\tilde{\tau}_1+\sum_{i=1}^{n-1}\tilde{\tau}_1\circ\theta_{\tau_i}
\]
and $B_n:=A_n+\sum_{i=1}^n\xi_i/\lambda_1^2$. Then by (\ref{ERe9}), we have
$A_n\le \tau_n\le B_n$. Thus
\[
A_n-\bar{\mathbb E}_\lambda A_n-Cn/\lambda^2
\le 
\tau_n-\alpha_n
\le 
B_n-\bar{\mathbb E}_\lambda B_n+Cn/\lambda^2.
\]
Hence, by \eqref{ERe20} and by taking $\beta>0$ small enough, we get
\begin{align*}
\bar{\mathbb P}_\lambda(\tau_n-\alpha_n> \rho\alpha_n)
&\le 
\bar{\mathbb P}_\lambda(B_n-\bar{\mathbb E}_\lambda B_n\ge \rho\alpha_n/2)\\
&\le 
\frac{\var B_n}{(\rho\alpha_n/2)^2},
\end{align*}
and
\begin{align*}
\bar{\mathbb P}_\lambda(\tau_n-\alpha_n<-\rho\alpha_n)
&\le 
\bar{\mathbb P}_\lambda(A_n-\bar{\mathbb E}_\lambda A_n\le -\rho\alpha_n/2)\\
&\le 
\frac{\var A_n}{(\rho\alpha_n/2)^2}.
\end{align*}
Since (recalling Proposition \ref{ERprop1})
\[
\var A_n=\var \tilde{\tau}_1+(n-1)\var\tilde{\tau}_1\circ\theta_{\tau_1}
\stackrel{(\ref{ERe14})}{\le }
Cn(\beta\lambda^2)^{-2}
\]
and 
\[
\var B_n
\le 
2\big(\var A_n+\var(\sum_{i=1}^n\xi/\lambda_1^2)\big)
=
2\var A_n+Cn/\lambda_1^4
\le 
Cn(\beta\lambda^2)^{-2},
\]
we conclude that
\[
\bar{\mathbb P}_\lambda(|\tau_n-\alpha_n|> \rho\alpha_n)
\le \frac{Cn(\beta\lambda^2)^{-2}}{(\rho\alpha_n/2)^2}
\le C/(n\rho^2).
\]
This completes the proof of \eqref{ERe23}.

Finally, combining (\ref{ERe21}), (\ref{ERe22}) and (\ref{ERe23}), we obtain
\[
\left|
\dfrac{\bar{\mathbb E}_\lambda X_{\tau_n}\cdot e_1}{\lambda\alpha_n}
-\dfrac{\bar{\mathbb E}_\lambda X_{\alpha_n}\cdot e_1}{\lambda\alpha_n}
\right|
\le 
C\rho+\frac{C}{\rho\sqrt{n}}.
\]
The lemma follows by taking $\rho=\frac{1}{n^{1/4}}$.\qed\\

\noindent{\it Proof of Theorem \ref{ER2}:}\\
First, we will show that when $\lambda\in(0,1)$ is small enough, for any $t\ge 1$,
\begin{equation}\label{ERe26}
\left|
\dfrac{\bar{\mathbb E}_\lambda X_{t/\lambda^2}\cdot e_1}{t/\lambda}
-\frac{v_\lambda\cdot e_1}{\lambda}
\right|
\le
\frac{C}{t^{1/5}}.
\end{equation}
Note that if $\ell\cdot e_1=0$,  then
$(X_n\cdot e_1)_{n=0}^\infty$ is a martingale and $\bar{\mathbb E}_\lambda X_n\cdot e_1=v_\lambda\cdot e_1=0$ for all $n$. 
Hence we only consider the non-trivial case $\ell\cdot e_1\neq 0$. Without loss of generality, assume $\ell\cdot e_1>0$. 

By Lemma \ref{ERl2}, the left side of \eqref{ERe26}
is uniformly bounded for all $t\ge 1$ and $\lambda\in (0,1)$. So it suffices to prove \eqref{ERe26} for all sufficiently large $t>0$ and sufficiently small $\lambda>0$.
When $t>0$ is sufficiently large and $\lambda>0$ is small enough, we let 
\begin{equation}\label{ERe24}
\beta=\beta(t)=t^{-1/5}
\end{equation} 
and set $n=n(t,\lambda)$ be the integer that satisfies
\[
	\alpha_n\le \frac{t}{\lambda^2}<\alpha_{n+1}.
\]
By \eqref{ERe20}, the existence of $n(t,\lambda)$ is guaranteed. Moreover,
\begin{equation}\label{ERe25}
n\ge Ct\beta= Ct^{-4/5}.
\end{equation}
Since
\begin{align*}
&\left|
\dfrac{\bar{\mathbb E}_\lambda X_{\alpha_n}\cdot e_1}{\lambda\alpha_n}-\dfrac{\bar{\mathbb E}_\lambda X_{t/\lambda^2}\cdot e_1}{t/\lambda}
\right|\\
&\le \frac{1}{\lambda\alpha_n}\bar{\mathbb E}_\lambda|X_{\alpha_n}-X_{t/\lambda^2}|+
\bar{\mathbb E}_\lambda|X_{t/\lambda}|(\frac{1}{\lambda\alpha_n}-\frac{1}{t})\\
&\le \frac{1}{\lambda\alpha_n}
\bar{\mathbb E}_\lambda
\big[\max_{\alpha_n\le s<\alpha_{n+1}}|X_{\alpha_n}-X_s|\big]+
\bar{\mathbb E}_\lambda[\max_{0\le s<\alpha_{n+1}}|X_s|]
\frac{\lambda \bar{\mathbb E}_\lambda[\tau_1\circ\theta_{\tau_n}]}{(\lambda\alpha_n)^2},
\end{align*}
by Lemma \ref{ERl2}, \eqref{ERe14} and \eqref{ERe20}, we obtain
\begin{equation*}
\left|
\dfrac{\bar{\mathbb E}_\lambda X_{\alpha_n}\cdot e_1}{\lambda\alpha_n}-
\dfrac{\bar{\mathbb E}_\lambda X_{t/\lambda^2}\cdot e_1}{t/\lambda}
\right|
\le 
\frac{C}{n}.
\end{equation*}
Combining Lemma \ref{ERl6}, Lemma \ref{ERthm3} and the above inequality, we conclude that
if $t$ is sufficiently large and $\lambda>0$ is sufficiently small,
then
\[
\left|
\dfrac{\bar{\mathbb E}_\lambda X_{t/\lambda^2}\cdot e_1}{t/\lambda}
-\frac{v_\lambda\cdot e_1}{\lambda}
\right|\le C\beta+\frac{C}{n^{1/4}}
\le
\frac{C}{t^{1/5}}.
\]
Here we used \eqref{ERe24} and \eqref{ERe25} in the last inequality. \eqref{ERe26} is proved.

The same equality for the remaining directions $e_2,e_3,\ldots,e_d$ can be
obtained using the same argument. Our proof of Theorem \ref{ER2} is complete. 
\qed






\appendix
\chapter{}
\label{Appendix}
\section{Discrete harmonic functions}
The purpose of this chapter is to present the proofs of the maximum
principle and the Harnack inequality (Theorem \ref{ERharnack}) for discrete
harmonic functions. The Harnack inequality was used in Section \ref{ERsecreg}
(in the proof of Lemma \ref{ERl4}) to construct the regeneration times.
These inequalities are due to Kuo and Trudinger \cite{KT}.
 
For the purpose of self-containedness, we will give the complete proofs of these
estimates. We follow the arguments in \cite{KT}, adding to it some extra
details.

Recall the definitions of $a$, $L_a$, $b$ and $b_0$ in Section \ref{SeTriid1}.
We consider discrete difference operates $L_a$
such that
\[\sum_y a(x,y)=1, \quad\forall x,\]
and $a(x,y)> 0$ only if $|x-y|=1$, denoted $x\sim y$.
We assume that $L_a$ is uniformly elliptic with constant
$\kappa\in(0,\frac{1}{2d}]$, that is,
\[
a(x,y)\ge \kappa \text{ for any $x, y$ such that $x\sim y$}.
\]
For $r>0, x\in\mathbb{R}^d$, let $B_r(x)=\{z\in\mathbb{Z}^d: |z-x|<r\}$. We also write
$B_r(o)$ as $B_r$. 

\subsection{Maximum principle}
For any bounded set $E\subset\mathbb{Z}^d$, let 
$\partial E=\{y\in E^c:x\sim y \text{ for some }x\in E\}$,
$\bar{E}=E\bigcup\partial E$ and $\diam E=\max\{|x-y|_\infty: x,y\in E\}$.
\begin{theorem}\cite[Theorem 2.1]{KT}\label{AMP}
Let $E\subset\mathbb{Z}^d$ be bounded and $u$ be a function on $\bar{E}$.
For $x\in E$, define 
 \[
I_u(x)=\{s\in\mathbb{R}^d: u(x)-s\cdot x\ge u(z)-s\cdot z, \forall z\in\bar{E}\}. 
 \]
 If 
\[L_a u(x)\ge -g(x)\]
for all $x\in E$ such that 	$I_u(x)=I_u(x,E,a)\neq \emptyset$, then
\[
\max_E u\le 
C\diam(\bar{E})
\big(\sum_{x\in E, I_u(x)\neq \emptyset}|g|^d\big)^{1/d}+\max_{\partial E}u,
\]
where $C$ is a constant determined by $d, \kappa$ and $b_0\diam E$.
\end{theorem}

\pf 
Without loss of generality, we assume $g\ge 0$ and 
\[
\max_E u=u(x_0)>\max_{\partial E}u
\]
for some $x_0\in E$. Otherwise, there is nothing to prove.

For $s\in \mathbb{R}^d$ such that 
\begin{align}\label{A*1}
|s|_\infty
&\le[u(x_0)-\max_{\partial E}u]/(d\diam\bar{E})\nonumber\\
&=:R=R(u,E),
\end{align}
we have
\[
u(x_0)-u(x)\ge s\cdot(x_0-x)
\]
for all $x\in\partial E$, which implies that $\max_{z\in\bar{E}}u(z)-s\cdot z$ is achieved in $E$.
Hence $s\in \bigcup_{x\in E}I_u(x)$ and the cube
\[
Q_R:=\{x:|x|_{\infty}< R\}\subset \bigcup_{x\in E}I_u(x).
\]

For any $p\in\mathbb{R}^d$, set
\[
f(p)=(|p|^{d/d-1}+\mu^{d/d-1})^{1-d},
\]
where $\mu>0$ is a constant to be fixed later.
Since for any $x\in E$, $I_u(x)\subset\mathbb{R}^d$ is bounded and closed, we
can choose
$p_x\in I_u(x)$ so that
\[
|p_x|=\min_{p\in I_u(x)}|p|.
\]
Then
\[
f(p_x)=\max_{p\in I_u(x)}f(p).
\]
Thus
\begin{equation}\label{A*2}
\int_{Q_R}f(s)\ud s\le \int_{\bigcup_{x\in E}I_u(x)}f(s)\ud s
\le
\sum_{x:I_u(x)\neq\emptyset}f(p_x)|I_u(x)|,
\end{equation}
where $|I_u(x)|$ denotes the Lebesgue measure of $I_u(x)$.

Further, 
we will show that, for any $x\in E$ with $I_u(x)\neq\emptyset$,
\begin{equation}\label{A*14}
|I_u(x)|\le (2/\kappa)^{d}[g(x)+b(x)p_x]^d.
\end{equation}
To this end, we fix an $x\in E$ with $I_u(x)\neq\emptyset$
and set
\[
w(z)=u(z)-p_x(z-x),  \quad\forall z\in\bar{E}.
\]
Then $w(x)\ge w(z)$ for all $z\in\bar{E}$ and
\begin{equation}\label{A*3}
I_u(x)=I_w(x)+p_x.
\end{equation}
Since for any $q\in I_w(x)$ and $i=1,\ldots,d$, 
\[
w(x)-w(x\pm e_i)\ge\mp q_i,
\]
we obtain (by ellipticity and by $w(x)\ge w(z)$, $\forall z\in\bar{E}$)
\begin{align*}
0\le \kappa|q|_\infty
&\le\sum_{y}a(x,y)(w(x)-w(y))\\
&=-L_a u+b(x)p_x\\
&\le g(x)+b(x)p_x.
\end{align*}
Hence 
\[
I_w(x)\subset \big[\frac{-g(x)-b(x)p_x}{\kappa},\frac{g(x)+b(x)p_x}{\kappa}\big]^d
\] and
\[
|I_u(x)|\stackrel{(\ref{A*3})}{=}|I_w(x)|\le (2/\kappa)^{d}[g(x)+b(x)p_x]^d.
\]
(\ref{A*14}) is proved.

(\ref{A*14}) and (\ref{A*2}) yield
\[
\int_{Q_R}f(s)\ud s
\le
(\dfrac{2}{\kappa})^d\sum_{x:I_u(x)\neq\emptyset}f(p_x)[g(x)+b(x)p_x]^d.
\]
Since by H\"{o}lder's inequality,
\[
g(x)+|b(x)||p_x|
\le
\big[(\frac{g(x)}{\mu})^d+|b(x)|^d\big]^{1/d} \big[\mu^{d/d-1}+|p_x|^{d/d-1}\big]^{(d-1)/d},
\]
we get
\begin{equation}\label{A*4}
\int_{Q_R}f(s)\ud s
\le 
(\frac{2}{\kappa})^d\sum_{x:I_u(x)\neq\emptyset}
\big[(\frac{g(x)}{\mu})^d+|b(x)|^d\big].
\end{equation}

On the other hand, by H\"{o}lder's inequality,
\[
f(s)=(|s|^{d/d-1}+\mu^{d/d-1})^{1-d}\ge 2^{2-d}(|s|^d+\mu^d)^{-1}.
\]
Thus
\begin{align}\label{A*5}
\int_{Q_R}f(s)\ud s
\ge
\int_{B_R}f(s)\ud s
&\ge 
2^{2-d}\int_{B_R}(|s|^d+\mu^d)^{-1}\ud s\nonumber\\
&=2^{2-d}\frac{\mathcal{O}_d}{d}\log[(\frac{R}{\mu})^d+1],
\end{align}
where $\mathcal{O}_d$ is the area of the unit sphere in $\mathbb{R}^d$.

Finally, combining (\ref{A*4}) and (\ref{A*5}) and putting 
\[
\mu:=[\sum_{x:I_u(x)\neq\emptyset}g(x)^d]^{1/d},\]
we conclude that
\[
\kappa^d 2^{2-2d}\frac{\mathcal{O}_d}{d}\log[(\frac{R}{\mu})^d+1]
\le 1+(b_0\diam\bar{E})^d.
\]
Recalling the definition of $R=R(u,E)$ in (\ref{A*1}), the theorem follows.
\qed

By the same argument as in the proof of Theorem \ref{Cmvi} (Section \ref{SeTriid1}), 
Theorem~\ref{AMP} and Lemma~\ref{Cmvilemma} imply
\begin{theorem}[Mean-value inequality]\label{Amvi}
For any function $u$ on $\bar{B}_R$ such that
\[
L_a u\ge 0, \quad x\in B_R
\]
and any $\sigma\in (0,1)$, $0<p\le d$, we have
\[
\max_{B_{\sigma R}}u\le C\norm{u^+}_{B_R,p},
\]
where $C$ depends on $\sigma, p, \kappa, d$ and $b_0R$.
\end{theorem}
\subsection{Harnack inequality}
\begin{theorem}[Harnack inequality]\cite[Corollary
4.5]{KT}\label{ERharnack}
 Let $u$ be a 
non-negative function on $B_R$, $R>1$. If
\[
L_a u=0
\]
in $B_R$,
then for any $\sigma\in (0,1)$ with $R(1-\sigma)>1$, we have
\[
\max_{B_{\sigma R}}u\le C\min_{B_{\sigma R}}u,
\]
where $C$ is a positive constant depending on $d, \kappa, \sigma$ and $b_0 R$.
\end{theorem}

\begin{lemma}\label{Ahlemma}
Suppose $u$ is a non-negative function on $\bar{B}_R$ that satisfies
\[ L_a u\le 0\]
in $B_R$. Then for any $\sigma\le\tau<1$,
\begin{equation}\label{Aehlemma}
\min_{B_{\tau R}}u\ge C\min_{B_{\sigma R}}u,
\end{equation}
where $C$ depends on $\kappa,d, \sigma,\tau$ and $b_0R$.
\end{lemma}

\pf
Recall the definition of $\eta=\eta_R(x)$ in Lemma \ref{Cmvilemma}. 
We will first show that there exists a constant $\beta=\beta(\sigma, b_0R, \kappa)$ such 
that 
\begin{equation}\label{A*10}
L_a \eta\ge -(2^\beta+\beta^3) R^{-3} \qquad\text{ in }B_R\setminus B_{\sigma R}.
\end{equation}
If $R-1\le |x|< R$, then $\eta(x)\le (2/R)^\beta\le 2^\beta R^{-3}$ for
$\beta\ge 3$.
Hence for $\beta\ge 3$, 
\[L_a\eta\ge -\eta\ge -2^\beta R^{-3}.\]
If $\sigma R\le |x|<R-1$, then $y\in B_R$ for all $y\sim x$.
For $i=1,\ldots, d$, the third derivative $D_i^3\eta$ of $\eta$ with respect to
$x_i$ satisfies
\begin{align*}
 |D_i^3\eta| &=\big|4\beta(\beta-1)x_iR^{-4}\eta^{1-3/\beta}
[3(1-|x|^2 /R^{2})-2(\beta-2)x_i^2/R^2]\big|\\
&\le 4\beta(\beta-1)(2\beta-1)R^{-3},
\end{align*}
and so, by Taylor's expansion, 
\[
\eta(x+e)-\eta(x)\ge \nabla\eta(x)\cdot e+\frac{1}{2}e^{T}D^2\eta(x)e-\frac{8}{6}\beta^3 R^{-3}.
\]
Thus
\begin{align*}
 L_a\eta(x)
 &=\sum_{e}a(x,e)(\eta(x+e)-\eta(x))\\
 &\ge \nabla\eta\cdot b(x)+\frac{1}{2}\sum_e a(x,x+e)e^{T}D^2\eta(x)e-\frac{4}{3}\beta^3 R^{-3}.
\end{align*}
Noting that 
\[
\nabla\eta\cdot b(x)
\stackrel{\eqref{Afifth}, \eta\le 1}{\ge}
-2(b_0R)\beta R^{-2}\eta^{1-2/\beta},
\]
and, for $\sigma R\le |x|<R-1$,
\begin{align*}
&\sum_e a(x,x+e)e^{T}D^2\eta(x)e\\
&=\sum_{i=1}^d(a(x,x+e_i)+a(x,x-e_i))D_{ii}\eta(x)\\
&=2\beta R^{-2}\eta^{1-2/\beta}\sum_{i=1}^d\big(a(x,x+e_i)+a(x,x-e_i)\big)\big(\frac{2(\beta-1)x_i^2}{R^2}-(1-\frac{|x|^2}{R^2})\big)\\
&\ge
2\beta R^{-2}\eta^{1-2/\beta}[4\kappa(\beta-1)\sigma^2-1],
\end{align*}
we have
\[
L_a\eta
\ge
[4\kappa(\beta-1)\sigma^2-1-2b_0R]R^{-2}\eta^{1-2/\beta}
-\frac{4}{3}\beta^3R^{-3}.
\]
Hence \eqref{A*10} also holds for $\sigma R\le |x|<R-1$ if we take
\[
\beta\ge 1+\frac{1+2b_0R}{4\kappa\sigma^2}.
\]
\eqref{A*10} is proved.

Next, let $m_\sigma:=\min_{B_{\sigma R}}u$ and $w:=m_\sigma\eta-u$.
Then
\begin{equation}\label{A*11}
\max_{B_{\tau R}}w\ge (1-\tau^2)^\beta m_\sigma-m_\tau.
\end{equation}
Since $w\le 0$ in $B_{\sigma R}\bigcup B_R^c$ and
\[
L_a w\stackrel{\eqref{A*10}}{\ge}
-(2^\beta+\beta^3)m_\sigma R^{-3} \qquad\text{in }B_R/B_{\sigma R},
\]
we get by the maximum principle that
\begin{equation}\label{A*12}
\max_{B_R}w
\le
C_1 m_\sigma R^{-1},
\end{equation}
where $C_1$ depends on $\kappa,d,\sigma$ and $b_0R$.
By \eqref{A*11} and \eqref{A*12},
\[
[(1-\tau^2)^\beta-\frac{C_1}{R}]m_\sigma\le m_\tau.
\]
Therefore, \eqref{Aehlemma} holds if $R$ satisfies
\[R>\frac{2C_1}{(1-\tau^2)^\beta}.\]
For $R\le\frac{2C_1}{(1-\tau^2)^\beta}$, it follows by iteration (noting $\kappa
u(x)\le u(y)$ for $x\sim y$) that
\[
\kappa^{2C_1(1-\tau^2)^{-\beta}}m_\sigma\le m_\tau.
\]
\eqref{Aehlemma} is proved.\qed\\

For any $z\in\mathbb{Z}^d$ and any
$n=(n_1,\ldots,n_d)\in\mathbb{N}^d$ , we let
\[
N(z,n):=(z+\prod_{i=1}^d[0,n_i-1])\cap\mathbb{Z}^d.
\]
We say that $N(z,n)$ is \textit{nice} if $n$ satisfies $\max_{i,j}|n_i-n_j|\le 1$.
Call $|n|_\infty$ the \textit{length} of the nice rectangle $N(z,n)$. Intuitively, a nice rectangle is ``nearly a cube". 

\begin{proposition}\label{Aprop0}
Let $u$ be a nonnegative function on $\bar{B}_R, R>0$ such that 
\[
L_a u\le 0 \quad\text{ in }B_R.
\]
Suppose $r\in(0, R/7\sqrt{d}]$ and 
$N=N(z,n)\subset Q_r$ is a nice rectangle in $Q_r$.
Then there exists a constant
$\delta=\delta(d,\kappa,b_0R)\in(0,1)$ such that, if
$\Gamma\subset B_R$ satisfies
\[|\Gamma\cap N|\ge \delta|N|,\]
then
\[
\min_{N'}u\ge C\min_{\Gamma}u,
\]
where 
$N'=(z+\prod_{i=1}^d[-n_i,2n_i-1])\cap\mathbb{Z}^d$
and $C$ depends on $\kappa,d, \sigma, \tau$ and $b_0R$.
\end{proposition}

\pf
When $|n|_\infty=1$, $N$ is a singleton, and the proposition follows by
iteration (noting that $u(x)\le\kappa u(y)$ for any $x\sim y$). So we only
consider the case when the length of $N$ is $\ge 2$.

\begin{figure}
\centering
\includegraphics[width=0.6\textwidth]{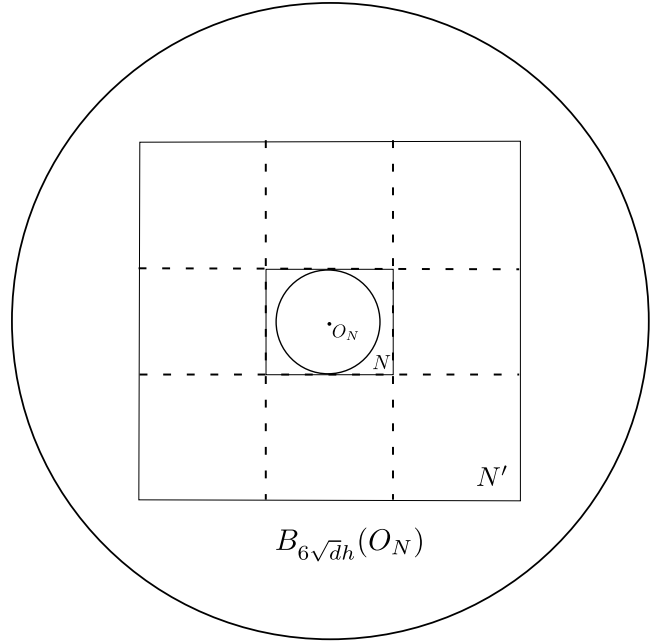}
\caption{$N$ is the rectangle in the center. $B_{h}(O_N)$ is the small circle.}
\label{Afig2}
\end{figure}

Denote the center of $N$ by $O_N=z+(\frac{n_1-1}{2},\cdots,\frac{n_d-1}{2})\in(\frac{1}{2}\mathbb{Z})^d$.
Setting $h=:\min_{i}n_i/2$, we have
\[
B_h(O_N)\subset N\subset B_{2\sqrt{d}h}(O_N).
\]
Since $h\ge\frac{|n|_\infty-1}{2}\ge\frac{|n|_\infty}{4}$, we have
\[
N'\subset B_{3|n|_\infty\sqrt{d}/2}(O_N)
\subset B_{6\sqrt{d}h}(O_N).
\]

Suppose for some $\delta\in(0,1)$,
\[|\Gamma\cap N|\ge \delta|N|.\]
Let $u_\Gamma=:\min_\Gamma u$ and \[v=:u_\Gamma-u,\] then $L_a v\le 0$ and $v^+|_\Gamma=0$.
By Theorem \ref{Amvi}, 
\begin{align*}
\max_{B_{h/2}(O_N)}v
&\le C\frac{1}{|B_{h}|}\sum_{B_h(O_N)}v^+\\
&\le C\frac{|B_h(O_N)\setminus\Gamma|}{|B_h|}\max_{B_h(O_N)}v\\
&\stackrel{|B_h|\ge C|N|}{\le}
C\frac{|N\setminus\Gamma|}{|N|}\max_{B_h(O_N)}v
\le C_2(1-\delta)\max_{B_h(O_N)}v,
\end{align*}
where $C_2$ depends on $\kappa,d$ and $b_0R$.
Taking $\delta=\delta(\kappa, d,b_0R)$ big enough such that $C_2(1-\delta)\le 1/2$, we get
\[\max_{B_{h/2}(O_N)}v
\le \frac{1}{2}\max_{B_h(O_N)}v.\]
Hence
\[
u_\Gamma-\min_{B_{h/2}(O_N)}u\le\frac{1}{2}(u_\Gamma-\min_{B_h(O_N)}u).
\]
Therefore, noting that (since $r\le R/7\sqrt{d}$) $B_{7\sqrt{d}h}(O_N)\subset B_R$,
\[
u_\Gamma\le 2\min_{B_{h/2}(O_N)}u
\stackrel{\text{Lemma \ref{Ahlemma}}}{\le}
C \min_{B_{6\sqrt{d}h}(O_N)}u\le C \min_{N'}u,
\]
with $C$ depending on
$\kappa,d$ and $b_0R$.\qed\\

\begin{lemma}\label{Ahlemma2}
Let $u$ be a nonnegative function on $\bar{B}_R, R>0$ such that 
\[
L_a u\le 0 \quad\text{ in }B_R.
\]
Let $r\in(0, R/7\sqrt{d}]$.
Then for any $\Gamma\subset Q_r$, 
there exists a subset $\Gamma_\delta\supset\Gamma$ of $Q_r$ such that
either $\Gamma_\delta=Q_r$ or $|\Gamma_\delta|>\delta^{-1}|\Gamma|$ holds,
and
\[
\min_{\Gamma_\delta}u\ge\gamma\min_\Gamma u.
\]
Here the constant $\gamma$ depends only on $\kappa, d$ and $b_0R$, and $\delta$ is the same as in Proposition \ref{Aprop0}.
\end{lemma}

\pf
We will construct $\Gamma_\delta$ through a cube decomposition procedure.

Observe that any nice rectangle with length $l\ge 2$ can be decomposed into (at most $2^d$) smaller
disjoint nice rectangles whose lengths are either $\lfloor \frac{l}{2}\rfloor$ or $\lfloor \frac{l}{2}\rfloor+1$. 
With abuse of terminology, we say that such a decomposition is \textit{nice}. Note that a nice decomposition may not
be unique. 

For any $\Gamma\subset Q_r$, set 
\[\mathcal{N}=\mathcal{N}(\Gamma):=\{N:N \text{ is nice and }|\Gamma\cap N|\ge \delta |N|\}.\]	
Now perform cube decompositions to $Q_r$ as follows. Assume that we have an imaginary ``bag".
In the first step, we put $Q_r$ into our ``bag" if $Q_r\in\mathcal{N}$, and decompose $Q_r$ nicely (into at most $2^d$ nice rectangles) if otherwise. 
In the second step, we repeat the same procedure on each of the remaining
rectangles, i.e.,
put a rectangle into our ``bag" if it is in $\mathcal{N}$, and
 decompose a rectangle (with lengths$\ge 2$) nicely if it is not in
$\mathcal{N}$. Repeat this procedure as often as necessary, and stop if there is
nothing to decompose or 
all the remaining rectangles are singletons in $Q_r\setminus\Gamma$.
The process will end within finite number of steps.
Denote the collection of the rectangles in our ``bag" by $\mathcal{N}_0$($\subset\mathcal{N}$).

For $N\in\mathcal{N}_0$ and $N\neq Q_r$, we denote by $N^{-1}$ its \textit{prior}, i.e, $N$ is obtained from a nice decomposition of $N^{-1}$ in the previous step. Set $Q_r^{-1}=Q_r$ and
\[
\Gamma_\delta:=\bigcup_{N\in\mathcal{N}_0}N^{-1}.
\]
Recall the definition of $N'$ in Proposition \ref{Aprop0}. 
For any $N\in\mathcal{N}_0$, since $|\Gamma\cap N|\ge \delta |N|$ and $N^{-1}\subset N'$,
 by the Proposition \ref{Aprop0} we have
\[
\min_{N^{-1}}u
\ge\min_{N'}u
\ge\gamma\min_\Gamma u.
\]
Hence,
\[
\min_{\Gamma_\delta}u\ge\gamma\min_\Gamma u.
\]

Moreover, note that $\Gamma_\delta=Q_r$ when $\mathcal{N}_0=\{Q_r\}$. Otherwise, if 
$\mathcal{N}_0\neq\{Q_r\}$, we have
\[
|\Gamma\cap N^{-1}|<\delta|N^{-1}|\quad\text{for all } N\in\mathcal{N}_0.
\]
Therefore, if $\mathcal{N}_0\neq\{Q_r\}$, \begin{align*}
|\Gamma|=\Abs{\bigcup_{N\in\mathcal{N}_0}(\Gamma\cap N)}
&\le \Abs{\bigcup_{N\in\mathcal{N}_0}(\Gamma\cap N^{-1})}\\
&<\sum_{N^{-1}:N\in\mathcal{N}_0}\delta\abs{N^{-1}}=\delta\abs{\Gamma_\delta}.
\end{align*}
Our proof is complete. \qed\\

\noindent{\it Proof of Theorem \ref{ERharnack}:}\\
We only consider the case when $\sigma<1/7\sqrt{d}$. 

For any $\Gamma\subset Q_{\sigma R}$, if $|Q_{\sigma R}|\le\delta^{-s}|\Gamma|$ for some $s\in\mathbb{N}$,
then we have 
\[
m:=\min_{Q_{\sigma R}}u\ge\gamma^s\min_\Gamma u
\]
by Lemma \ref{Ahlemma2} and iteration.
Hence for $t\ge 0$, putting $\Gamma^t:=\{x\in Q_{\sigma R}: u(x)\ge t\}$, we get
\begin{equation}\label{A*13}
m\ge\gamma^{\lceil\log_\delta(|\Gamma^t|/|Q_{\sigma R}|)\rceil}t
\ge \left(\frac{|\Gamma^t|}{|Q_{\sigma R}|}\right)^{\log_\delta\gamma}\gamma t.
\end{equation}
Note that $q:=\log_\gamma \delta>0$, since $\gamma,\delta\in(0,1)$. 
Therefore, for any $p\in(0,q)$,
\begin{align*}
\frac{1}{|Q_{\sigma R}|}\sum_{Q_{\sigma R}}u^p
&=m^p+\frac{1}{|Q_{\sigma R}|}\sum_{Q_{\sigma R}}\int_m^\infty pt^{p-1}1_{u\ge t}\ud t\\
&=m^p+\int_m^\infty pt^{p-1}\frac{|\Gamma^t|}{|Q_{\sigma R}|}\ud t\\
&\stackrel{(\ref{A*13})}{\le}
m^p+\int_m^\infty pt^{p-1}(\frac{m}{\gamma t})^q\ud t\le C m^p,
\end{align*}
where $C$ depends on $\kappa, d$ and $b_0R$. Combining this and Theorem \ref{Amvi}, the Harnack inequality for $\sigma<1/7\sqrt{d}$ is proved. 

The case $\sigma\ge 1/7\sqrt{d}$ then follows by a chaining argument.
\qed
\end{document}